\documentclass{CUP-JNL-BLG}

\usepackage{latexsym}
\usepackage{graphicx}
\usepackage{multicol,multirow}
\usepackage{amsmath,amssymb,amsfonts}
\usepackage{mathrsfs}
\usepackage{amsthm}
\usepackage{rotating}
\usepackage{appendix}
\usepackage[numbers,super,sort&compress,round]{natbib}
\usepackage{ifpdf}
\usepackage[T1]{fontenc}
\usepackage{times}
\usepackage{newtxmath}
\usepackage{textcomp}%
\usepackage{xcolor}%
\usepackage{hyperref}
\usepackage{lipsum}

\usepackage{import}
\usepackage{amsmath,graphicx}
\usepackage{amssymb}
\usepackage{algorithm}
\usepackage{algorithmic}
\usepackage{mathtools}
\usepackage{caption}
\usepackage{subcaption}
\usepackage{array}
\usepackage{booktabs}
\usepackage{tikz}
\usepackage{rotating}
\usepackage{cleveref}
\usepackage{multicol}
\usepackage{dsfont}
\usepackage{adjustbox}
\usepackage{pgfplots}
\usepackage{changes}
\usepackage{epstopdf}
\usetikzlibrary{arrows,shapes,shadows.blur,spy,fadings,shadings}
\tikzset{verre/.style={draw=SkyBlue,fill=SkyBlue!30}}
\tikzset{punkt/.style={rectangle,rounded corners,draw=black, very thick,text width=2cm,minimum height=2em,text centered}}
\tikzset{punkt2/.style={rectangle,rounded corners,draw=mygreen, very thick,text width=1.5cm,minimum height=2em,text centered}}
\tikzset{double arrow/.style args={#1 colored by #2 and #3}{
         -stealth,line width=#1,#2, % first arrow
          postaction={draw,-stealth,#3,line width=(#1)/3,
          shorten <=(#1)/3,shorten >=2*(#1)/3},}} % second arrow%-- Videos
\usepackage[framemethod=tikz]{mdframed}
\definecolor{myblue}{rgb}{0.122, 0.435, 0.798}
\newmdenv[innerlinewidth=0.5pt, roundcorner=4pt,linecolor=myblue,innerleftmargin=6pt,
innerrightmargin=6pt,innertopmargin=6pt,innerbottommargin=6pt]{mybox}

\def\x{{\mathbf x}}
\def\lspace{$\ell_0$\;}	% l0_
\def\CL{COL0RME}
\def\Phie{\Phi_{{\scriptscriptstyle \mathtt{CEL0}}}} % CEl0 penalty in ND
\DeclareMathOperator*{\argmin}{arg\,min}

\newtheorem{remark}{Remark}

\DeclareMathOperator{\EX}{\mathbb{E}}% expected value
\definecolor{mygreen}{rgb}{0, 0.5, 0.1}

\newcommand{\mdf}[1]{\textcolor{black}{#1}}
\newcommand{\mdfsec}[1]{\textcolor{black}{#1}}
\newcommand{\mdffirst}[1]{\textcolor{black}{#1}}
\newcommand{\mdfsecrev}[1]{\textcolor{black}{#1}}

\def\thebibliography#1{\section{References}\footnotesize\list
 {[\arabic{enumi}]}{\settowidth\labelwidth{[#1]}\leftmargin\labelwidth
 \advance\leftmargin\labelsep
 \usecounter{enumi}}
 \def\newblock{\hskip .11em plus .33em minus .07em}
 \sloppy\clubpenalty4000\widowpenalty4000
 \sfcode`\.=1000\relax}

\jname{Biological Imaging}
\jyear{2022}
\jvol{2}
\jdoi{10.1017/S2633903X22000010}
\begin{document}

\begin{Frontmatter}

\title[Article Title]{COL0RME: Super-resolution microscopy based on sparse blinking/fluctuating fluorophore localization and intensity estimation}

\author*[1]{Vasiliki Stergiopoulou}\email{vasiliki.stergiopoulou@i3s.unice.fr}\orcid{0000-0002-0799-2028}
\author[1]{Luca Calatroni}
\author[2]{José Henrique de Morais Goulart}
\author[3]{Sébastien Schaub}
\author[1]{Laure Blanc-Féraud}

\address[1]{Université Côte d’Azur, CNRS, INRIA, I3S,  \orgaddress{\city{Sophia Antipolis}, \postcode{06900}, \country{France}}}

\address[2]{IRIT, Université de Toulouse, Toulouse INP, \orgaddress{\city{Toulouse}, \postcode{31071}, \country{France}}}

\address[3]{Sorbonne Université, CNRS, LBDV, \orgaddress{\city{Villefranche-sur-Mer}, \postcode{06230}, \country{France}}}

\received{02 August 2021}
\revised{26 January 2021}
\accepted{27 January 2022}

\authormark{Vasiliki Stergiopoulou et al.}

\keywords{Super-Resolution, Fluorescence microscopy, Sparse Optimization, SOFI method}

\abstract{ 
To overcome the physical barriers caused by light diffraction, super-resolution techniques are often applied in fluorescence microscopy. State-of-the-art approaches require specific and often demanding acquisition conditions to achieve adequate levels of both spatial and temporal resolution. Analyzing the stochastic fluctuations of the fluorescent molecules provides a solution to the aforementioned limitations, as sufficiently high spatio-temporal resolution for live-cell imaging can be achieved by using common microscopes and conventional fluorescent dyes. Based on this idea, we present \CL, a method for COvariance-based \lspace super-Resolution Microscopy with intensity Estimation, which achieves good spatio-temporal resolution by solving a sparse optimization problem in the covariance domain and discuss automatic parameter selection strategies. The method is composed of two steps: the former where both the emitters' independence and the sparse distribution of the fluorescent molecules are exploited to provide an accurate localization; the latter where real intensity values are estimated given the computed support. The paper is furnished with several numerical results both on synthetic and real fluorescence microscopy images and several comparisons with state-of-the art approaches are provided. Our results show that COL0RME outperforms competing methods exploiting analogously temporal fluctuations; in particular, it achieves better localization, reduces background artifacts and avoids fine parameter tuning.}

\begin{policy}[Impact Statement]
This research paper describes a super-resolution method improving the spatial resolution of images acquired by common fluorescence microscopes and conventional \mdfsec{blinking/fluctuating} fluorophores. The problem is formulated in terms of a sparse and convex/non-convex optimization problem in the covariance domain for which a well-detailed algorithmic and numerical description are provided. It is addressed to an audience working at the interface between applied mathematics and biological image analysis. The proposed approach is validated on several synthetic datasets and shows promising results also when applied to real data, thus paving the way for new future research directions.
\end{policy}

\end{Frontmatter}

% Introduction
\section{Introduction}
%Many biological entities of very small size are of great importance and their detailed imaging gives valuable information for the biological processes at the cellular and subcellular level. 
In the field of fluorescence (or, more generally, light) microscopy, the main factor characterizing the microscope resolution is the limit imposed by the diffraction of light: structures with size smaller than the diffraction barrier (typically around 250nm in the lateral direction) cannot be well distinguished nor localized. The need to investigate small sub-cellular entities thus led to the implementation of a plethora of super-resolution methods.

A large and powerful family of imaging techniques achieving nanometric resolution are the ones often known as Single Molecule Localization Microscopy (SMLM) techniques, see, e.g. \cite{smlm,SR_fight_club} for a review. Among them, methods such as Photo-Activated Localization Microscopy (PALM) \cite{PALM} and STochastic Optical Reconstruction Microscopy (STORM) \cite{STORM} are designed so as to create a super-resolved image (achieving around $20$nm of resolution) by activating and precisely localizing only a few molecules in each of thousands of acquired frames at a time. For their use, these methods need 
%specific photo-activable fluorophores, 
\mdf{specific photoactivatable, photoswitchable, and binding-activated fluorophores, among others\cite{switchable_fluorophores},}
as well as, a large number (typically thousands) of sparse acquired frames leading to a poor temporal resolution and large exposure times which can significantly damage the sample. A different technique improving spatial resolution is well-known under the name of STimulated Emission Depletion (STED) microscopy \cite{Hell:94}. Similarly to SMLM, STED techniques are based on a time-consuming and possibly harmful acquisition procedure requiring special equipment. In STED microscopy, the size of the point spread function (PSF) is reduced as
\mdf{a depletion beam of light will induce stimulated emission from molecules outside the region of interest and thus switch them off.}
%by using a stimulated emission laser and by switching off some out-of-focus fluorescent molecules. 
Structured Illumination Microscopy (SIM) \cite{sim} methods use patterned illumination to excite the sample; differently from the aforementioned approaches, images here can be recovered with high temporal-resolution via high speed acquisitions that
%do not damage 
\mdf{cause comparatively little damage to}
the sample, but at the cost of a relatively low spatial resolution and, more importantly, the requirement of a specific illumination setup.
%\deleted{Note that the majority of the super-resolution approaches commonly used are grid-based, i.e.~they formalize the super-resolution problem},
\mdf{Note that in this paper we address grid-based super-resolution approaches, i.e.~the ones that formalize the super-resolution problem as the task of retrieving a well-detailed image on a fine grid from coarse measurements.}
More recently, off-the-grid super-resolution approaches have started to be studied in the literature, such as the one of Candès \textit{et al.} \cite{Cands2012SuperResolutionFN}, with applications to SMLM data in Denoyelle \textit{et al.}\cite{Denoyelle_2019}, as well as DAOSTORM\cite{DAOSTORM}, a high-density super-resolution microscopy algorithm. The great advantage of the gridless approaches is that there are no limitations imposed by the size of the discrete grid considered. However, both the theoretical study of the problem and its numerical realization become very hard due to the infinite-dimensional and typically non-convex nature of the optimization.

During the last decade, a new approach taking advantage of the independent stochastic temporal fluctuations/blinking of conventional fluorescent emitters appeared in the literature. A  stack of images is acquired at a high temporal rate, typically $20-100$ images/s, by means of common microscopes (such as widefield, confocal or Total Internal Reflection Fluorecence (TIRF) ones) using standard fluorophores, and then their independent fluctuations/blinking are exploited so as to compute a  super-resolved image. Note that no specific material is needed here, neither for the illumination setup nor for fluorophores. Several methods exploiting the sequence of images have been proposed over the last years. Due to standard acquisition settings, temporal resolution properties are drastically improved.
To start with, Super-resolution Optical Fluctuation Imaging (SOFI) \cite{sofi} is a powerful technique where second and/or higher-order statistical analysis is performed, leading to a significant reduction of the size of the PSF. \mdfsec{An extension of SOFI that combines several cumulant orders and achieves better resolution levels than SOFI is the method bSOFI\cite{bSOFI}.}
%The method bSOFI allows higher order cumulants to be used as the brightness and blinking response is linearized while in SOFI it is not the case.}
%More than SOFI, the method bSOFI allows higher order cumulants to be used due to the extraction of meaningful parameters like the on-time ratio, the brightness and the density of the fluorophores, and as a result achieves better resolution levels.
However, spatial resolution still cannot reach the same levels of PALM/STORM. Almost the same behavior has been noticed in Super-Resolution Radial Fluctuations (SRRF) \cite{srrf} microscopy, where super-resolution is achieved by calculating the degree of local symmetry at each frame. Despite its  easy manipulation and broad applicability, SRRF creates significant reconstruction artifacts which may limit its use in view of accurate analysis. \mdfsec{Other methods which belong to the same category and are worth mentioning are: the method 3B \cite{3B}, which uses Bayesian analysis and takes advantage of the blinking and bleaching events of standard fluorescent molecules, the method Entropy-based Super-resolution Imaging (ESI) \cite{ESI} that computes entropy values pixel-by-pixel, weighted with higher order statistics and the method Spatial COvariance REconstructive (SCORE) \cite{SCORE} that analyzes intensity statistics, similarly to SOFI, but further reduces noise and computational cost by computing only a few components that have a significant contribution to the intensity variances of the pixels.} %that takes into account the covariance between all pairs of pixels and by sorting the eigen modes, efficiently resolves the emitter distribution from a reduced linear subspace of principal components rather than the original image pixel space with higher dimension. = SCORE needs rephrasing
In addition, the approach SPARsity-based super-resolution COrrelation Microscopy (SPARCOM) \cite{SPARCOMold,SPARCOM} exploits, as SOFI, both the lack of correlation between distinct emitters as well as the sparse distribution of the fluorescent molecules via the use of an $\ell_1$ regularization defined on the emitters' covariance matrix. Along the same lines, a deep-learning method exploiting algorithmic unfolding, called Learned SPARCOM (LSPARCOM) \cite{LSPARCOM}, has recently been introduced. Differently from plain SPARCOM, the advantage of LSPARCOM is that neither previous knowledge of the PSF nor any heuristic choice of the regularization parameter for tuning the sparsity level is required. As far as the reconstruction quality is concerned, both SPARCOM and LSPARCOM create some artifacts under challenging imaging conditions, for example when the noise and/or background level are relatively high.
Finally, without using higher order statistics, a constrained tensor modeling approach that estimates a map of local molecule densities and their overall intensities, as well as, a matrix-based \mdf{formulation} that promotes structure sparsity via an $\ell_0$ type regularizer, are available in \cite{goulart}. These approaches can achieve excellent temporal resolution levels, but the spatial resolution is limited.

\paragraph{Contribution}
In this paper, we propose a method for live-cell super-resolution imaging based on the sparse analysis of the stochastic fluctuations of molecule intensities. The proposed approach provides a good level of both temporal and spatial resolution, thus allowing for both precise molecule localization and intensity estimation at the same time, while relaxing the need for special equipment (microscope, fluorescent dyes) typically encountered in state-of-the art super-resolution methods such as, e.g., SMLM. The proposed method is called COL0RME, which stands for COvariance-based super-Resolution Microscopy with intensity Estimation. Similarly to SPARCOM \cite{SPARCOM}, COL0RME enforces signal sparsity in the covariance domain by means of sparsity-promoting terms, either of convex ($\ell_1$, TV) or non-convex ($\ell_0$-based)-type. Differently from SPARCOM, COL0RME allows also for an accurate estimation of the noise variance in the data and is complemented with an automatic selection strategy of the model hyperparameters. Furthermore, and more importantly, COL0RME allows for the estimation of both signal and background intensity, which are relevant pieces of information for biological studies. By exploiting information on the estimated noise statistics, the parameter selection in this step is also made fully automatic, based on the standard discrepancy principle. We remark that an earlier version of  COL0RME has been already introduced by the authors in \cite{ISBI_COL0RME}. Here, we consider an extended formulation combined with automatic parameter selection strategies which allows for the analysis of more challenging data having, e.g., spatially varying background.
%By providing an extensive presentation of numerical results, we show that the method has a broad applicability to data obtained under different imaging conditions. 
The method is validated on simulated and tested on challenging real data. Our results show that COL0RME outperforms competing methods in terms of localization precision, parameter tuning and removal of background artifacts.   
%and to the best of our knowledge, \CL\ is the only super-resolution method exploiting temporal fluctuations which is capable of retrieving this information. 
%Simulated and experimental results show that \CL\ performs well in terms of molecule localization and can retrieve accurate intensity information.

%\txtr{(Not the correct place)} As regards the method COL0RME, as well as the majority of the super-resolution approaches, the reconstructed image lies on a finer grid. Therefore, the resolution of the super-resolved image is limited, due to the chosen discrete grid. To overcome this limitation, grid-less super-resolution approaches can be found in the literature in a theoretical level \cite{Cands2012SuperResolutionFN}, but also as an application addressing SMLM data \cite{Denoyelle_2019}. However, with the grid absence, the computations become very hard due to the "infinite dimensions", making it more preferable to lose some accuracy by being limited by the grid.
%----------------------------------------------------------%

\begin{figure}
    \centering
    \includegraphics[width=0.95\linewidth]{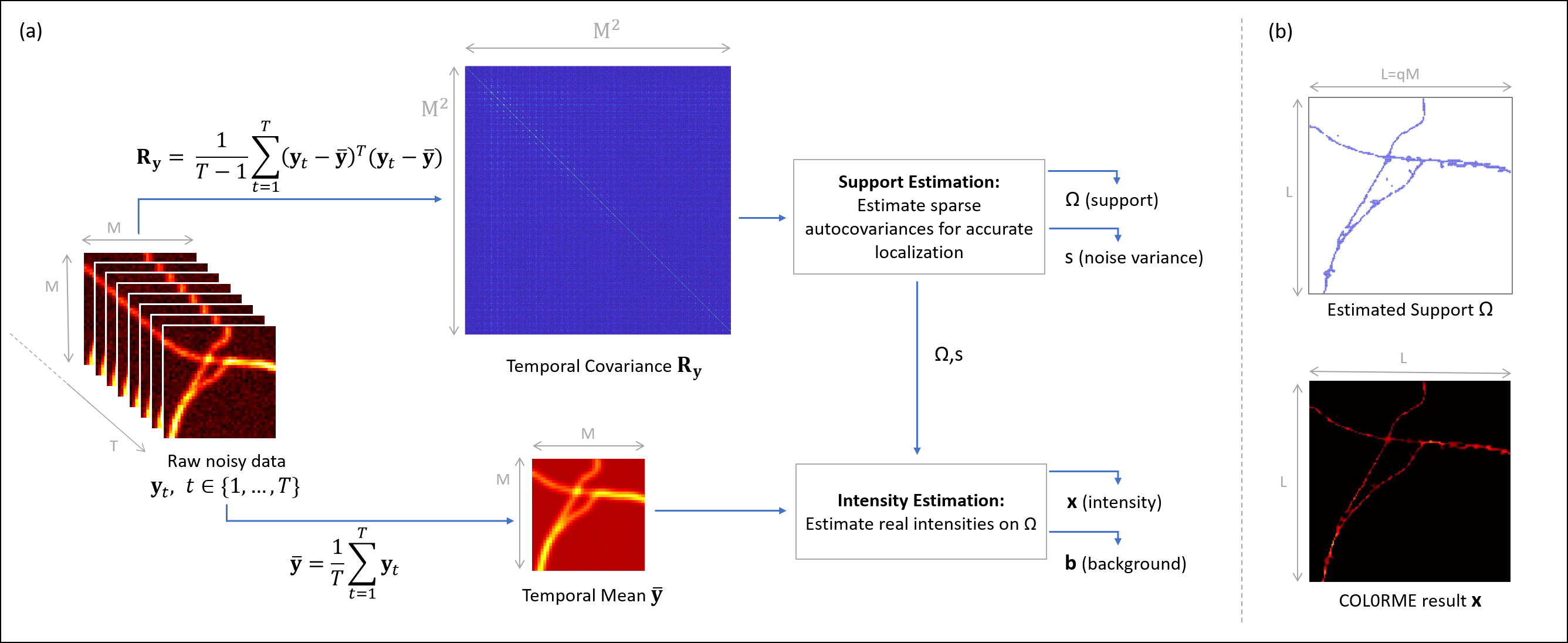}
    \caption{Principles of COL0RME. (a) An overview of the two steps (Support Estimation and Intensity Estimation) by visualizing the inputs/outputs of each, as well as the interaction between them. (b) The two main outputs of COL0RME are: the support $\Omega \subset \mathbb{R}^{L^2}$ containing the locations of the fine-grid pixels with at least one fluorescent molecule, and the intensity $\mathbf{x} \in \mathbb{R}^{L^2}$ whose non-null values are estimated only on $\Omega$}
    %Note that they are both defined on a $q$-times finer grid, with $q=4$ in this example 
    \label{fig:my_label}
\end{figure}

% Mathematical Modeling
\section{Mathematical Modeling}  \label{sec:math_mod}

For real scalars $T, M>1$ and $t\in\left\{1, 2, \ldots,T\right\}$, let $\mathbf{Y}_t \in \mathbb{R}^{M\times M}$ be the blurred, noisy and down-sampled image frame acquired at time $t$. We look for a high-resolution image $\mathbf{X} \in \mathbb{R}^{L\times L}$  being defined as $\mathbf{X} = \frac{1}{T} \sum_{t=1}^{T} \mathbf{X}_t$ with $L = qM$ and defined on a $q$-times finer grid, with $q\in\mathbb{N}$.  Note that in the following applications we typically set $q=4$. The image formation model describing the acquisition process at each  $t$ can be written as:
\begin{equation} \label{eq:model}
    \mathbf{Y}_t = {\cal{M}}_q ({\cal{H}}( \mathbf{X}_t)) + \mathbf{B} +\mathbf{N_t},
\end{equation}
where ${\cal{M}}_q:\mathbb{R}^{L \times L} \rightarrow \mathbb{R}^{M \times M}$ is a down-sampling operator summing every $q$ consecutive pixels in both dimensions, ${\cal{H}}:\mathbb{R}^{L \times L} \rightarrow \mathbb{R}^{L \times L}$ is a convolution operator defined by the PSF of the optical imaging system and $\mathbf{B}\in \mathbb{R}^{M\times M}$ models the background, which collects the contributions of the out-of-focus (and the ambient) fluorescent molecules. Motivated by experimental observations showing that the \mdfsec{blinking/fluctuating} behaviour of the out-of-focus molecules is not visible after convolution with wide de-focused PSFs, we assume that the background is temporally constant ($\mathbf{B}$ does not depend on $t$), while we allow it to smoothly vary in space. Finally, $\mathbf{N}_t\in \mathbb{R}^{M\times M}$ describes the presence of noise modeled here as a matrix of independent and identically distributed (i.i.d.) Gaussian random variables with zero mean and variance $s\in\mathbb{R}^+$ \mdf{taking into account both the underlying electronic noise and the noise bias induced by $\mathbf{B}$ (see Remark \ref{rem:poisson model} for more details on the approximation considered).}
%, as we are limited to the chosen grid, 
We assume that the molecules are located at the center of each pixel and that there is no displacement of the specimen during the imaging period, which is a reasonable assumption whenever short time acquisitions are considered.

\mdf{\begin{remark}  \label{rem:poisson model}
% The real model is 2, but approximations lead
A more appropriate model taking also into account the presence of signal-dependent Poisson noise in the data would be the following:
\begin{equation}\label{eq:poisson_model}
    \mathbf{Y}_t = P\left( {\cal{M}}_q \left({\cal{H}}\left( \mathbf{X}_t\right)\right) + \mathbf{B}\right) +\mdf{\mathbf{E_t}} = P\left( {\cal{M}}_q \left({\cal{H}}\left( \mathbf{X}_t\right)\right) \right) + P\left( \mathbf{B} \right) + \mathbf{E}_t ,\qquad \forall t=1,2,\ldots,T,
\end{equation}
where, for $\mathbf{W}\in\mathbb{R}^{M\times M}$, $P(\mathbf{W})$ represents the realization of a multivariate Poisson variable of parameter $\mathbf{W}$ and $\mathbf{E_t} \in \mathbb{R}^{M \times M}$ models electronic noise with a matrix of i.i.d. Gaussian entries of zero mean and constant variance $\sigma^2 \in \mathbb{R}_+$. Note that the second equality in \eqref{eq:poisson_model} holds due to the independence between ${\cal{M}}_q \left({\cal{H}}\left( \mathbf{X}_t\right)\right)$ and $\mathbf{B}$. Model \eqref{eq:poisson_model} is indeed the one we used for the generation of the simulated data, see Section \ref{sec: Simulated Data}. 
% Here explain why we remove B outside of the parenthesis and rewrite the model with the modifications: s=\sigma^2+b instead of s=\sigma^2
However, to simplify the reconstruction process, 
%we considered in our modeling \deleted{the} \mdf{an} additive Gaussian noise component $\mathbf{N}_t$, with zero mean and variance $s\in\mathbb{R}_+$. 
we simplified \eqref{eq:poisson_model} by assuming that $\mathbf{B}$ has sufficiently large entries, so that $P(\mathbf{B})$ can be approximated as $P(\mathbf{B}) \approx \hat{\mathbf{B}}$ with $\hat{\mathbf{B}}_{i,j}\sim\mathcal{N}(\mathbf{B}_{i,j},\mathbf{B}_{i,j})$, where $(i,j) \in \{1,\dots,M \}^2$, thus considering:
\begin{equation}  \label{eq:poisson_gauss_mod}
\mathbf{Y}_t =  P\left( {\cal{M}}_q \left({\cal{H}}\left( \mathbf{X}_t\right)\right) \right) + \hat{\mathbf{B}} + \mathbf{E}_t,\qquad \forall t=1,2,\ldots,T.
\end{equation}
By now further approximating the variance of $\hat{\mathbf{B}}$ with a constant $b\in\mathbb{R}_+$ to be interpreted as the average of $\mathbf{B}$, we have that by simple manipulations:
\[
\hat{\mathbf{B}} + \mathbf{E}_t = \mathbf{B} + \mathbf{N}_t,
\]
where the independence between $\hat{\mathbf{B}}$ and $\mathbf{E}_t$ has been exploited. We can thus retrieve \eqref{eq:model} from \eqref{eq:poisson_gauss_mod} by neglecting the Poisson noise dependence in $P\left({\cal{M}}_q \left({\cal{H}}\left( \mathbf{X}_t\right)\right)\right)$ and that the variance of every entry of the random term $\mathbf{N}_t$ is $s = \sigma^2 + b$. A more detailed and less approximated modelling taking into account the signal-dependent nature of the noise in the data could represent a very interesting area of future research.
%the random term $\mathbf{N}_t$ has to be thought as the sum $s = \sigma^2 + b$.\\
 \end{remark}}
 
In vectorized form, model (\ref{eq:model}) reads:
\begin{equation} \label{eq:model_vec}
    \mathbf{y}_t = \mathbf{\Psi} \mathbf{x}_t + \mathbf{b} +\mathbf{n}_t,
\end{equation}
where $\mathbf{\Psi} \in \mathbb{R}^{M^2 \times L^2}$ is the matrix representing the composition ${\cal{M}}_q \circ {\cal{H}}$, while $\mathbf{y}_t\in\mathbb{R}^{M^2}$, $\mathbf{x}_t\in\mathbb{R}^{L^2}$, $\mathbf{b}\in\mathbb{R}^{M^2}$ and $\mathbf{n}_t\in\mathbb{R}^{M^2}$  are the column-wise vectorizations of $\mathbf{Y}_t$, $\mathbf{X}_t$, $\mathbf{B}$ and $\mathbf{N}_t$  in \eqref{eq:model}, respectively. 

For all $t$ and given $\mathbf{\Psi}$ and  $\mathbf{y}_t$, the problem can thus be formulated as
\[
\text{find }\quad \mathbf{x}=\frac{1}{T}\sum_{t=1}^T \mathbf{x}_t\in\mathbb{R}^{L^2}, \mathbf{b}\in\mathbb{R}^{M^2}\text{ and }s>0\quad\text{s.t.}\quad \mathbf{x}_t\quad \text{ solves }\eqref{eq:model_vec}.
\]
%a super-resolved image $\mathbf{x}$, defined as $\mathbf{x}=\frac{1}{T}\sum_{t=1}^T \mathbf{x}_t$, from the many $\mathbf{y}_t$ acquisitions and in estimating $\mathbf{b}$ and $s$.
In order to exploit the statistical behavior of the fluorescent emitters, we reformulate the model in the covariance domain. This idea was previously exploited by the SOFI approach \cite{sofi} and was shown to significantly reduce the full-width-at-half-maximum (FWHM) of the PSF. In particular, the use of second-order statistics for a Gaussian PSF corresponds to a reduction factor of the FWHM of $\sqrt2$.

To formulate the model, we consider the frames 
$(\mathbf{y}_t)_{t \in \{1,\dots,T\} }$ as $T$ realizations of a random variable $\mathbf{y}$ with covariance matrix defined by:
\begin{equation}{\label{eq:cov_mat}}
    \mathbf{R_y} = \EX_{\mathbf{y}}\{(\mathbf{y} - \EX_{\mathbf{y}}\{\mathbf{y}\})(\mathbf{y} - \EX_{\mathbf{y}}\{\mathbf{y}\})^\intercal\},
\end{equation}
where $\EX_{\mathbf{y}}\{\cdot\}$ denotes the expected value computed w.r.t. to the unknown law of ${\mathbf{y}}$. We estimate $\mathbf{R_y}$ by computing the empirical covariance matrix, i.e.:
\begin{equation*}{\label{eq:cov_mat_emperical}}
    \mathbf{R_y} \approx
    \frac{1}{T-1}\sum_{t=1}^T (\mathbf{y}_t-\overline{\mathbf{y}})(\mathbf{y}_t-\overline{\mathbf{y}})^\intercal,
\end{equation*}
where $\overline{\mathbf{y}}=\frac{1}{T}\sum_{t=1}^T \mathbf{y}_t$ denotes the empirical temporal mean.
From (\ref{eq:model_vec}) and (\ref{eq:cov_mat}), we thus deduce the relation:
\begin{equation}{\label{eq:cov_model}}
    \mathbf{R_y} = \mathbf{\Psi} \mathbf{R_x} \mathbf{\Psi}^\intercal + \mathbf{R_n},
\end{equation}
where $\mathbf{R_x} \in \mathbb{R}^{L^2 \times L^2}$ and $\mathbf{R_n} \in \mathbb{R}^{M^2 \times M^2}$ are the covariance matrices of $(\mathbf{x}_t)_{t \in \{1,\dots,T\} }$ and $(\mathbf{n}_t)_{t \in \{1,\dots,T\} }$,
%$\left\{ \mathbf{x}_t \right\}$ and $\left\{ \mathbf{n}_t\right\}$,
respectively. As the background is stationary by assumption, the covariance matrix of $\mathbf{b}$ is zero. Recalling now that the emitters are uncorrelated by assumption, we deduce that $\mathbf{R_x}$ is diagonal. We thus set $\mathbf{r_x} := \text{ diag}(\mathbf{R_x})\in\mathbb{R}^{L^2}$. Furthermore, by the i.i.d.~assumption on $\mathbf{n}_t$, we have that $\mathbf{R_n} =  s \mathbf{I_{M^2}}$, where $s \in \mathbb{R}_+$ and $\mathbf{I_{M^2}}$ is the identity matrix in $\mathbb{R}^{M^2 \times M^2}$.
%{Furthermore, regarding the noise component $\mathbf{n}_t$ that has uncorrelated entities, the associated covariance matrix $\mathbf{R_n}$ is diagonal with the main diagonal being equal to the vector $\mathbf{s} \in \mathbb{R}^{M^2}_+$. For each component of $\mathbf{s}$, it stands: $\mathbf{s}_i = \sigma^2 + \mathbf{b}_i, \forall{i} \in {1,...,M^2}$. In order to provide a parameter selection strategy, see sub-section \ref{sec: DP}, we approximate only here the background $\mathbf{b}$ with a constant, interpreted as the mean value of $\mathbf{b}$, as we aim to retrieve a constant noise variance $s$. Fianlly the matrix $\mathbf{R_n}$, can be simplified as follows: $\mathbf{R_n} =  s \mathbf{I_{M^2}}$, where $s \in \mathbb{R}_+$ and $\mathbf{I_{M^2}}$ is the identity matrix in $\mathbb{R}^{M^2 \times M^2}$.}
Note that the model in equation \eqref{eq:cov_model} is similar to the SPARCOM one presented in \cite{SPARCOM}, with the difference that here we consider also noise contributions by including in the model the diagonal covariance matrix $\mathbf{R_n}$. Finally, 
the vectorized form of the model in the covariance domain can thus be written as:
\begin{equation*}  \label{eq:support_model}
    \mathbf{r_y} = (\mathbf{\Psi} \odot \mathbf{\Psi}) \mathbf{r_x} + s \mathbf{v_I},
\end{equation*}
where  $\odot$ denotes the Khatri–Rao (column-wise Kronecker) product, $\mathbf{r_y} \in\mathbb{R}^{M^4}$ is the column-wise vectorization of $\mathbf{R_y}$ and $\mathbf{v_I} = \text{vec}(\mathbf{I_{M^2}})$. 

%--------------------------------------------------------------------%

% COL0RME - Support estimation
\section{COL0RME, step I: support estimation for precise  molecule localization} \label{sec:step1}

Similarly to SPARCOM \cite{SPARCOM}, our approach makes use of the fact that the solution $\mathbf{r_x}$ is sparse, while including further the estimation of $s>0$ for dealing with more challenging scenarios. In order to compare specific regularity \emph{a-priori} constraints on the solution, we make use of different regularization terms, whose importance is controlled by a regularization hyperparameter $\lambda>0$. By further introducing some non-negativity constraints for both variables $\mathbf{r_x}$ and $s$, we thus aim to solve:
\begin{equation}  
    \label{eq:support_mini}
    \argmin\limits_{\mathbf{r_x} \geq 0,~ s \geq 0}~ {\cal{F}}(\mathbf{r_x},s) + {\cal{R}}(\mathbf{r_x};\lambda),
\end{equation}
where the data fidelity term is defined by:
\begin{equation}  
{\cal{F}}(\mathbf{r_x},s) = \frac12 \| \mathbf{r_y} -(\mathbf\Psi \odot \mathbf\Psi) \mathbf{r_x} - s \mathbf{v_I} \|_2^2,
\end{equation}
and $\cal{R}(\cdot;\lambda)$ is a sparsity-promoting penalty. Ideally, one would like to make use of the \lspace norm  to enforce sparsity. However, as it is well-known, solving the resulting non-continuous, non-convex and combinatorial  minimization problem is an NP-hard problem. A way to circumvent this difficulty consists in using the continuous exact relaxation of the \lspace norm (CEL0) proposed by Soubies \textit{et al.} in \cite{CELO}. The CEL0 regularization is continuous, non-convex and preserves the global minima of the original $\ell_2-\ell_0$ problem while removing some local ones. It is defined as follows:
\begin{equation}  \label{eq:CEL0}
     {\cal{R}}(\mathbf{r_x};\lambda) =\Phie{(\mathbf{r_x};\lambda)} =  \sum\limits_{i=1}^{L^2} \lambda - \frac{\|\mathbf{a}_i\|^2}{2}\left( |(\mathbf{r_x})_i| - \frac{\sqrt{2\lambda}}{\|\mathbf{a}_i\|} \right)^2 \mathds{1} _{\{|({\mathbf{r_x}})_i| \leq \frac{\sqrt{2 \lambda}}{\|\mathbf{a}_i\|}\}},
\end{equation}
where $\mathbf{a}_i = (\mathbf\Psi \odot \mathbf\Psi)_i$ denotes the $i$-th column of the operator $\mathbf{A}:=\mathbf\Psi \odot \mathbf\Psi$. 
%The functional including the CEL0 penalty is continuous and non-convex, but has the same minimizers as the corresponding $\ell_2-\ell_0$ problem.\\ 

A different, convex way of favoring sparsity consists in taking as regularizer the $\ell_1$ norm, that is:
\begin{align}  
    \label{eq:L1}
    {\cal{R}}(\mathbf{r_x};\lambda) = \lambda \|\mathbf{r_x}\|_1.
\end{align}
Besides convexity and as it is well-known, the key difference between using the $\ell_0$ and the $\ell_1$-norm is that the $\ell_0$ provides a correct interpretation of sparsity by counting only the number of the non-zero coefficients, while the $\ell_1$ depends also on the magnitude of the coefficients. 
However, its use as a sparsity-promoting regularizer is nowadays well-established (see, e.g., \cite{rwl1}) and also used effectively in other microscopy applications, such as SPARCOM \cite{SPARCOM}.
%Nevertheless, as regards the images we deal with and the fact that we want to avoid very sparse and dotted reconstructions, the $\ell_1$-norm is an advantageous choice. 

\begin{figure}[H]
\centering
    \includegraphics[width=0.17\textwidth]{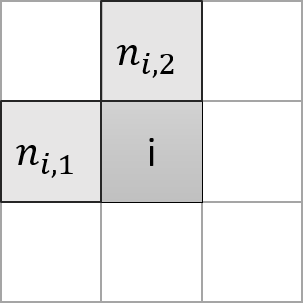}
    \caption{The one-sided nearest horizontal and vertical neighbors of the pixel $i$ used to compute the gradient discretization in \eqref{eq:TV}}
    \label{fig:neighbours} 
\end{figure}

Finally, in order to model situations where piece-wise constant structures are considered, we consider a different regularization term favoring gradient-sparsity by using the Total Variation (TV) regularization defined  in a discrete setting as follows:
\begin{align}    \label{eq:TV}
    {\cal{R}}(\mathbf{r_x};\lambda) = \lambda TV(\mathbf{r_x}) = \lambda \sum\limits_{i=1}^{L^2} \left(|({\mathbf{r_x}})_i - ({\mathbf{r_x}})_{n_{i,1}}|^2+|({\mathbf{r_x}})_i-({\mathbf{r_x}})_{n_{i,2}}|^2 \right)^{\frac12},
\end{align}
where $(n_{i,1},n_{i,2}) \in \{1,\dots,L^2 \}^2$ indicate the locations of the horizontal and vertical nearest neighbor pixels of pixel $i$, as shown in Figure \ref{fig:neighbours}. For the computation of the TV penalty, Neumann boundary conditions have been used.

To solve \eqref{eq:support_mini} we use the Alternate Minimization algorithm  between $s$ and $\mathbf{r_x}$ \cite{attouch}, see the pseudo-code reported in Algorithm \ref{Algorithm:AMA_support}. Note that, at each $k\geq 1$, the update for the variable $s$ can be efficiently computed through the following explicit expression:
\begin{equation*}
s^{k+1} = \frac{1}{M^2} \mathbf{v_I}^\intercal ( \mathbf{r_y} -(\mathbf\Psi \odot \mathbf\Psi) \mathbf{r_x}^k). 
\end{equation*}

%\begin{equation*}
%s^{k+1} = \frac{1}{M^2} \sum_{i=1}^{M^2} ( \mathbf{R_y}(i,i) - \{(\mathbf\Psi \odot \mathbf\Psi) \mathbf{r_x}^k\}(i,i)) \text{ - simpler than that}
%\end{equation*}

Concerning the update of $\mathbf{r_x}$, different algorithms were used depending on the choice of the regularization term in \eqref{eq:CEL0}, \eqref{eq:L1} and \eqref{eq:TV}. For the CEL0 penalty \eqref{eq:CEL0} we used the iteratively reweighted $\ell_1$ algorithm (IRL1) \cite{WRL1}, following Gazagnes et al. \cite{Gazagnes} with Fast Iterative Shrinkage-Thresholding Algorithm (FISTA) \cite{FISTA} as inner solver. 
%With more details,...
If the $\ell_1$ norm \eqref{eq:L1} is chosen, FISTA is used.
%explain a bit
Finally, when the TV penalty \eqref{eq:TV} is employed, the Primal-Dual Splitting Method in \cite{Primal-Dual} was considered. %\txtr{(Should I say more for these 3 algorithms?)}
% maybe more details

\begin{algorithm}[H]
\caption{COL0RME, Step I: Support Estimation}
\label{Algorithm:AMA_support}
\begin{algorithmic}
\REQUIRE $\mathbf{r_y}\in\mathbb{R}^{M^4}, \mathbf{r_x}^0\in\mathbb{R}^{L^2}, \lambda>0$
\REPEAT 
\STATE $s^{k+1} = \argmin\limits_{s \in \mathbb{R}_+} {\cal{F}}(\mathbf{r_x}^{k},s)$
\STATE $\mathbf{r_x}^{k+1} = \argmin\limits_{\mathbf{r_x} \in \mathbb{R}_+^{L^2}}{\cal{F}}(\mathbf{r_x},s^{k+1})+  {\cal{R}}(\mathbf{r_x};\lambda)$
\UNTIL convergence
\RETURN $\Omega_{\mathbf{x}}, s$
\end{algorithmic}
\end{algorithm}

Following the description provided by Attouch \textit{et al.} in \cite{attouch},  convergence of Algorithm \ref{Algorithm:AMA_support} can be guaranteed only if an additional quadratic term is introduced in the objective function of the second minimization sub-problem. Nonetheless, empirical convergence was observed also without such additional terms. % A formal proof of convergence is left for future work.

\begin{figure}[H]
\centering
\setlength\tabcolsep{1.5pt}
\begin{tabular}{ccccc}
    & $\bar{\mathbf{y}}$ + GT & \hspace{-0.7cm}CEL0 result & \hspace{-0.7cm}$\ell_1$ result & \hspace{-0.7cm}TV result \\
(a) & \adjustbox{valign=m,vspace=1pt}{\includegraphics[width=0.28\textwidth]{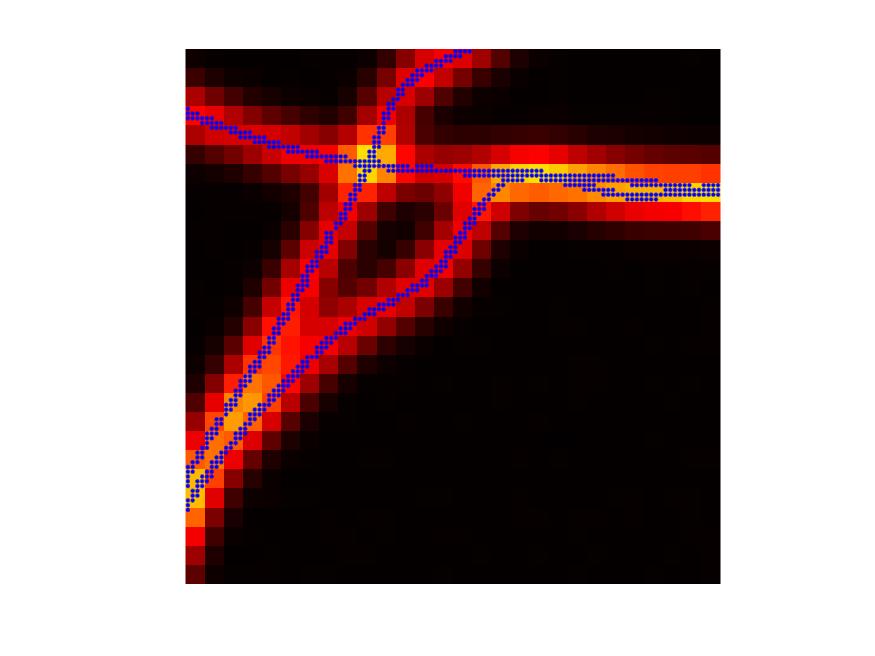}} & \hspace{-0.7cm}\adjustbox{valign=m,vspace=1pt}{\includegraphics[width=.28\textwidth]{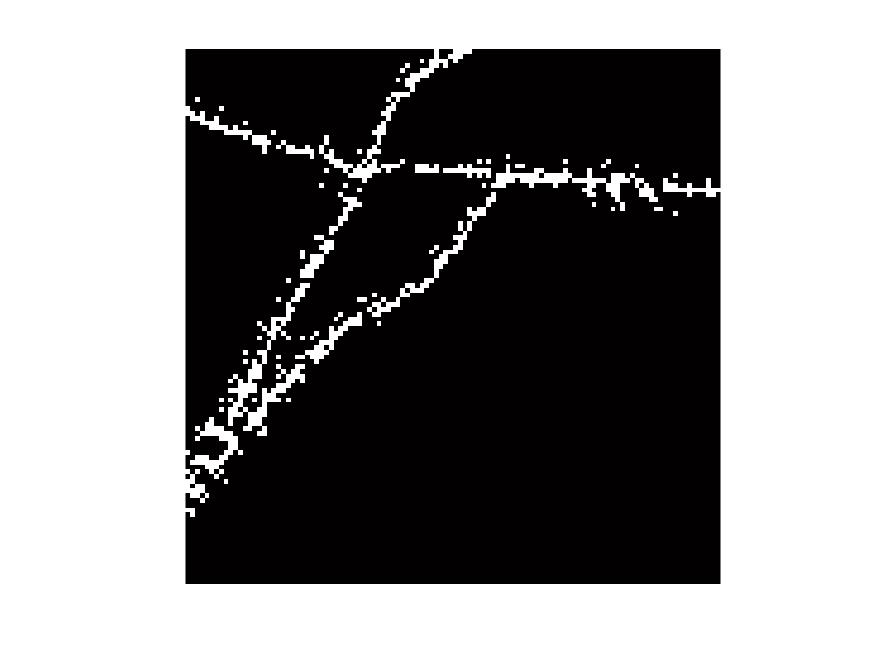}} & \hspace{-0.7cm}\adjustbox{valign=m,vspace=1pt}{\includegraphics[width=.28\textwidth]{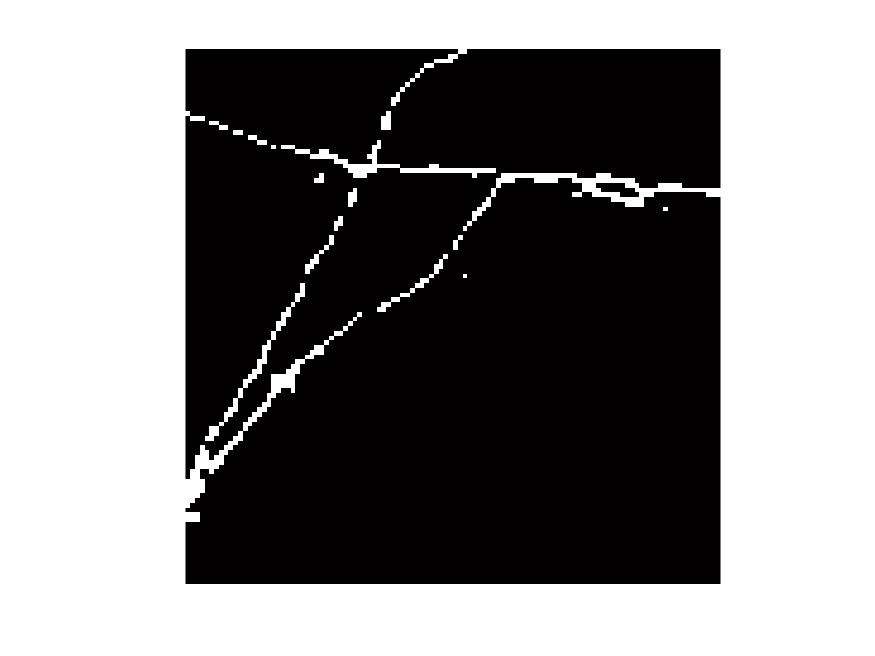}} & \hspace{-0.7cm}\adjustbox{valign=m,vspace=1pt}{\includegraphics[width=.28\textwidth]{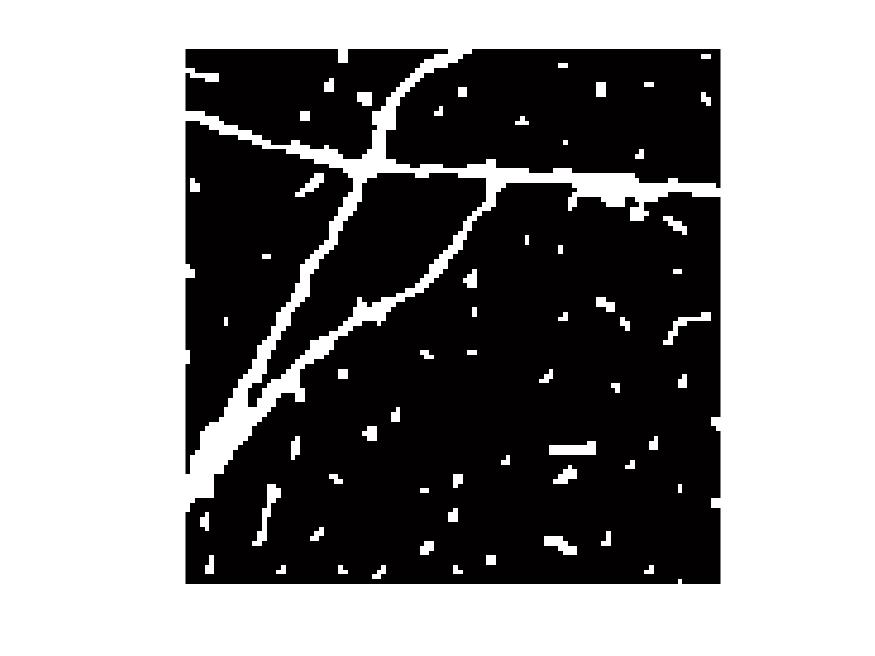}}\\
(b) & \adjustbox{valign=m,vspace=1pt}{\includegraphics[width=0.28\textwidth]{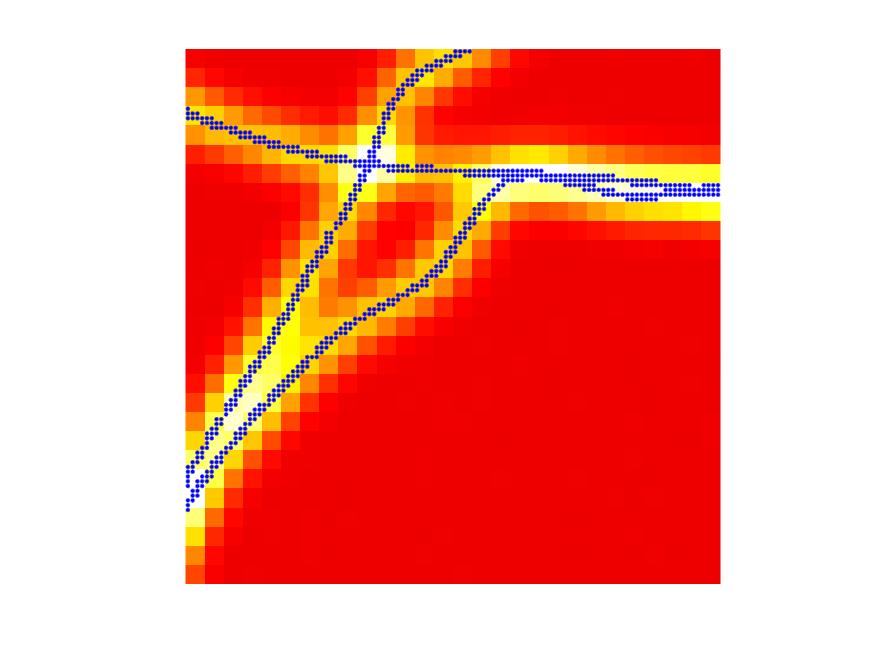}} & \hspace{-0.7cm}\adjustbox{valign=m,vspace=1pt}{\includegraphics[width=.28\textwidth]{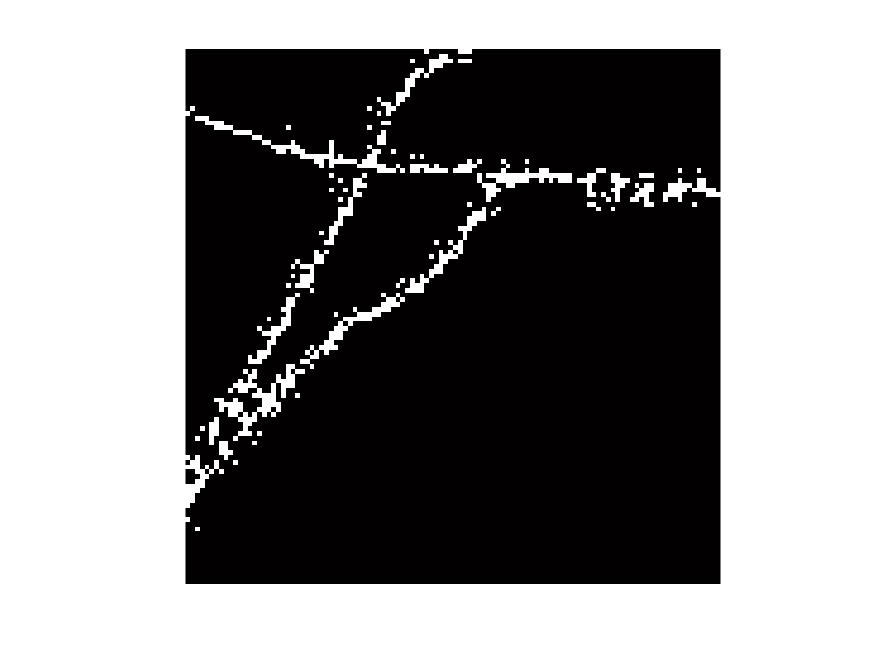}} & \hspace{-0.7cm}\adjustbox{valign=m,vspace=1pt}{\includegraphics[width=.28\textwidth]{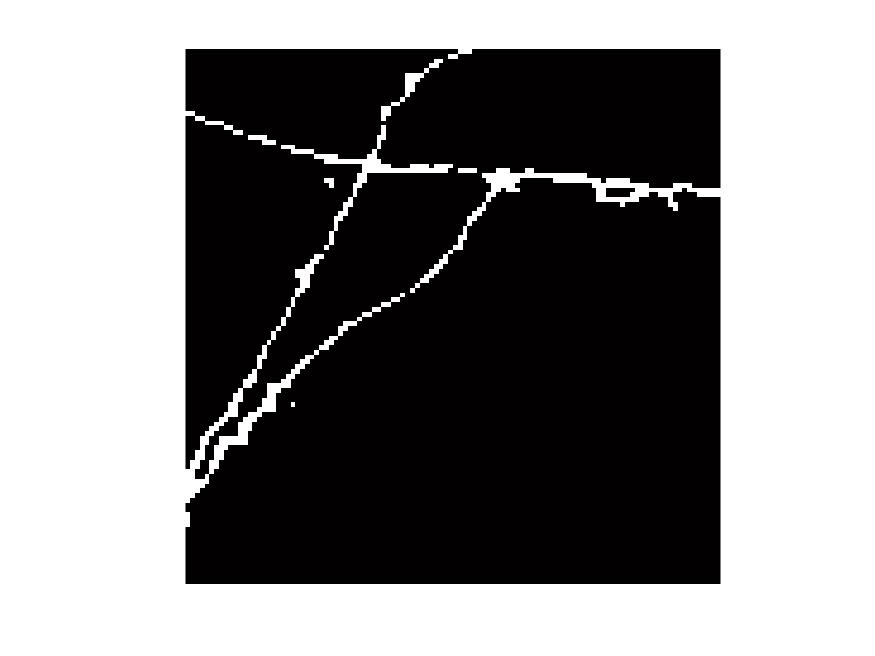}} & \hspace{-0.7cm}\adjustbox{valign=m,vspace=1pt}{\includegraphics[width=.28\textwidth]{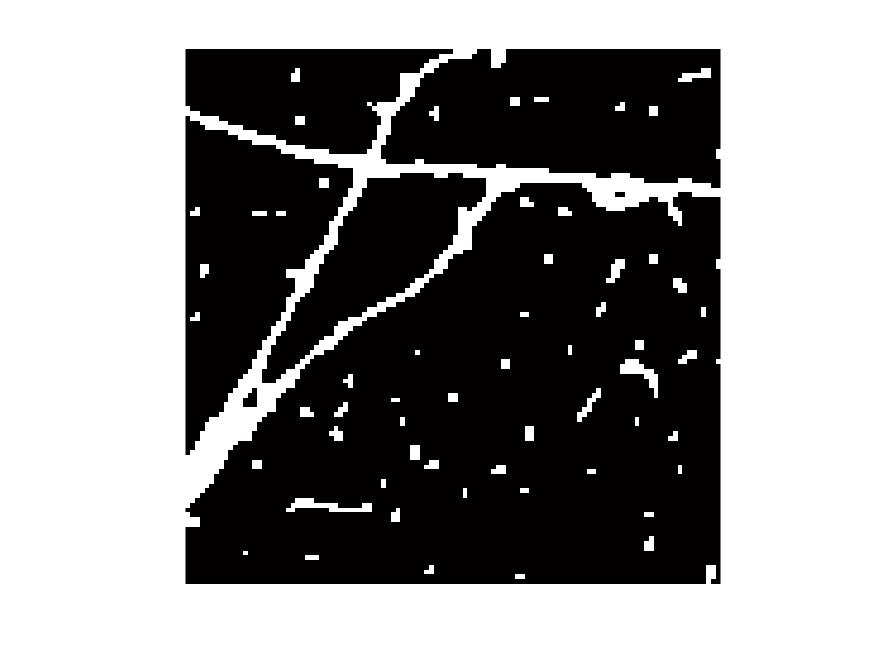}}\\
\end{tabular}
 \caption{(a) Noisy simulated dataset with low-background (LB) and stack size: $T=500$ frames, (b) Noisy simulated high-background (HB) dataset, with $T=500$ frames. From left to right: Superimposed diffraction limited image (temporal mean of the stack) with 4x zoom on ground truth support (blue), CEL0 reconstruction, $\ell_1$ reconstruction and TV reconstruction}
    \label{support_fig}
 \end{figure}

To evaluate the performance of the first step of the method COL0RME using the different regularization penalties described above, we created two noisy simulated datasets, with low background (LB) and high background (HB), respectively and used them to apply COL0RME and estimate the desired sample support. More details on the two datasets are available in the following sub-section \ref{sec: Simulated Data}. The results obtained by using the three different regularizers are reported in Figure \ref{support_fig}. In this example we chose the regularization parameter $\lambda$ heuristically, while more details about the selection of the parameter are given in the subsection \ref{parameter_lambda}.

Despite its continuous and smooth reconstruction, we observe that the reconstruction obtained by the TV regularizer does not provide precise localization results . For example, the separation of the two filaments on the top-right corner is not visible and while the junction of the other two filaments on the bottom-left should appear further down, we clearly see that those filaments are erroneously glued together. %For this reason, we will not further consider in the following the TV-regularizer but rather consider the other two regularizers by which more
 %precise 
 %\mdf{accurate}
 %localization \mdf{in specific areas of interest} can be achieved. 
 \mdf{Nonetheless, the choice of an appropriate regularizer tailored to favor fine structures as the ones observed in the GT image constitutes a challenging problem that should be addressed in future research. }
  
 The Jaccard indices (JI) of both the results obtained when using the CEL0 and $\ell_1$ regularizer, that allow for more precise localization, 
 %are presented in the Figure \ref{fig: JI}. 
 have been computed.
 The Jaccard index, is a quantity in the range $[0,1]$ computed as the ratio between correct detections (CD) and the sum of correct detections, false positives (FP) and false negatives (FN), that is $\text{JI}:=CD/ (CD+FN+FP)$, up to a tolerance $\delta>0$, measure in nm. \mdf{A correct detection occurs when one pixel at most $\delta$ nm away from a ground truth pixel is added to the support. In order to match the pixels from the estimated support to the ones from the ground truth, we employ the standard Gale–Shapley algorithm \cite{stable_maraige}. Once the matching has been performed, we can simply count the number of ground truth pixels which have not been detected (false negatives) and also the number of pixels in the estimated support which have not been matched to any ground truth pixel (false positives).}
 
\begin{figure}[H]
     \centering
     \begin{subfigure}[b]{0.44\textwidth}
         \centering
         \includegraphics[width=\textwidth]{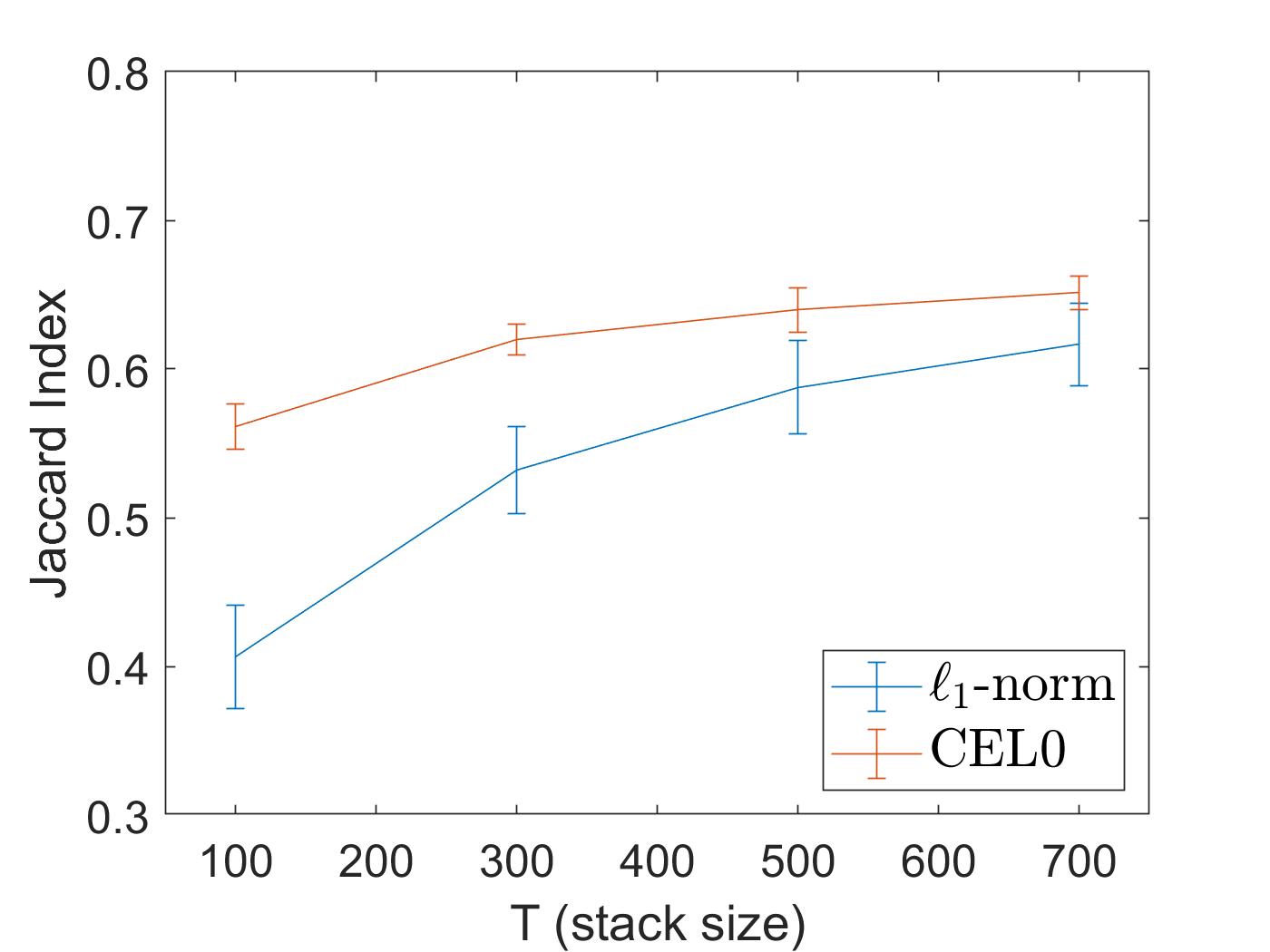}
         \caption{LB dataset}
     \end{subfigure}
     \hfill
     \begin{subfigure}[b]{0.44\textwidth}
         \centering
         \includegraphics[width=\textwidth]{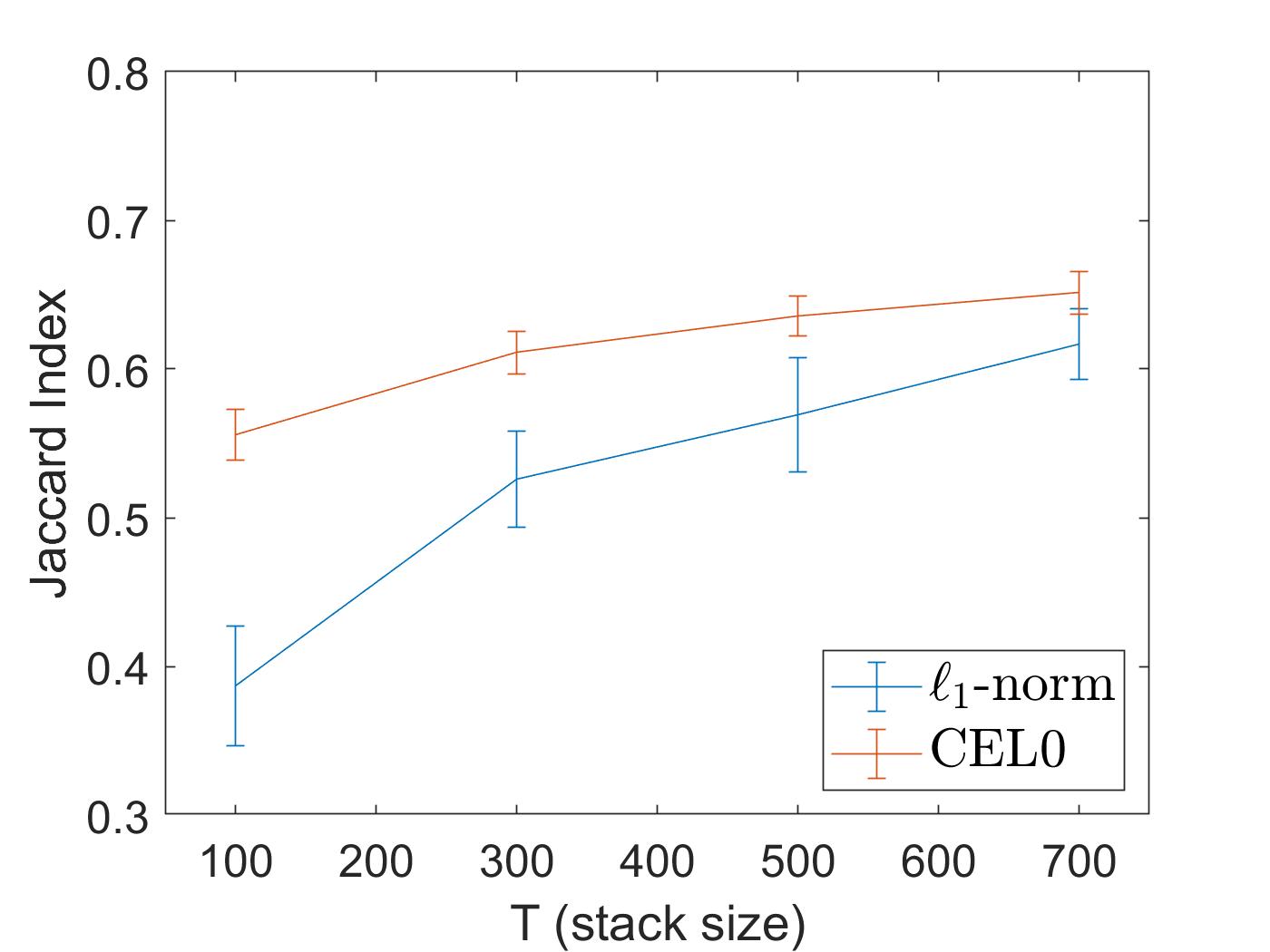}
         \caption{HB dataset}
     \end{subfigure}
     \caption{Jaccard Index values with tolerance $\delta = 40 nm$ for the low-background (LB) and high-background (HB) dataset, for different stack sizes and regularization penalty choices. The tolerance, $\delta = 40$ nm, is set so that we allow the correct detections, that needed to be counted for the computation of the Jaccard Index, to be found not only in the same pixel but also to any of the 8-neighbouring pixels}
     \label{fig: JI}
\end{figure}

 The Figure  \ref{fig: JI} reports the average Jaccard index computed from 20 different noise realizations, as well as, an error bar (vertical lines) that represent the standard deviation, for several stack sizes. According to the figure, a slightly better Jaccard index is obtained when the CEL0 regularizer is being used, while an increase in the number of frames, when both regularizers being used, leads to better Jaccard index, hence better localization. \mdfsec{As the reader may notice, such quantitative assessment could look inconsistent with the visual results reported in Figure \ref{support_fig}. By definition, the JI tends to assume higher values whenever more CD are found even in presence of more FP (as it happens for the CEL0 reconstruction), while it gets more penalized when FN happen, as they affect the computation "twice", reducing the numerator and increasing the denominator. }

\subsection{Accurate noise variance estimation}  \label{sec:noise_variance}
%\txtb{1. A graph showing how well we estimation $\sigma^2$ for different stack sizes. Compromise photobleaching - more information (more frames). For the best result, which is the stack size, how many molecules have bleached?\\
%2. An image while estimating the variance of the noise and while not - Big difference COL0RME+SPARCOM}

%  Along with the estimations of the emitter's temporal sparse covariance matrix, the estimation of the noise variance in the joint model \eqref{eq:support_mini} allows for much more precise results even in challenging acquisition conditions. In Figure \ref{fig:noise_var} we show how the relative error between the computed noise variance $s$ and the ground-truth one $\sigma^2$ used to produce simulated no-background (NB) data (more details in \ref{sec: Simulated Data}), decays as the number of temporal frames increases. In this example $15$ dB of Gaussian noise were used, while the value of \deleted{$\sigma$} \mdf{$\sigma^2$} in average used for the different stack sizes is equal to $7.11 \times 10^5$. Note that, in general, the estimation of the noise variance  obtained by COL0RME is very precise.

% \begin{figure}[H]
%      \centering
%          \centering
%          \includegraphics[width=0.5\textwidth]{Figures/error_noBg.jpg}
%         \caption{No-background dataset: the relative error in noise variance estimation, defined as: Error = $\frac{|s - \sigma^2|}{|\sigma^2|}$. The Error is computed for 20 different noise realizations, presenting in the graph the mean and the standard deviation (error bars).}
%         \label{fig:noise_var}
% \end{figure}

\begin{figure}[H]
\centering
     \begin{subfigure}[b]{0.4\textwidth}
         \centering
         \includegraphics[width=\textwidth]{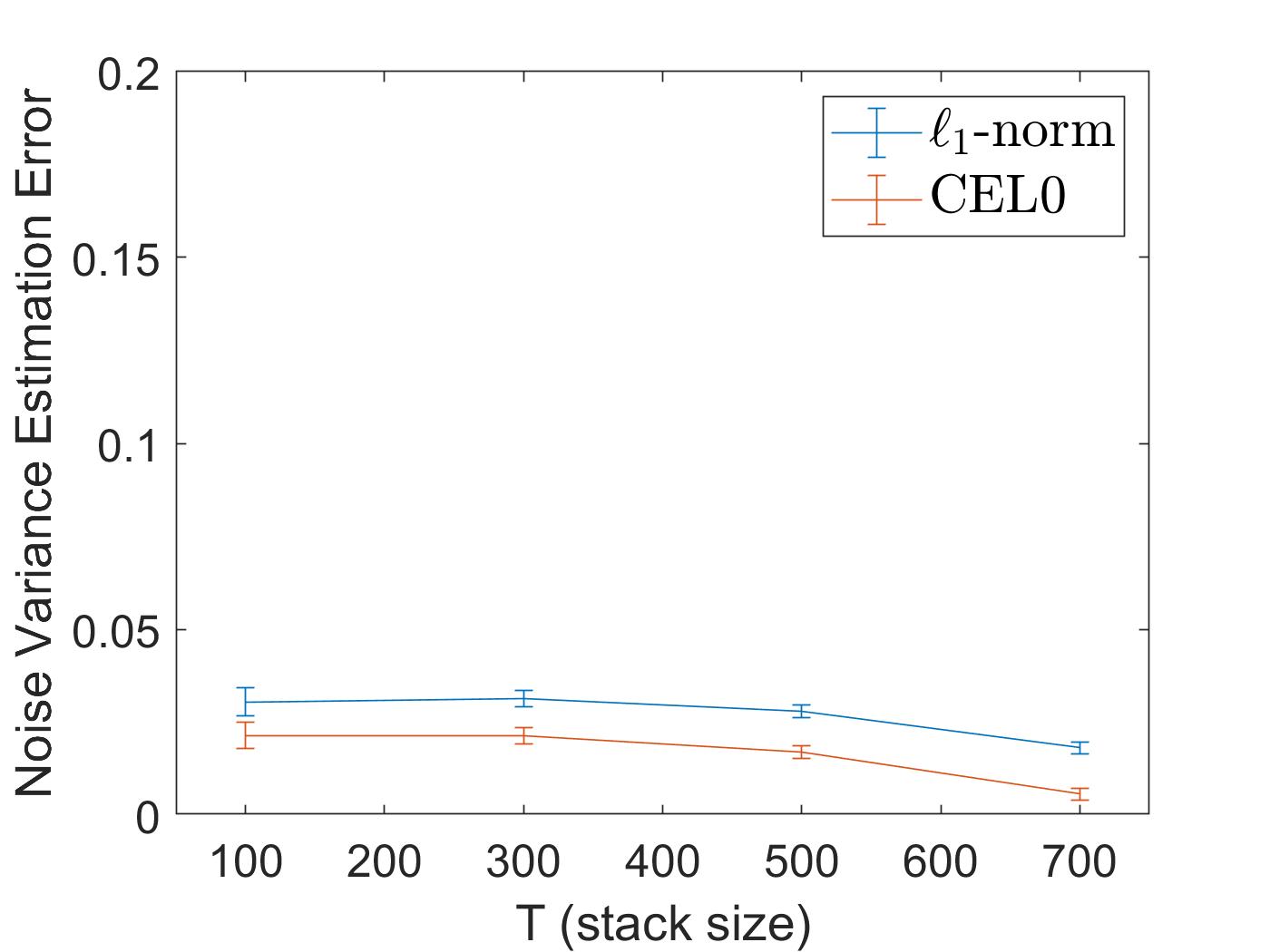}
         \caption{LB dataset}
     \end{subfigure}
     \begin{subfigure}[b]{0.4\textwidth}
         \centering
         \includegraphics[width=\textwidth]{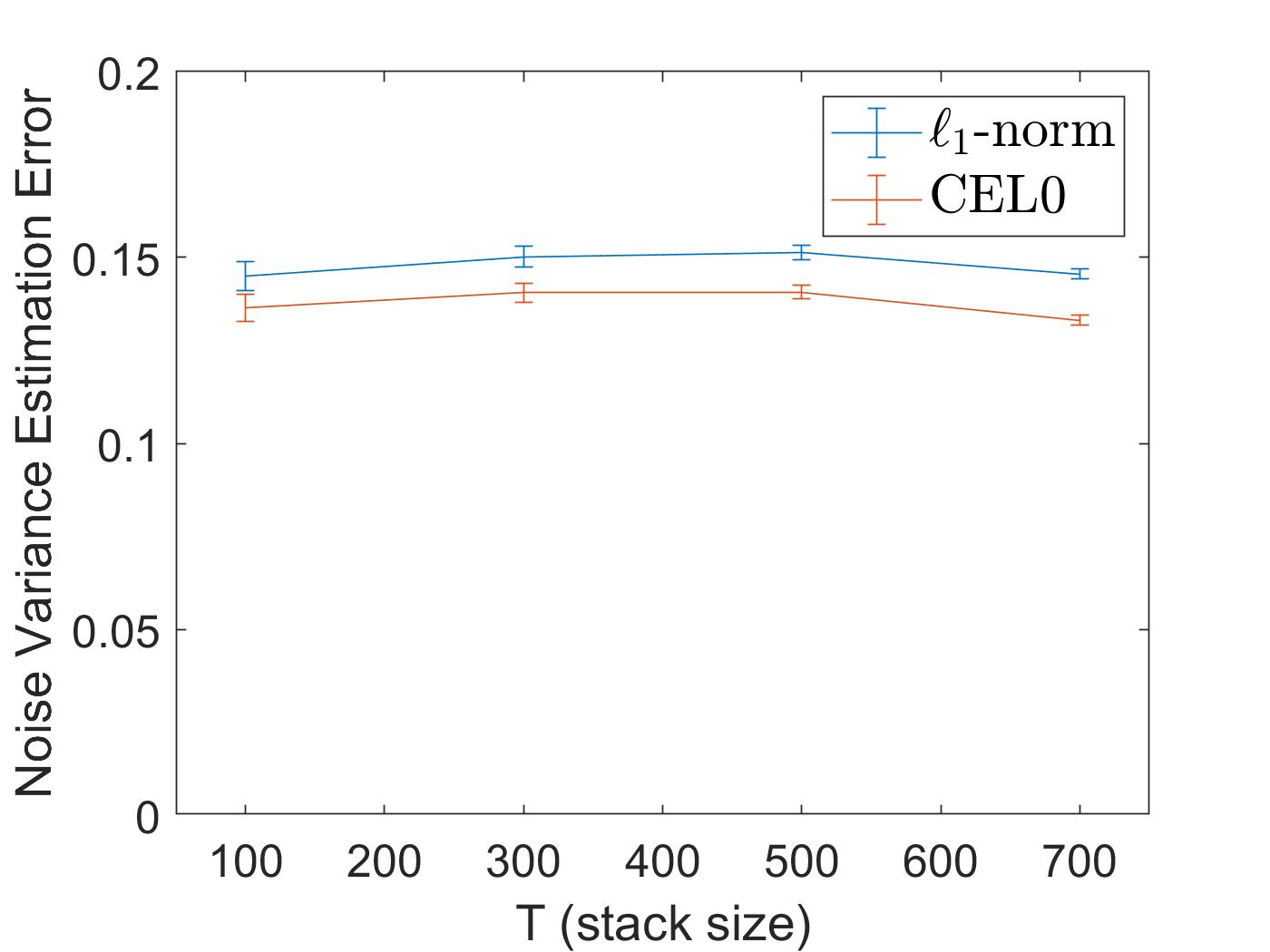}
         \caption{HB dataset}
     \end{subfigure}
        \caption{The relative error in noise variance estimation, defined as: Error = $\frac{|s - \sigma^2|}{|\sigma^2|}$, where $\sigma^2$ is the constant variance of the electronic noise. The Error is computed for 20 different noise realizations, presenting in the graph the mean and the standard deviation (error bars)}
        \label{fig:noise_var}
\end{figure}

 Along with the estimations of the emitter's temporal sparse covariance matrix, the estimation of the noise variance in the joint model \eqref{eq:support_mini} allows for much more precise results even in challenging acquisition conditions. \mdf{In Figure \ref{fig:noise_var} we show the relative error between the computed noise variance $s$ and the constant variance of the electronic noise $\sigma^2$ used to produce simulated low-background (LB) and high-background (HB) data. The relative error is higher in the case of the HB dataset, something that is expected, as in our noise variance estimation $s$ there is a bias coming from the background (see Remark \ref{rem:poisson model}). In the case of the LB dataset, as the background is low, the bias is sufficiently small so that it is barely visible in the error graph. In our experiments, a Gaussian noise with a corresponding SNR of approximately 16 dB is being used, while the value of \mdf{$\sigma^2$} is in average equal to $7.11 \times 10^5$ for the LB dataset and  $7.13 \times 10^5$ for the HB dataset. Note that, in general, the estimation of the noise variance $s$ obtained by COL0RME is very precise.}

%----------------------------------------------------------%

% COL0RME - Intensity estimation
\section{COL0RME, step II: Intensity estimation}  \label{sec:stepII}

From the previous step, we obtain a sparse estimation of $\mathbf{r_x}\in\mathbb{R}^{L^2}$. Its support, i.e. the location of non-zero variances, can thus be deduced. This is denoted in the following by  $\Omega := \left\{i: {(\mathbf{r_x})}_i\neq 0 \right\}\subset \left\{ 1,\ldots,L^2\right\}$. Note that this set corresponds indeed to the support of the desired $\mathbf{x}$, hence in the following we will use the same notation to denote both sets.

We are now interested in enriching COL0RME with an additional step where intensity information of the signal $\mathbf{x}$   can be retrieved in correspondence with the estimated support $\Omega$. 
To do so, we thus propose an intensity estimation procedure for $\mathbf{x}$ restricted only to the pixels of interest. Under this modeling assumption, it is thus reasonable to consider a regularization term favoring smooth intensities on $\Omega$, in agreement to the intensity typically found in real images. %Note that this choice further has the advantage of discouraging the appearance of isolated points resulting from possible errors in the support estimation step. 

In order to take into account the modeling of blurry and out-of-focus fluorescent molecules, we further include in our model \eqref{eq:model_vec} a regularization term for smooth background estimation.
%that describes , 
%a regularization term is needed to get the well-posedness of the problem and to enforce smoothness. 
We can thus consider the following joint minimization problem:
\begin{equation}{\label{eq:intensity_constraint}}
    \argmin\limits_{\mathbf{x}\in\mathbb{R}_+^{|\Omega|},~ \mathbf{b}\in\mathbb{R}_+^{M^2}} ~\frac12 \|\mathbf{\Psi_\Omega} \mathbf{x} - (\overline{\mathbf{y}} - \mathbf{b})\|_2^2 + \frac{\mu}{2} \|\nabla_{\Omega}\mathbf{x}\|_2^2 + \frac{\beta}{2} \|\nabla\mathbf{b}\|_2^2,
\end{equation}
where the data term models the presence of Gaussian noise, 
%$\overline{\mathbf{y}} = \sum_{t=1}^T\mathbf{y_t}$ 
%is the temporal sampling average of the acquired data   % we said this already!
$\mu,\beta>0$ are regularization parameters 
%The parameter $\beta$ is typically chosen big enough to enforce  background smoothness.
%\txtb{i think that depending on the data $\beta$ could well be small. maybe we should not say anything here and comment later on the choice of $\beta$ in the parameter estimation section?}
and the operator $\mathbf{\Psi_\Omega} \in\mathbb{R}^{M^2\times |\Omega|}$ is a matrix whose $i$-th column is extracted from $\mathbf\Psi$ for all indexes $i \in \Omega$. Finally, the regularization term on $\mathbf{x}$ is the squared norm of the discrete gradient restricted to $\Omega$, i.e.:
\begin{equation*}
    \|\nabla_{\Omega} \mathbf{x}\|_2^2 := \sum\limits_{i \in \Omega} \sum\limits_{j \in \mathcal{N}(i)\cap\Omega} (x_i - x_j)^2,
\end{equation*}
where $\mathcal{N}(i)$ denotes the 8-pixel neighborhood of $i\in\Omega$. Note that, according to this definition, $\nabla_\Omega\mathbf{x}$ denotes a (redundant) isotropic discretization of the gradient of $\mathbf{x}$ evaluated for each pixel in the support $\Omega$. Note that this definition coincides with the standard one for $\nabla \mathbf{x}$ restricted to points in the support $\Omega$.

The non-negativity constraints on $\mathbf{x}$ and $\mathbf{b}$ as well as the one restricting the estimation of $\mathbf{x}$  on $\Omega$ can be relaxed by using suitable \mdffirst{smooth} penalty terms, so that, finally, the following optimization problem can be addressed:
\begin{equation}\label{eq:intensity_penalized}
    \argmin\limits_{\mathbf{x}\in\mathbb{R}^{L^2},~ \mathbf{b}\in\mathbb{R}^{M^2}} ~\frac12 \| \mathbf{\Psi} \mathbf{x} - (\overline{\mathbf{y}} - \mathbf{b})\|_2^2 + \frac{\mu}{2} \|\nabla\mathbf{x}\|_2^2 + \frac{\beta}{2} \|\nabla\mathbf{b}\|_2^2 + \frac{\alpha}{2} \left(\|\mathbf{I_\Omega x}\|_2^2 + \sum_{i=1}^{L^2}\ [\phi(\mathbf{x}_i)]^2 + \sum_{i=1}^{M^2}\ [\phi(\mathbf{b}_i)]^2\right),
\end{equation}
where the parameter $\alpha\gg 1$ can be chosen arbitrarily high to enforce the constraints, $\mathbf{I_\Omega}$ is a diagonal matrix acting as characteristic function of $\Omega$, i.e. defined as: 
    \[ \mathbf{I_\Omega}(i,i) = \begin{cases} \mbox{0} & \mbox{if } i \in \Omega, \\ \mbox{1} & \mbox{if } i \not\in \Omega \end{cases},\qquad \forall i \in \{1, ... , L^2\}, \]
    and $\phi: \mathbb{R}\to\mathbb{R}$ is used to penalize negative entries, being defined as:
    \begin{equation}
        \phi(z) := 
        \begin{cases}
          0 & \text{if } z \geq 0,\\
          z & \text{if } z < 0
        \end{cases},\qquad\forall z\in\mathbb{R}.
        \label{phi}
    \end{equation}
\mdffirst{We anticipate here that considering the unconstrained problem \eqref{eq:intensity_penalized} instead of the original, constrained, one \eqref{eq:intensity_constraint}, will come in handy for the design of an automatic parameter selection strategy, as we further detail in Section \ref{sec: DP}. }

To solve the joint-minimization problem \eqref{eq:intensity_penalized} we use the Alternate Minimization algorithm, see  Algorithm~\ref{Algorithm:AMA_intensity}. In the following subsections, we provide more details on the solution of the  two minimization sub-problems.

\begin{algorithm}[!h]
\caption{COL0RME, Step II: Intensity Estimation}
\label{Algorithm:AMA_intensity}
\begin{algorithmic}
\REQUIRE $\overline{\mathbf{y}}\in\mathbb{R}^{M^2}, \mathbf{x}^0\in\mathbb{R}^{L^2},\mathbf{b}^0\in\mathbb{R}^{M^2}, \mu,\beta>0$, $\alpha\gg 1$
\REPEAT 
\STATE $\mathbf{x}^{k+1} = \argmin\limits_{\mathbf{x}\in \mathbb{R}^{L^2}} \frac12\|\mathbf{\Psi x - (\overline{\mathbf{y}} - \mathbf{b}^{k})} \|_2^2
+\frac{\mu}{2} \|\nabla\mathbf{x}\|_2^2 
+\frac{\alpha}{2} \left(\|\mathbf{I_\Omega x}\|_2^2 
+\sum_{i=1}^{L^2}\ [\phi(\mathbf{x}_i)]^2\right)$
\STATE $\mathbf{b}^{k+1} = \argmin\limits_{\mathbf{b} \in \mathbb{R}^{M^2}}\frac12 \| \mathbf{b} - (\overline{\mathbf{y}}-\mathbf{\Psi} \mathbf{x}^{k+1}) \|_2^2 +\frac{\beta}{2} \|\nabla\mathbf{b}\|_2^2 +\frac{\alpha}{2}\sum_{i=1}^{M^2}\ [\phi(\mathbf{b}_i)]^2$
\UNTIL convergence
\RETURN $\mathbf{x},\mathbf{b}$
\end{algorithmic}
\end{algorithm}

\subsection{First sub-problem: update of $\mathbf{x}$}

In order to find at each $k\geq 1$ the optimal solution $\mathbf{x}^{k+1} \in \mathbb{R}^{L^2}$ for the first sub-problem, we need to solve a minimization problem of the form:
\begin{equation}
     \mathbf{x}^{k+1} = \argmin_{\mathbf{x}\in \mathbb{R}^{L^2}} ~g(\mathbf{x};\mathbf{b}^k) + h(\mathbf{x}),
     \label{min_x}
\end{equation} 
where, for $\mathbf{b}^k\in\mathbb{R}^{M^2}$ being fixed at each iteration $k\geq 1$, $g(
\cdot; \mathbf{b}^k): \mathbb{R}^{M^2} \rightarrow \mathbb{R}_+$ is a proper and convex function with Lipschitz gradient,  defined as:
\begin{equation}  \label{eq:def_g}
    g(\mathbf{x}; \mathbf{b}^k) := \frac12\|\mathbf{\Psi x} - (\overline{\mathbf{y}} - \mathbf{b}^{k}) \|_2^2 
+\frac{\mu}{2} \|\nabla \mathbf{x}\|_2^2, 
\end{equation}
and where the function $h : \mathbb{R}^{L^2} \rightarrow \mathbb{R}$ encodes the penalty terms:
\begin{equation}
    h(\mathbf{x}) = \frac{\alpha}{2} \left(\|\mathbf{I_\Omega x}\|_2^2 
+\sum_{i=1}^{L^2}\ [\phi(\mathbf{x}_i)]^2\right).
\label{eq: h}
\end{equation}

Solution of \eqref{min_x} can be obtained iteratively, using, for instance, the proximal gradient descent algorithm, whose iteration can be defined as follows :
\begin{equation}
    \mathbf{x}^{n+1} = \textbf{\text{prox}}_{h, \tau}(\mathbf{x}^{n} - \tau \nabla g(\mathbf{x}^{n})), \quad n=1,2,.. ,
    \label{x_m}
\end{equation}
where $\nabla g(\cdot)$ denotes the gradient of $g$, $\tau \in (0, \frac{1}{L_{g}}]$ is the algorithmic step-size chosen inside a range depending on the  Lipschitz constant of $\nabla g$, here denoted by $L_g$, to guarantee convergence. The proximal update in \eqref{x_m} can be computed explicitly using the computations reported in Appendix \ref{appendixA}. One can show in fact that, for each $\mathbf{w}\in\mathbb{R}^{L^2}$ there holds element-wise:
\begin{equation}  \label{eq:prox}
    \left( \textbf{\text{prox}}_{h, \tau}(\mathbf{w}) \right)_i = {\text{prox}}_{h, \tau}(\mathbf{w}_i) = 
    \begin{cases}
      \frac{\mathbf{w}_i}{1 + \alpha \tau \mathbf{I_\Omega}(i,i)} & \text{if } {\mathbf{w}_i} \geq 0,\\
      \frac{\mathbf{w}_i}{1 + \alpha \tau ( \mathbf{I_\Omega}(i,i)+1)} & \text{if } \mathbf{w}_i < 0.
    \end{cases}       
\end{equation}

%The quantity $\nabla g(\mathbf{x})$ is equal to:
%\begin{equation*}
%    \nabla g(\mathbf{x}) = (\mathbf\Psi^\intercal \mathbf\Psi +\mu \mathbf{D}^\intercal \mathbf{D}) \mathbf{x} - \mathbf\Psi^\intercal %(\overline{\mathbf{y}} - \mathbf{b})
%\end{equation*}
%
%and the step $t \in (0, \frac{1}{Lip}]$, where the Lipschitz constant: $Lip = \| \mathbf\Psi^\intercal \mathbf\Psi +\mu \mathbf{D}^\intercal %\mathbf{D}\|_2$ 
 
\begin{remark}  \label{rem:proximal}
As the reader may have noted, we consider the proximal gradient descent algorithm \eqref{x_m} for solving \eqref{min_x}, even though both functions $g$ and $h$ in \eqref{eq:def_g} and \eqref{eq: h} respectively, are smooth and convex, hence, in principle, (accelerated) gradient descent algorithms could be used. Note, however, that the presence of the large penalty parameter $\alpha\gg 1$ would significantly slow down convergence speed in such case as the step size $\tau$ in this case would be constrained to the smaller range $(0, 
\frac{1}{L_g + \alpha}]$. By considering the penalty contributions in terms of their proximal operators, this limitation doesn't affect the range of $\tau$ and convergence is still guaranteed\cite{Combettes2005} in a computationally fast way through the update \eqref{eq:prox}.
 \end{remark}
 
 %\txtb{Maybe add here a reference to the paper of Combettes Wajs for convergence of proximal gradient?}

\subsection{Second sub-problem: update of $\mathbf{b}$}

As far as the estimation of the background is concerned, the minimization problem we aim to solve at each $k\geq 1$ takes the form:
\begin{equation}
    \mathbf{b}^{k+1} = \argmin_{\mathbf{b} \in \mathbb{R}^{M^2}}~ r(\mathbf{b};\mathbf{x}^{k+1}) + q(\mathbf{b}),
    \label{eq: min_varback}
\end{equation}
where:
\begin{equation*}
    r(\mathbf{b};\mathbf{x}^{k+1}) := \frac12 \| \mathbf{b} - (\overline{\mathbf{y}}-\mathbf{\Psi} \mathbf{x}^{k+1}) \|_2^2 +\frac{\beta}{2} \|\nabla\mathbf{b}\|_2^2,\qquad     q({\mathbf{b}}) := \frac{\alpha}{2}\sum_{i=1}^{M^2}\ [\phi(\mathbf{b}_i)]^2.
\end{equation*}
Note that $r(\cdot;\mathbf{x}^{k+1}):\mathbb{R}^{M^2}\to\mathbb{R}_+$ is a convex function with $L_r$-Lipschitz gradient and $q:\mathbb{R}^{M^2}\to\mathbb{R}_+$ encodes (large, depending on $\alpha\gg 1$) penalty contributions. Recalling Remark \ref{rem:proximal},  we thus use again the proximal gradient descent algorithm for solving \eqref{eq: min_varback}. The desired solution $\hat{\mathbf{b}}$ at each $k\geq 1$ can thus be found by iterating:
\begin{equation}
    \mathbf{b}^{n+1} = \textbf{\text{prox}}_{q, \delta}(\mathbf{b}^{n} - \delta \nabla r(\mathbf{b}^{n})), \quad n=1,2,.. ,
    \label{b}
\end{equation}
for $\delta\in(0,\frac{1}{L_r}]$.
The proximal operator $\textbf{\text{prox}}_{q, \delta}(\cdot)$, has an explicit expression and it is defined element-wise for $i=1,\ldots,M^2$ as:
\begin{equation}
  \left( \textbf{\text{prox}}_{q, \delta}(\mathbf{d}) \right)_i = {\text{prox}}_{q, \delta}(\mathbf{d}_i) = 
\begin{cases}
  \mathbf{d}_i & \text{if } \mathbf{d}_i \geq 0,\\
  \frac{\mathbf{d}_i}{1 + \alpha \delta} & \text{if } \mathbf{d}_i < 0.
\end{cases}
\end{equation}
%\begin{equation*}
%    \nabla r(\mathbf{b}) = \mathbf{\Psi x} -(\overline{\mathbf{y}} - \mathbf{b}) + \beta \mathbf{D}^\intercal \mathbf{D} \mathbf{b},
%\end{equation*}

%with $\text{Lip} = \| \mathbf{I_{M^2}} + \beta \mathbf{D}^\intercal \mathbf{D}\|_2$ being the Lipschitz constant.% We can compute, fast: $t = \frac{1}{\|\mathbf{I_{N^2}}\|_2 + \beta \|\mathbf{D}^\intercal \mathbf{D}\|_2 } = \frac{1}{1 + \beta \|\mathbf{D}^\intercal \mathbf{D}\|_2 }$, result which we can obtain using the triangular inequality.

\subsection{Intensity and Background estimation results}

\begin{figure}
\centering
\setlength\tabcolsep{1.5pt}
\begin{tabular}{cccc}
 & \hspace{-0.3cm} $\bar{\mathbf{y}} (LB)$ & \hspace{-0.4cm}$\bar{\mathbf{y}} (HB)$ & \hspace{-0.4cm}$\mathbf{x}^{GT}$ \\
    &\hspace{-0.1cm}\adjustbox{valign=m,vspace=1pt}{\includegraphics[width=0.28\textwidth]{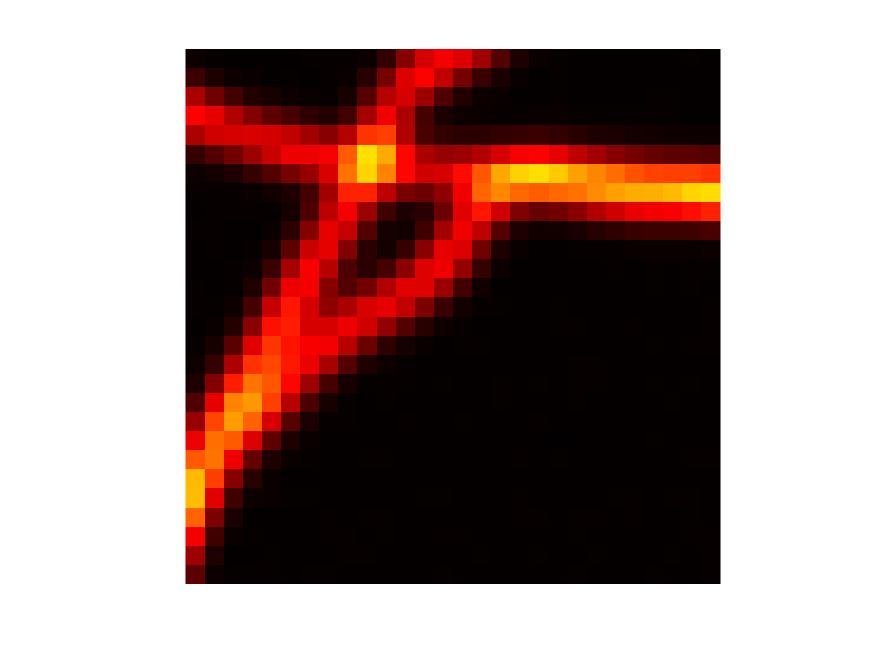}} &\hspace{-0.4cm}\adjustbox{valign=m,vspace=1pt}{\includegraphics[width=0.28\textwidth]{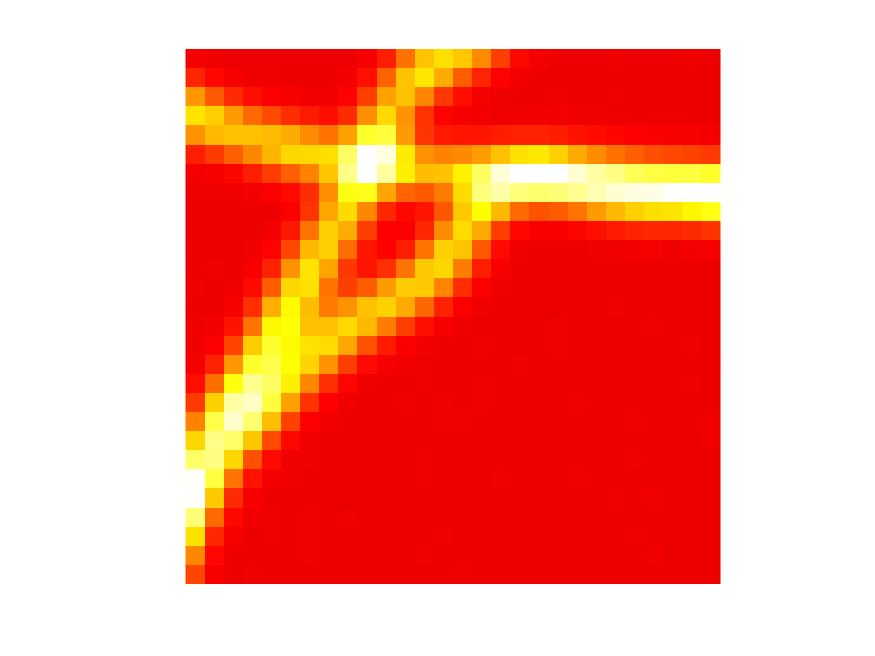}} 
    &\hspace{-0.4cm}\adjustbox{valign=m,vspace=1pt}{\includegraphics[width=0.28\textwidth]{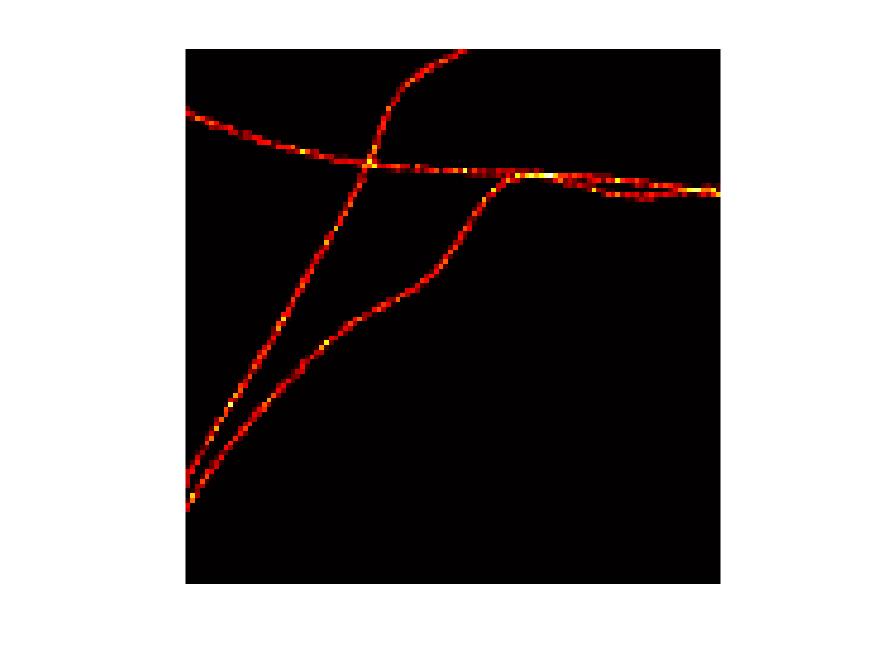}}\\
    & \multicolumn{2}{c}{\adjustbox{valign=m,vspace=1pt}{\includegraphics[width=.28\textwidth]{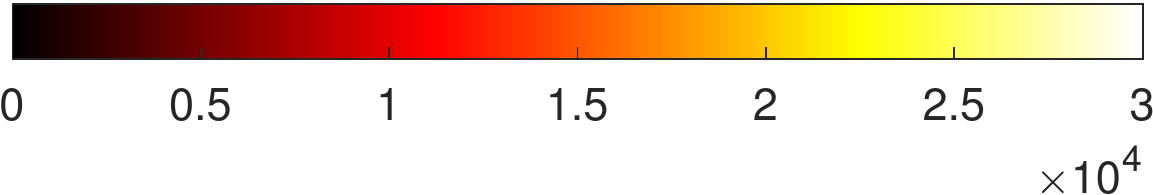}} }& \hspace{-0.2cm}\adjustbox{valign=m,vspace=1pt}{\includegraphics[width=.28\textwidth]{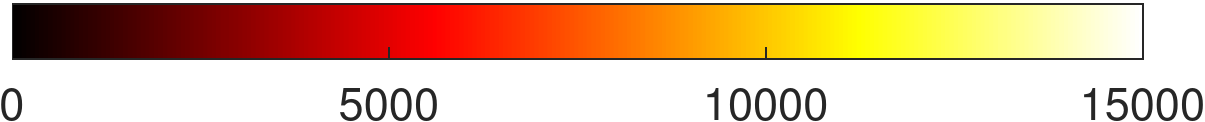}}\\
    & \hspace{-0.1cm} COL0RME-CEL0 & \hspace{-0.4cm}COL0RME-$\ell_1$ & \hspace{-0.4cm}COL0RME-TV \\
(a) & \hspace{-0.1cm}\adjustbox{valign=m,vspace=1pt}{\includegraphics[width=.28\textwidth]{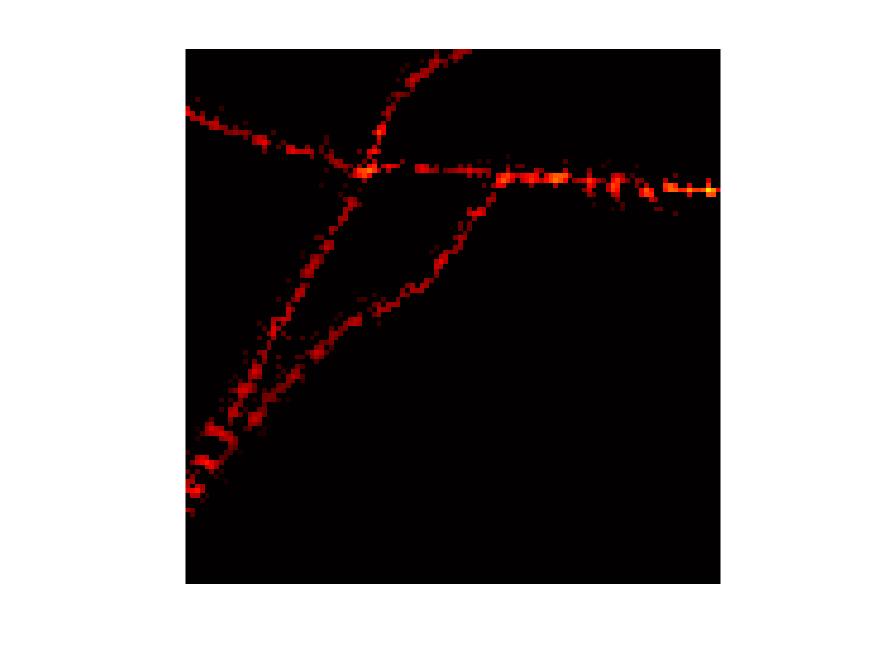}} 
& \hspace{-0.4cm}\adjustbox{valign=m,vspace=1pt}{\includegraphics[width=.28\textwidth]{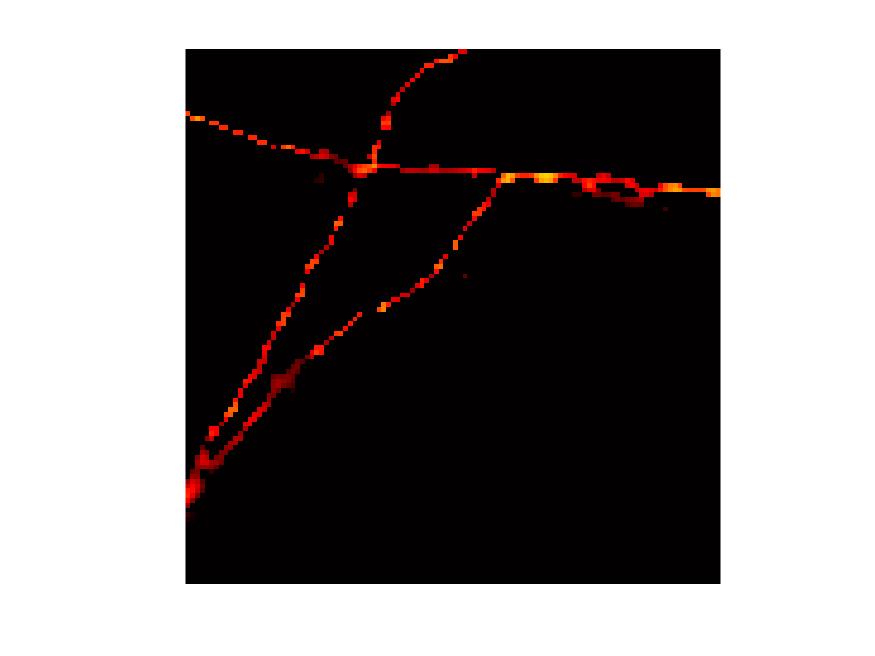}} 
& \hspace{-0.4cm}\adjustbox{valign=m,vspace=1pt}{\includegraphics[width=.28\textwidth]{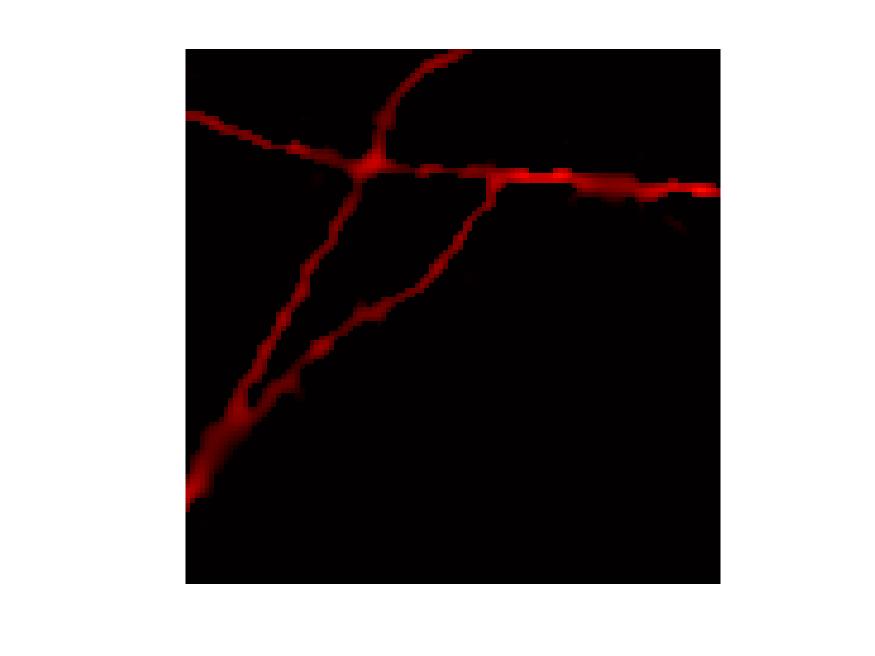}} \\
(b) &\hspace{-0.1cm}\adjustbox{valign=m,vspace=1pt}{\includegraphics[width=.28\textwidth]{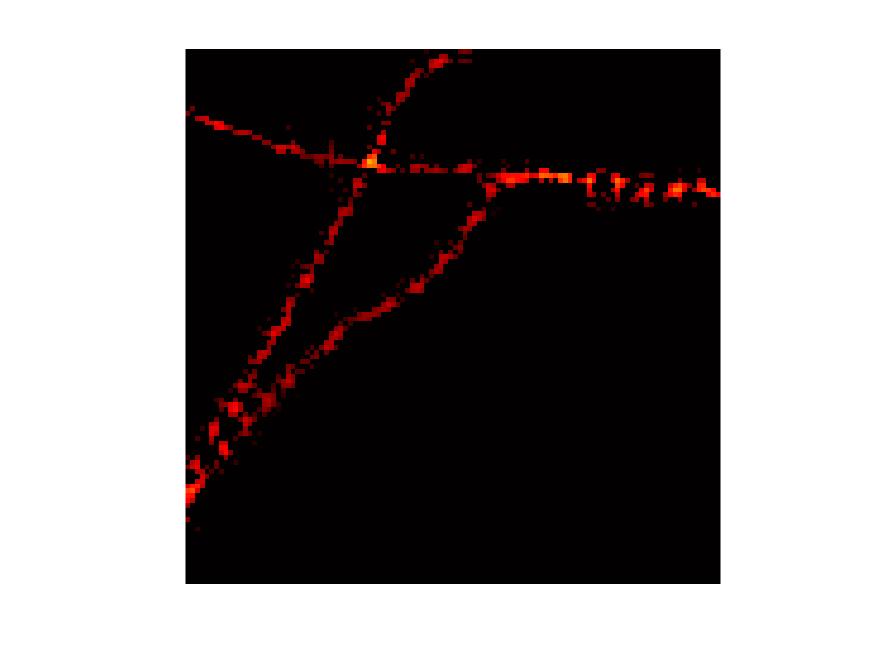}} & \hspace{-0.4cm}\adjustbox{valign=m,vspace=1pt}{\includegraphics[width=.28\textwidth]{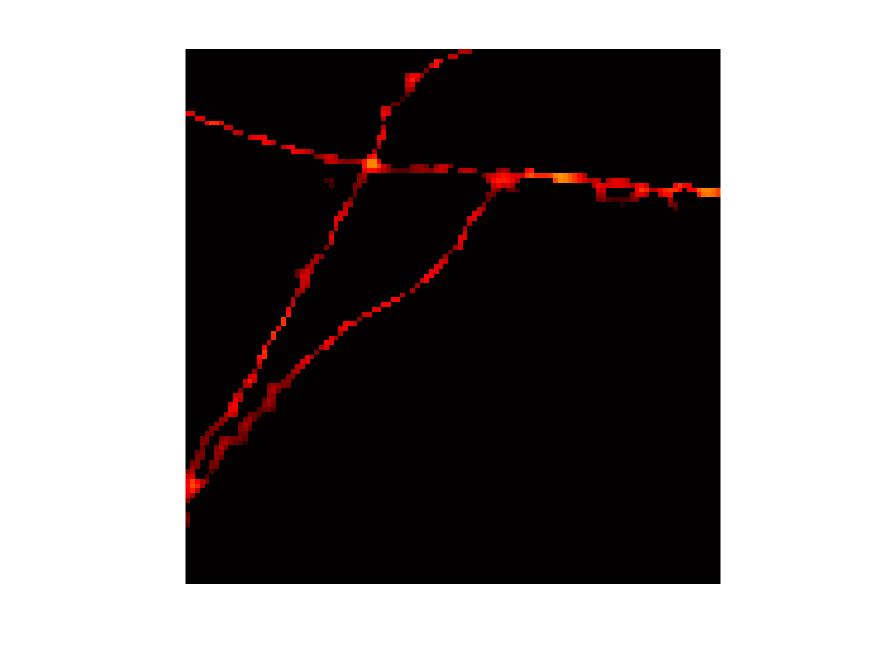}} & \hspace{-0.4cm}\adjustbox{valign=m,vspace=1pt}{\includegraphics[width=.28\textwidth]{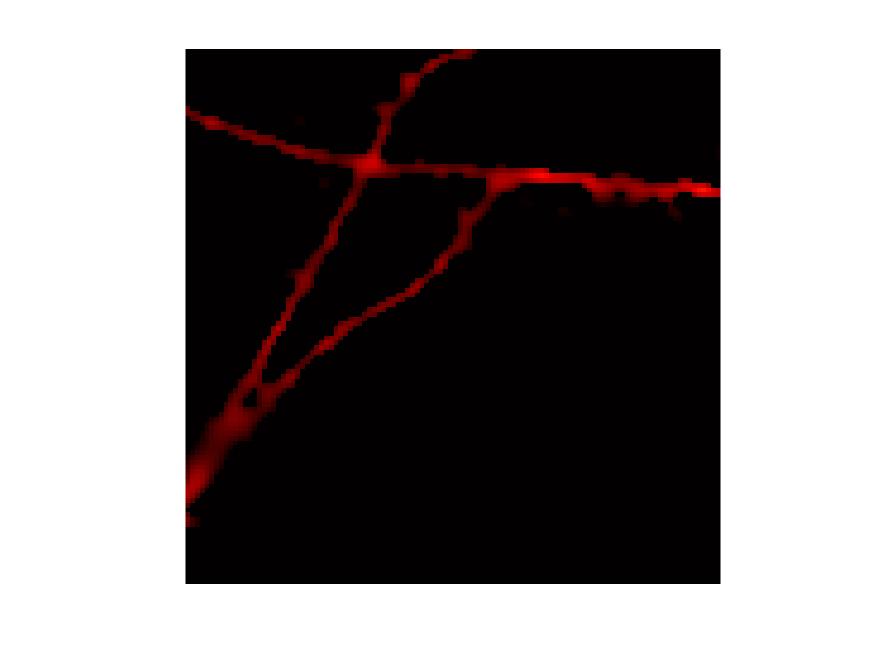}}\\
& &  \hspace{-0.2cm}\adjustbox{valign=m,vspace=1pt}{\includegraphics[width=.28\textwidth]{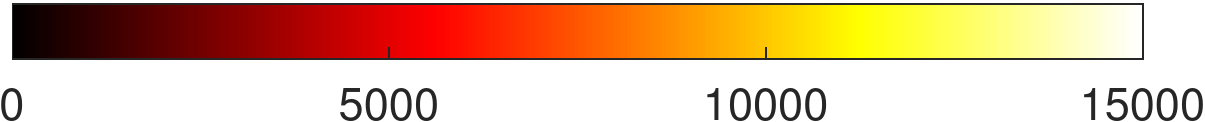}}&\\
\\
\end{tabular}
 \caption{On top: Diffraction limited image $\bar{\mathbf{y}}=\frac{1}{T}\sum_{t=1}^T \mathbf{y_t}$, with T=500, (4x zoom) for the low-background (LB) dataset and for the high-background (HB) dataset, Ground truth (GT) intensity image. (a) Reconstructions for the noisy simulated dataset with low-background (LB), (b) Reconstruction for the noisy simulated dataset with high-background (HB). From left to right: intensity estimation result on estimated support using CEL0 regularization, $\ell_1$ regularization and TV regularization. For all COL0RME intensity estimations, the same colorbar, presented at the bottom of the figure, has been used} 
 %(see Figure  \ref{fig:colorbars} for colorbars)
 
    \label{intensity_fig}
 \end{figure}

Intensity estimation results can be found in Figure \ref{intensity_fig} where \eqref{eq:intensity_penalized} is used for intensity/background estimation on the supports $\Omega_{\cal{R}}$ estimated from the first step of COL0RME using ${\cal{R}}=$ CEL0, ${\cal{R}}=\ell_1$ and ${\cal{R}}=$ TV. We are referring to them as COL0RME-CEL0, COL0RME-$\ell_1$ and COL0RME-TV, respectively. The colormap ranges are different for the coarse-grid and fine-grid representations, as explained in section \ref{sec: Simulated Data}
%, and their colorbars can be found in Figure \ref{fig:colorbars}. 
The result on $\Omega_{TV}$, even after the second step does not allow for the observation of a few significant details (e.g. the separation of the two filament on  the bottom left corner) and that is why it will not further discussed.
 
\begin{figure}[ht]
     \centering
     \begin{subfigure}[b]{0.44\textwidth}
         \centering
         \includegraphics[width=\textwidth]{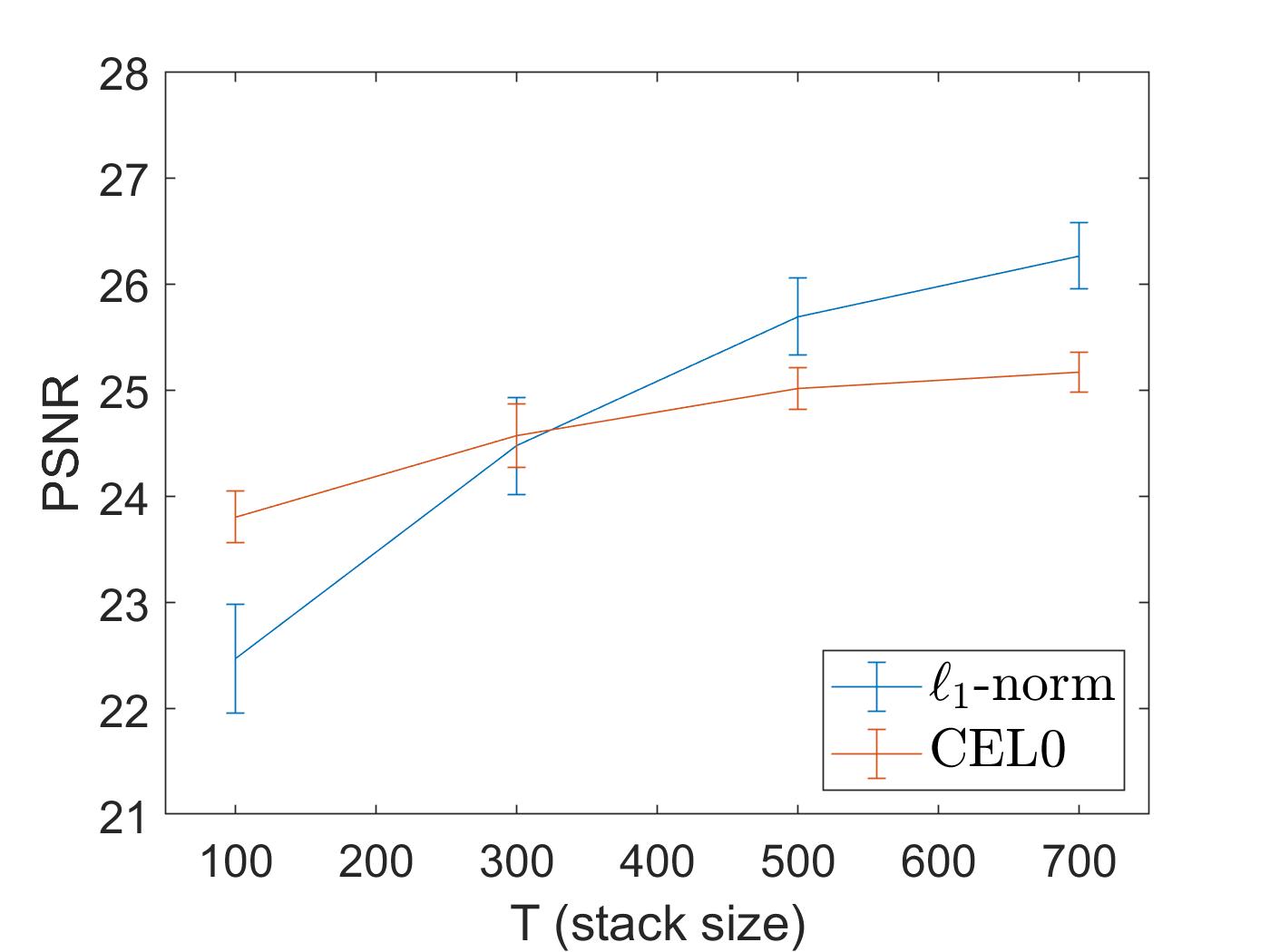}
         \caption{LB dataset}
     \end{subfigure}
     \hfill
     \begin{subfigure}[b]{0.44\textwidth}
         \centering
         \includegraphics[width=\textwidth]{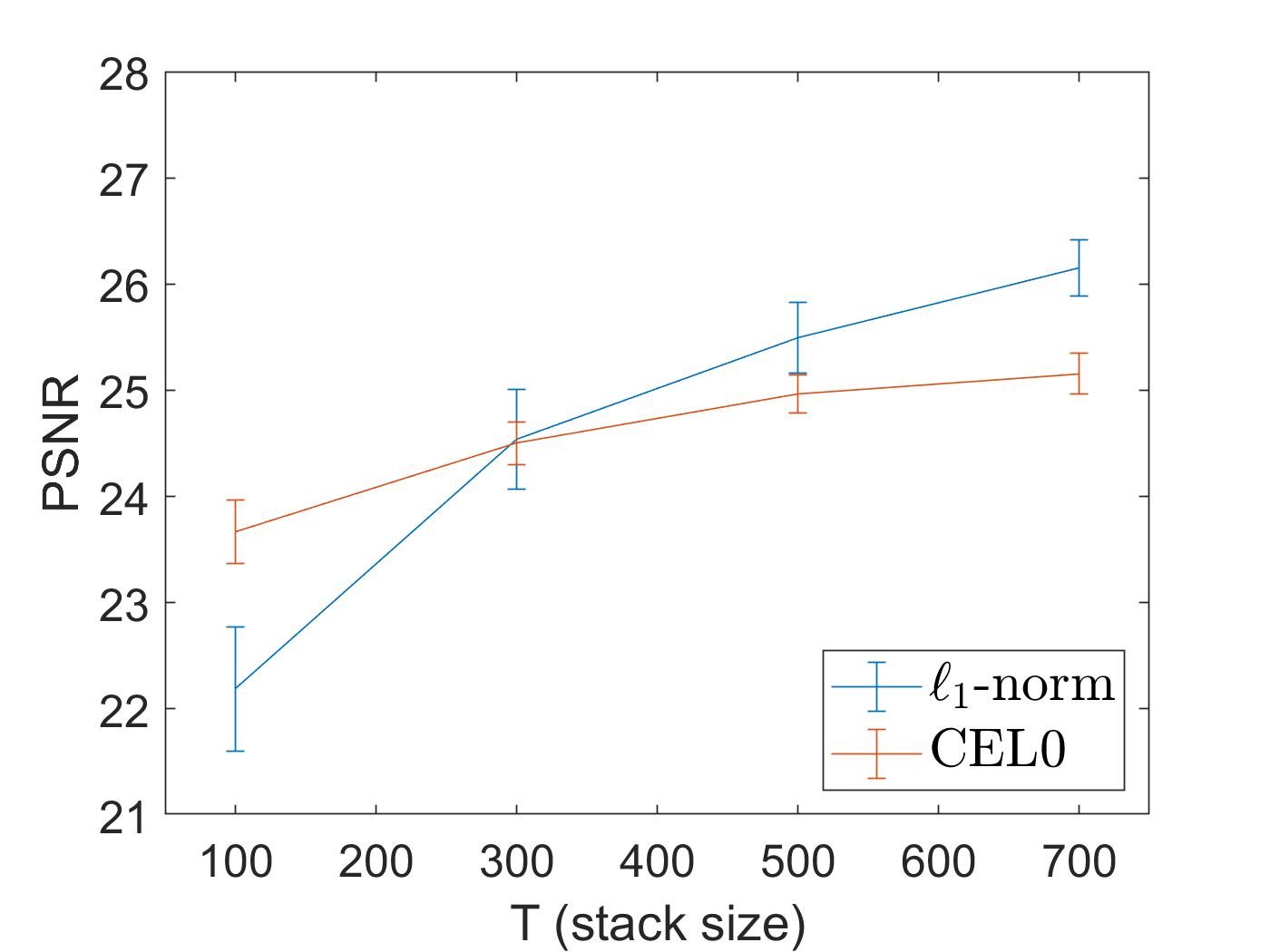}
         \caption{HB dataset}
     \end{subfigure}
     \caption{COL0RME PSNR values for two different datasets (low-background and high-background dataset), stack sizes and regularization penalty choices. The mean and the standard deviation of 20 different noise realizations are presented}
     \label{psnr}
\end{figure}

 A quantitative assessment for the other two regularization penalty choices, $\Omega_{CEL0}$ and $\Omega_{\ell_1}$, is available in Figure \ref{psnr}. More precisely we compute the Peak-Signal-to-Noise-Ratio (PSNR), given the following formula:
 \begin{equation}
     \text{PSNR}_\text{dB} = 10 \log_{10} \left( \frac{\text{MAX}^2_\mathbf{R}}{\text{MSE}}\right), \qquad \text{MSE} = \frac{1}{L^2}\sum\limits_{i=1}^{L^2}\left(\mathbf{R}_i - \mathbf{K}_i\right)^2,
 \end{equation}
where $\mathbf{R}\in \mathbb{R}^{L^2}$ is the reference image, $\mathbf{K}\in \mathbb{R}^{L^2}$ the image we want to evaluate using the PSNR metric and $\text{MAX}_\mathbf{R}$ the maximum value of the image $\mathbf{R}$. In our case, the reference image is the ground truth intensity image: $\mathbf{x}^{GT} \in \mathbb{R}^{L^2}$. The higher the PSNR, the better the quality of the reconstructed image.  

According to Figures \ref{intensity_fig} and \ref{psnr}, when only a few frames are considered (eg. $T = 100$ frames, high temporal resolution), the method performs better by using the CEL0 penalty for the support estimation. However, when longer temporal sequences are available (e.g. $T=500$ or $T = 700$ frames) the method performs better by using the $\ell_1$-norm instead. In addition to this, for both penalizations, PSNR improves as the number of temporal frames increases.

\begin{figure}
\centering
\setlength\tabcolsep{1.5pt}
\begin{tabular}{ccccc}
    & $\bar{\mathbf{y}}$  & \hspace{-0.3cm} {\small{COL0RME-CEL0}} ($\mathbf{b}$) & \hspace{-0.3cm}{\small{COL0RME-$\ell_1$}} ($\mathbf{b}$) & \hspace{-0.3cm}$\mathbf{b}^{GT}$ \\
(a) & \adjustbox{valign=m,vspace=1pt}{\includegraphics[width=0.25\textwidth]{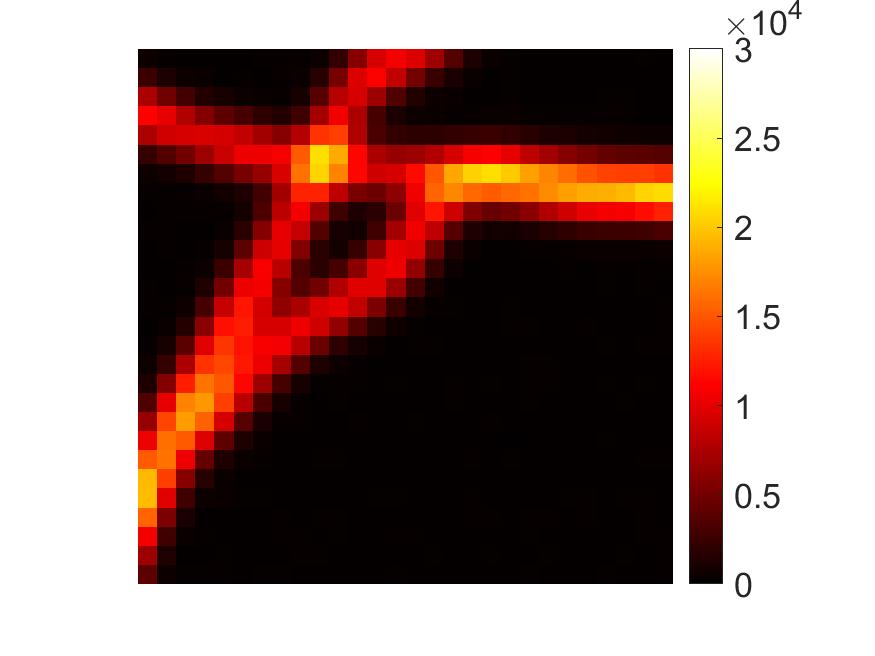}}  & \hspace{-0.3cm}\adjustbox{valign=m,vspace=1pt}{\includegraphics[width=.25\textwidth]{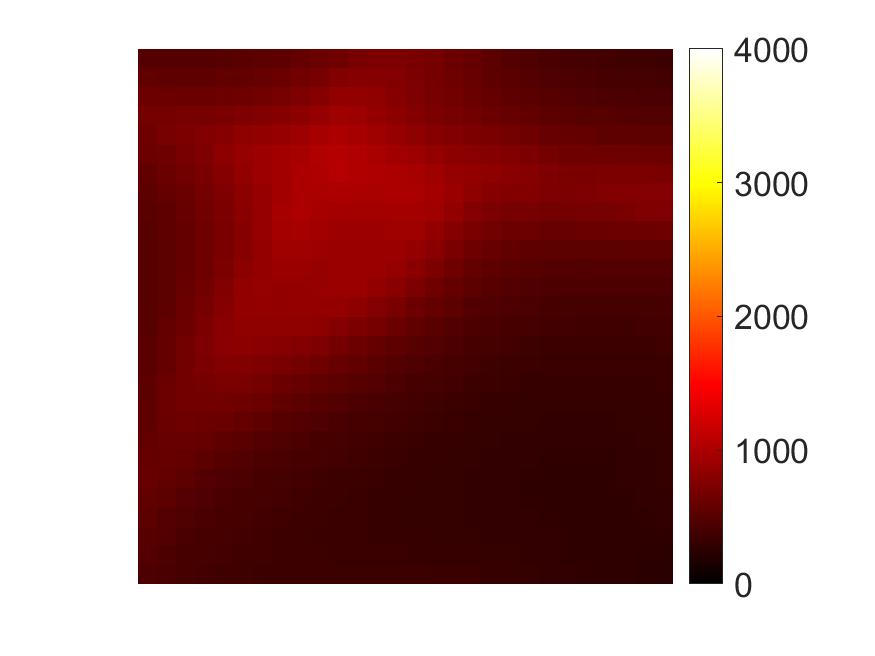}} & \hspace{-0.3cm}\adjustbox{valign=m,vspace=1pt}{\includegraphics[width=.25\textwidth]{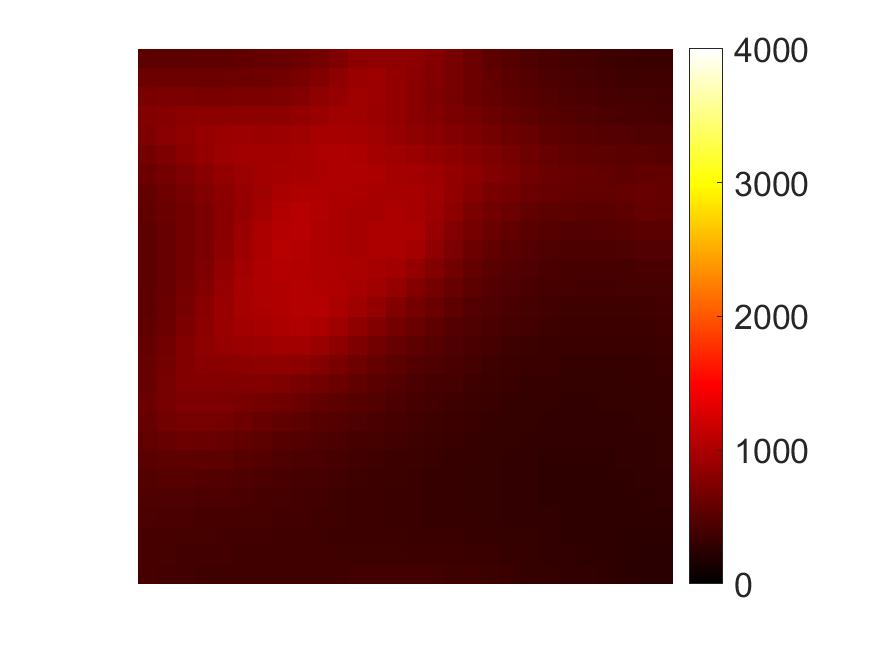}} & \hspace{-0.3cm}\adjustbox{valign=m,vspace=1pt}{\includegraphics[width=.25\textwidth]{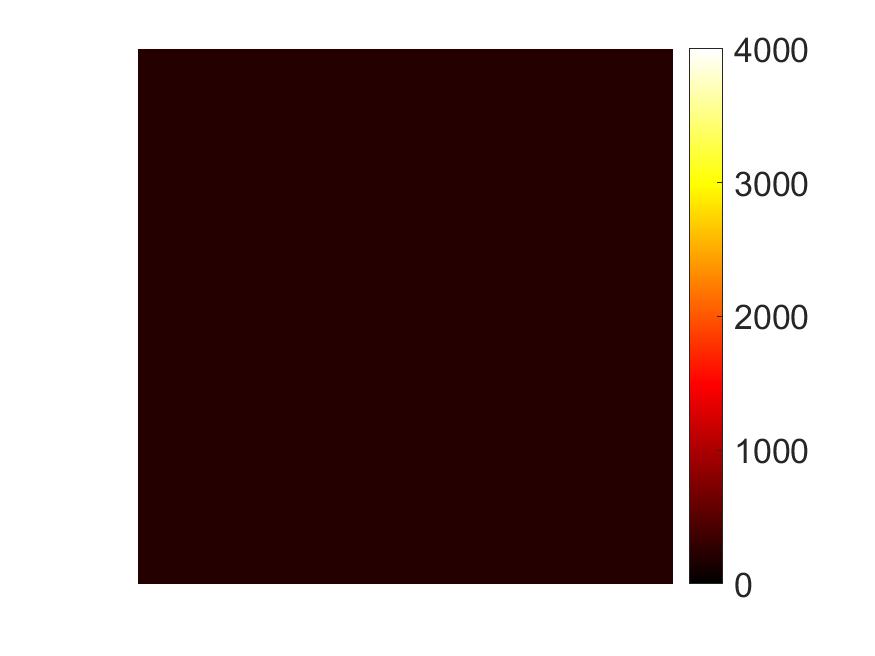}}\\
(b) & \adjustbox{valign=m,vspace=1pt}{\includegraphics[width=0.25\textwidth]{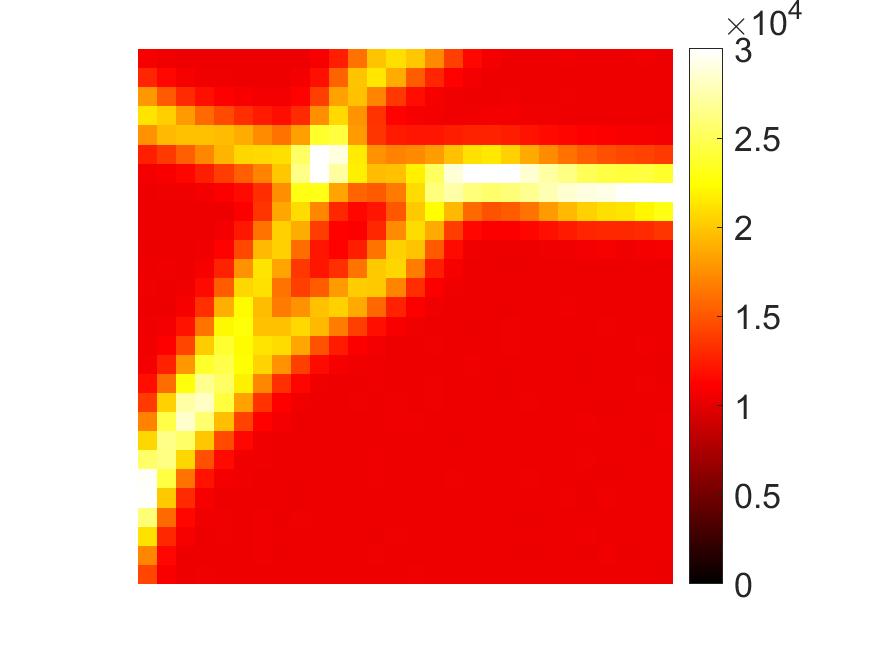}} & \hspace{-0.3cm}\adjustbox{valign=m,vspace=1pt}{\includegraphics[width=.25\textwidth]{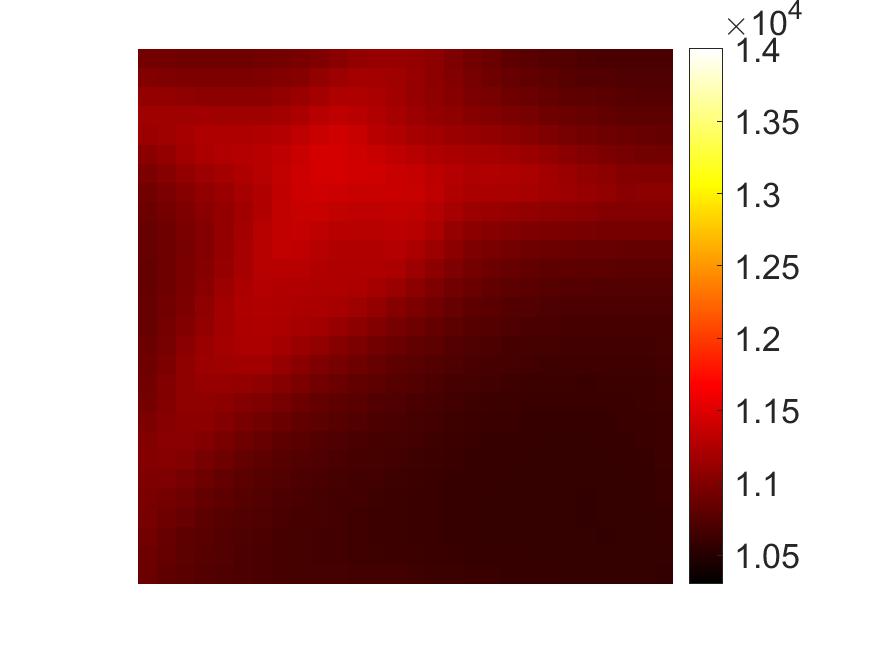}} & \hspace{-0.3cm}\adjustbox{valign=m,vspace=1pt}{\includegraphics[width=.25\textwidth]{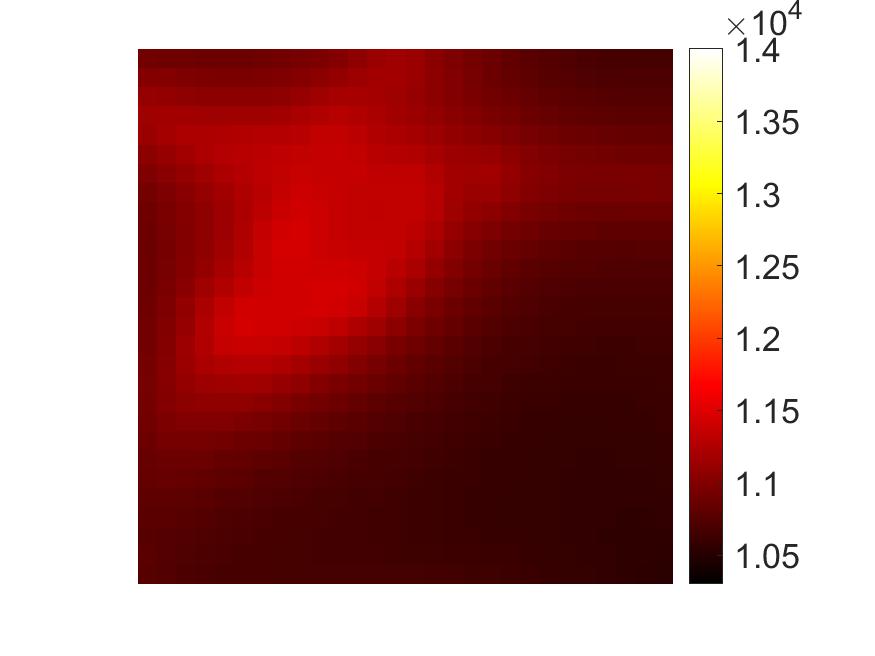}} & \hspace{-0.3cm}\adjustbox{valign=m,vspace=1pt}{\includegraphics[width=.25\textwidth]{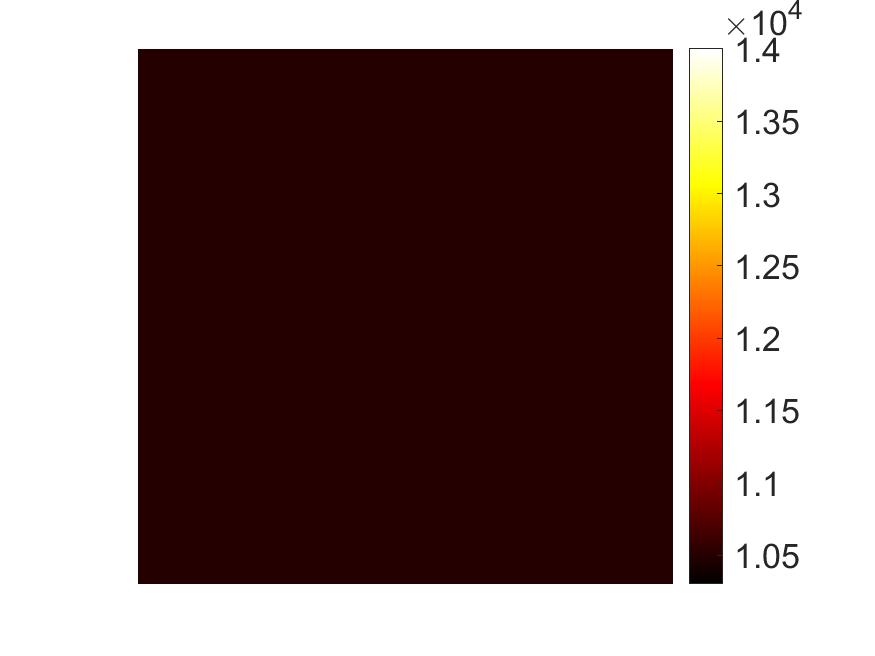}} \\
\end{tabular}
 \caption{(a) Low-background (LB) dataset: Diffraction limited image $\bar{\mathbf{y}}=\frac{1}{T}\sum_{t=1}^T \mathbf{y_t}$ with T=500 (4x zoom), Background estimation result on estimated support using CEL0 and $\ell_1$ regularization, Ground truth (GT) background image.
 (b)High-background (HB) dataset: Diffraction limited image $\bar{\mathbf{y}}=\frac{1}{T}\sum_{t=1}^T \mathbf{y_t}$ with T=500 (4x zoom), Background estimation result on estimated support using CEL0 and $\ell_1$ regularization, Ground truth (GT) background image. Please note the different scales between the diffraction limited and background images for a better visualization of the results}
 \label{fig:back}
 \end{figure}

\mdf{Background estimation results are available in Figure \ref{fig:back} where \eqref{eq:intensity_penalized} is used for intensity/background estimation on the supports $\Omega_{\cal{R}}$, with ${\cal{R}}=$ CEL0 and ${\cal{R}}=\ell_1$, that have been already estimated in the first step. In the figure there is also the constant background generated by the SOFI Simulation Tool\cite{SOFItool}, the software we used to generate our simulated data (more details in Section \ref{sec: Simulated Data}). Although the results look different due to the considered space-variant regularisation on $\mathbf{b}$, the variations are very little. The estimated background is smooth, as expected, while higher values are estimated near the simulated filaments and values closer to the true background are found away from them.}

\section{Automatic selection of regularization parameters}

We describe in this section two parameter selection strategies addressing the problem of estimating the regularization parameters $\lambda$ and $\mu$ appearing in the COL0RME support estimation problem \eqref{eq:support_mini} and intensity estimation one \eqref{eq:intensity_constraint}, respectively. The other two regularization parameters $\beta$ and $\alpha$ do not need fine tuning. They are both chosen arbitrary high, so as with large enough $\beta$ to allow for a very smooth background and with very high $\alpha$ to respect the required constraints (positivity for both intensity and background and restriction to the predefined support only for the intensity estimation). 

\subsection{Estimation of support regularization parameter $\lambda$}
\label{parameter_lambda}
%\txtb{Present the idea of restarting, Present the support results \\ 1. Result's presentation will look like the Figures 1 and 2 (they are not the correct results), Comments: make the images squared, interpolation of the Diffraction liminited image and superimpose green or blue crosses to present the GT support\\
%2. Compute the Jaccard Index or the FRC}\\

The selection of the regularization parameter value $\lambda$ in \eqref{eq:support_mini} is critical, as it determines the sparsity level of the support of the emitters. For its estimation, we start by computing a reference value $\lambda_{max}$, defined as the smallest regularization parameter for which the identically zero solution is found. It is indeed possible to compute such a $\lambda_{max}$ for both regularization terms CEL0 and $\ell_1$ (see \cite{soubiesPhD} and \cite{koulouri}). Once such values are known, we thus need to find a fraction $\gamma\in(0,1)$ of $\lambda_{max}$ corresponding to the choice $\lambda = \gamma \lambda_{max}$. For the CEL0 regularizer the expression for $\lambda_{max}$ (see Proposition 10.9 in \cite{soubiesPhD}) is:
\begin{equation}
    \lambda^{CEL0}_{max} := \max_{1 \leq i \leq L^2} \frac{\langle\mathbf{a}_i, \mathbf{r_y}\rangle^2}{2\|\mathbf{a}_i\|^2},
\end{equation}
where $\mathbf{a}_i = (\mathbf\Psi \odot \mathbf\Psi)_i$ denotes the $i$-th column of the operator $\mathbf{A}:=\mathbf\Psi \odot \mathbf\Psi$. Regarding the $\ell_1$-norm regularization penalty, $\lambda_{max}$ is given as follows:
\begin{equation}
    \lambda^{\ell_1}_{max} := \| \mathbf{A}^\intercal \mathbf{r_y}\|_\infty = \max_{1 \leq i \leq L^2} \langle\mathbf{a}_i, \mathbf{r_y}\rangle.
\end{equation}

%\txtb{I think that switching the description from $\lambda$ to $\gamma$ may be a bit confusing if someone looks only at this section. Can we leave the description as above but then talk only about $\lambda$, which is indeed the parameter appearing in the model?}

\begin{figure}[H]
     \centering
     \begin{subfigure}[b]{0.44\textwidth}
         \centering
         \includegraphics[width=\textwidth]{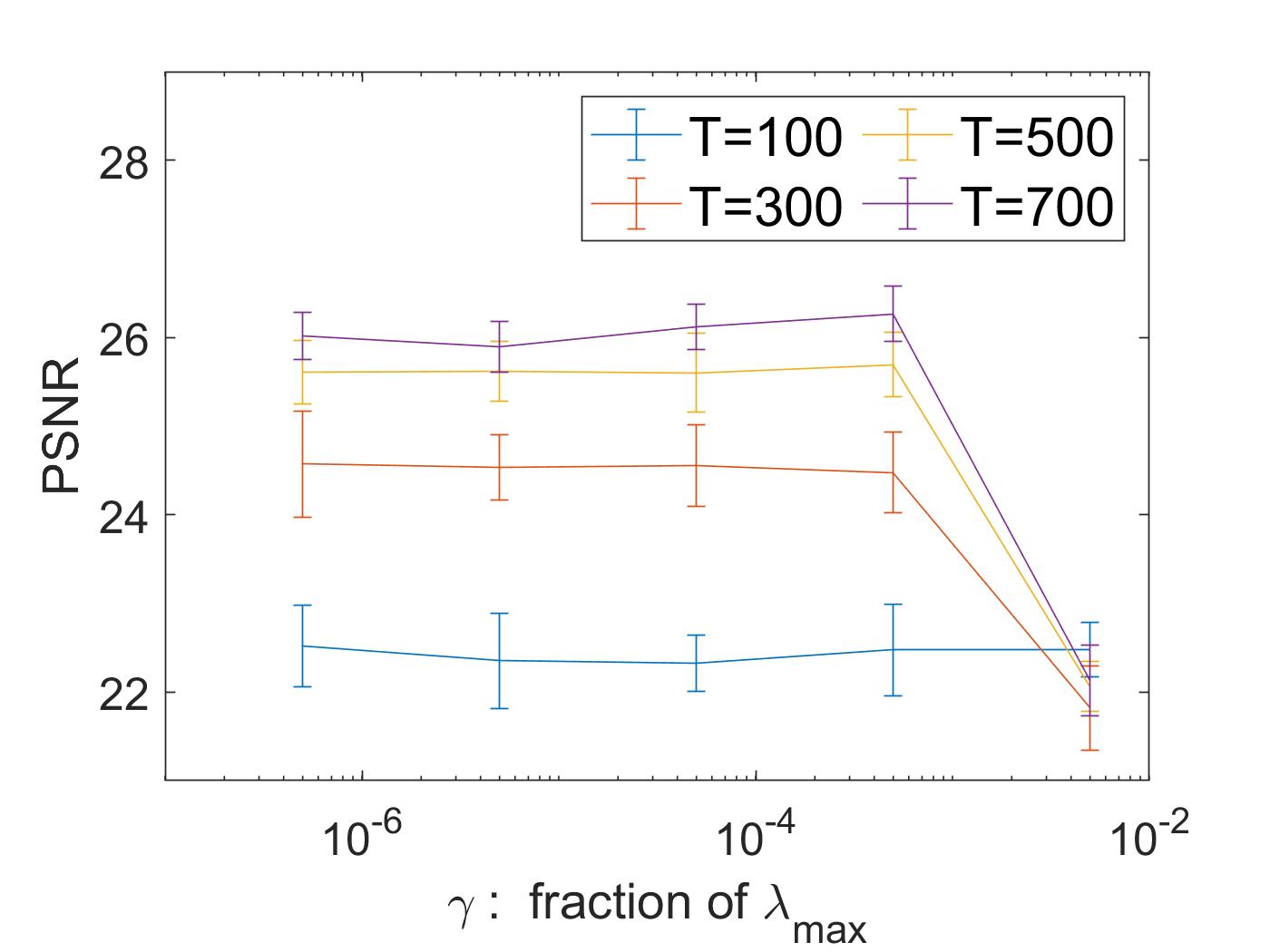}
         \caption{LB dataset}
     \end{subfigure}
     \hfill
     \begin{subfigure}[b]{0.44\textwidth}
         \centering
         \includegraphics[width=\textwidth]{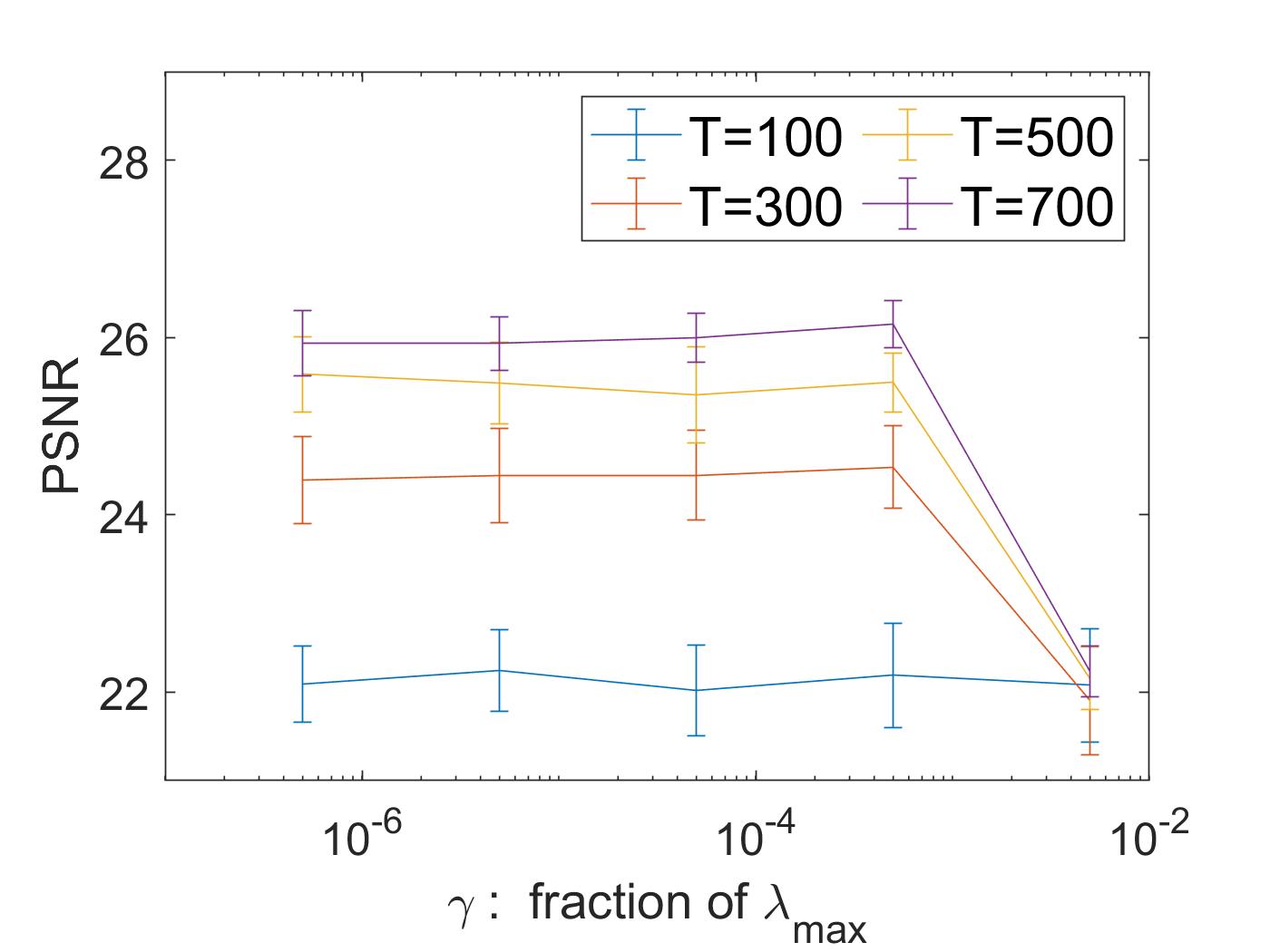}
         \caption{HB dataset}
     \end{subfigure}
     \caption{The PSNR value of the final COL0RME image, using the $\ell_1$-norm regularizer for support estimation, for different $\gamma$ values, evaluating in both the low-background (LB) and high-background (HB) dataset. The mean and the standard deviation of 20 different noise realization are presented}
     \label{fig:lambdas}
\end{figure}

As far as $\ell_1$ is used as regularization term in \eqref{eq:support_mini}, we report in Figure \ref{fig:lambdas} a graph showing how the PSNR value of the final estimated intensity image (i.e. after the application of the second step of COL0RME) varies for the two datasets considered depending on $\lambda$. It can be observed that for a large range of values $\lambda$, the final PSNR remains almost the same. Although this may look a bit surprising at a first sight, we remark that such a robust result is due, essentially, to the second step of the algorithm where false localizations related to an underestimation of $\lambda$ can be corrected through the intensity estimation step. Note, however, that in the case of an overestimation of $\lambda$, points contained in the original support are definitively lost so no benefit is obtained from the intensity estimation step, hence the overall PSNR decreases.

When the CEL0 penalty is used for support estimation, a heuristic parameter selection strategy can be used to improve the localization results but also to avoid the fine parameter tuning. More specifically, the non-convexity of the model can be used by considering an algorithmic restarting approach to improve the support reconstruction quality. In short, a value of $\lambda$ can be fixed, typically $\lambda = \gamma \lambda_{max}^{CEL0}$ with $\gamma \approx 5\times 10^{-4}$, so as to achieve a very sparse reconstruction. Then, the support estimation algorithm can be run and iteratively repeated with a new initialization (that is, restarted) several times. While keeping $\lambda$ fixed along this procedure, a wise choice of the initialization depending, but not being equal to the previous output can be used to enrich the support, see Appendix \ref{appendixc} for more details. Non-convexity is here exploited by changing, for a fixed $\lambda$, the initialization at each algorithmic restart, so that new local minimizers (corresponding to possible support points) can be computed. The final support image can thus be computed as the superposition of the different solutions computed at each restarting. In such a way, a good result for a not-finely-tuned value of $\lambda$ can be computed.

\subsection{Estimation of intensity regularization parameter $\mu$ by discrepancy principle}
\label{sec: DP}

In this section we provide some details on the estimation of the parameter $\mu$ in \eqref{eq:intensity_constraint}, which is crucial for an accurate intensity estimation. Recall that the problem we are looking at in this second step is
\begin{equation}
    \text{find}\quad\mathbf{x}\in\mathbb{R}^{L^2}\quad\text{s.t.}\quad\overline{\mathbf{y}} = \mathbf{\Psi x} + \mathbf{b} + \overline{\mathbf{n}},
    \label{eq: Model_intensity}
\end{equation}
where the quantities correspond to the temporal averages of the vectorized model in \eqref{eq:model_vec}, so that  $\overline{\mathbf{n}} = \frac{1}{T}\sum_{t=1}^T\mathbf{n_t}$. The temporal realizations $\mathbf{n_t}$ of the random vector $\mathbf{n}$ follow a normal distribution with zero mean and covariance matrix $s\mathbf{I}_{M^2}$, where $s$ has been estimated in the first step of the algorithm, see Section \ref{sec:noise_variance}. Consequently, the vector $\overline{\mathbf{n}}$ follows also a normal distribution with zero mean and  covariance matrix equal to $\frac{s}{T}\mathbf{I}_{M^2}$. As both $s$ and $T$ are known, we can use the discrepancy principle, a well-known a-posteriori parameter-choice strategy (see, e.g., \cite{DiscInvProb_Hansen,Gfrerer}), to efficiently estimate the hyper-parameter $\mu$. To detail how the procedure is applied to our problem, we write $\mathbf{x}_\mu$ in the following to highlight the dependence of $\mathbf{x}$ on $\mu$. According to the discrepancy principle strategy, the regularization parameter $\mu$ is chosen so that the residual norm of the regularized solution satisfies:
\begin{equation}
    \|\overline{\mathbf{y}} - \Psi \hat{\mathbf{x}}_\mu - \hat{\mathbf{b}}\|_2^2 = \nu_{DP}^2\|\overline{\mathbf{n}}\|_2^2,
    \label{DiscrPrinc}
\end{equation}
where $\hat{\mathbf{x}}_\mu \in \mathbb{R}^{L^2}$ and $\hat{\mathbf{b}} \in \mathbb{R}^{M^2}$ are the  solutions of \eqref{eq:intensity_constraint}. The expected value of $\|\overline{\mathbf{n}}\|_2^2$ is:
\begin{equation}
    \mathop{{}\mathbb{E}}\{\|\overline{\mathbf{n}}\|_2^2\} =M^2 \frac{s}{T},
\end{equation}
which can be used as an approximation of $\|\overline{\mathbf{n}}\|_2^2$ for $M^2$ big enough.
%something that happens in our case, we can approximate $\|\overline{\mathbf{n}}\|_2^2 \approx M^2 \frac{s}{K} $. 
The scalar value $\nu_{DP} \approx 1$ is a 'safety factor' that plays an important role in the case when a good estimate of $\|\overline{\mathbf{n}}\|_2$ is not available. In such situations a value $\nu_{DP}$ closer to $2$ is used. As detailed in Section \ref{sec:noise_variance}, the estimation of $s$ is rather precise in this case, hence we fix $\nu_{DP} = 1$ in the following.

%As we have already computed in the first step of the method with a great precision the variance of the noise, we can efficiently estimate the hyper-parameter $\mu$ using the discrepancy principle, a well-known parameter-choice strategy. Knowing the variance of the noise, we can compute 
%Due to the steep part of the L-curve we see that 
%( An underestimate of $\|\mathbf{e}\|_2$ is likely to produce an underregularized solution 
%with a very large (semi)norm 
%and on the other hand, an overestimate of $\|\mathbf{e}\|_2$ produces an overregularized solution with too large regularization error. So, an adaptation of the "safety factor" in certain cases is crucial. )
% When a good estimate for $\|\mathbf{e}\|_2$ is known, as it happens in our case, this method yields a good regularization parameter, corresponding to a regularized solution immediately to the right of the L-curve’s corner.

We can now define the function $f(\mu): \mathbb{R}_+ \rightarrow \mathbb{R}$ as:  
\begin{equation}
    f(\mu) = \frac12\|\overline{\mathbf{y}} - \mathbf\Psi \hat{\mathbf{x}}_\mu - \hat{\mathbf{b}}\|_2^2 - \frac{\nu_{DP}^2}{2}\|\overline{\mathbf{n}}\|_2^2.
    \label{f}
\end{equation}
We want to find the value $\hat{\mu}$ such that $f(\hat{\mu}) = 0$. This can be done iteratively, using the Newton's method whose iterations read:
\begin{equation}
    \mu_{n+1} = \mu_n - \frac{f(\mu_n)}{f'(\mu_n)}, \quad n=1,2,.. .
\end{equation}

In order to be able to compute easily the values  $f(\mu)$ and $f '(\mu)$, the values $\hat{\mathbf{x}}_\mu \in \mathbb{R}^{L^2}$, $\hat{\mathbf{b}} \in \mathbb{R}^{M^2}$ and $\hat{\mathbf{x}}'_\mu = \frac{\partial}{\partial \mu} \hat{\mathbf{x}}_\mu  \in \mathbb{R}^{L^2}$ need to be computed, as it can be easily noticed by writing the expression of $f'(\mu)$ which reads:
\begin{align}
    f'(\mu) &= \frac{\partial}{\partial \mu}\{\frac12\|\overline{\mathbf{y}} - \mathbf\Psi \hat{\mathbf{x}}_\mu - \hat{\mathbf{b}}\|_2^2\} = (\hat{\mathbf{x}}'_\mu)^\intercal \mathbf\Psi^\intercal(\overline{\mathbf{y}} - \mathbf\Psi \hat{\mathbf{x}}_\mu - \hat{\mathbf{b}}).
    \label{f_til}
\end{align}

The values $\hat{\mathbf{x}}_\mu$ and $\hat{\mathbf{b}}$ can be found by solving the minimization problem \eqref{eq:intensity_constraint}. As far as $\hat{\mathbf{x}}'_\mu$ is concerned, we report in Appendix \ref{appendixB} the steps necessary for its computation. \mdffirst{We note here, however, that in order to compute such a quantity, the relaxation of the support/non-negativity constraints by means of the smooth quadratic terms discussed above is fundamental.} One can show that $\hat{\mathbf{x}}'_\mu$ is the solution of the following minimization problem:
\begin{equation}\label{min x'}
     \hat{\mathbf{x}}'_\mu = \argmin_{\mathbf{x}\in \mathbb{R}^{L^2}} \frac12\|\mathbf{\Psi x} \|_2^2 + \frac{\mu}{2} \|\nabla \mathbf{x} + \mathbf{c}\|_2^2+ \frac{\alpha}{2} \left( \|\mathbf{I_\Omega x}\|_2^2 + \| \mathbf{I}_{\mathbf{\hat{x}}_\mu}\mathbf{x}\|^2_2\right),
\end{equation}
where $\mathbf{c}$ is a known quantity defined by $\mathbf{c} = \frac{1}{\mu} \nabla \hat{\mathbf{x}}_\mu$, and the diagonal matrix $\mathbf{I}_{\mathbf{\hat{x}}_\mu} \in \mathbb{R}^{L^2 \times L^2}$ identifies the support of $\mathbf{\hat{x}_\mu}$ by:
\[
\mathbf{I}_{\mathbf{\hat{x}}_\mu}(i,i) = \begin{cases}
 0 & \text{if ${(\hat{\mathbf{x}}_\mu)_i} \geq 0$},\\
          1 & \text{if ${(\hat{\mathbf{x}}_\mu)_i} < 0$}  .
\end{cases}
\]

\begin{figure}
    \centering
    \includegraphics[width=0.8\linewidth]{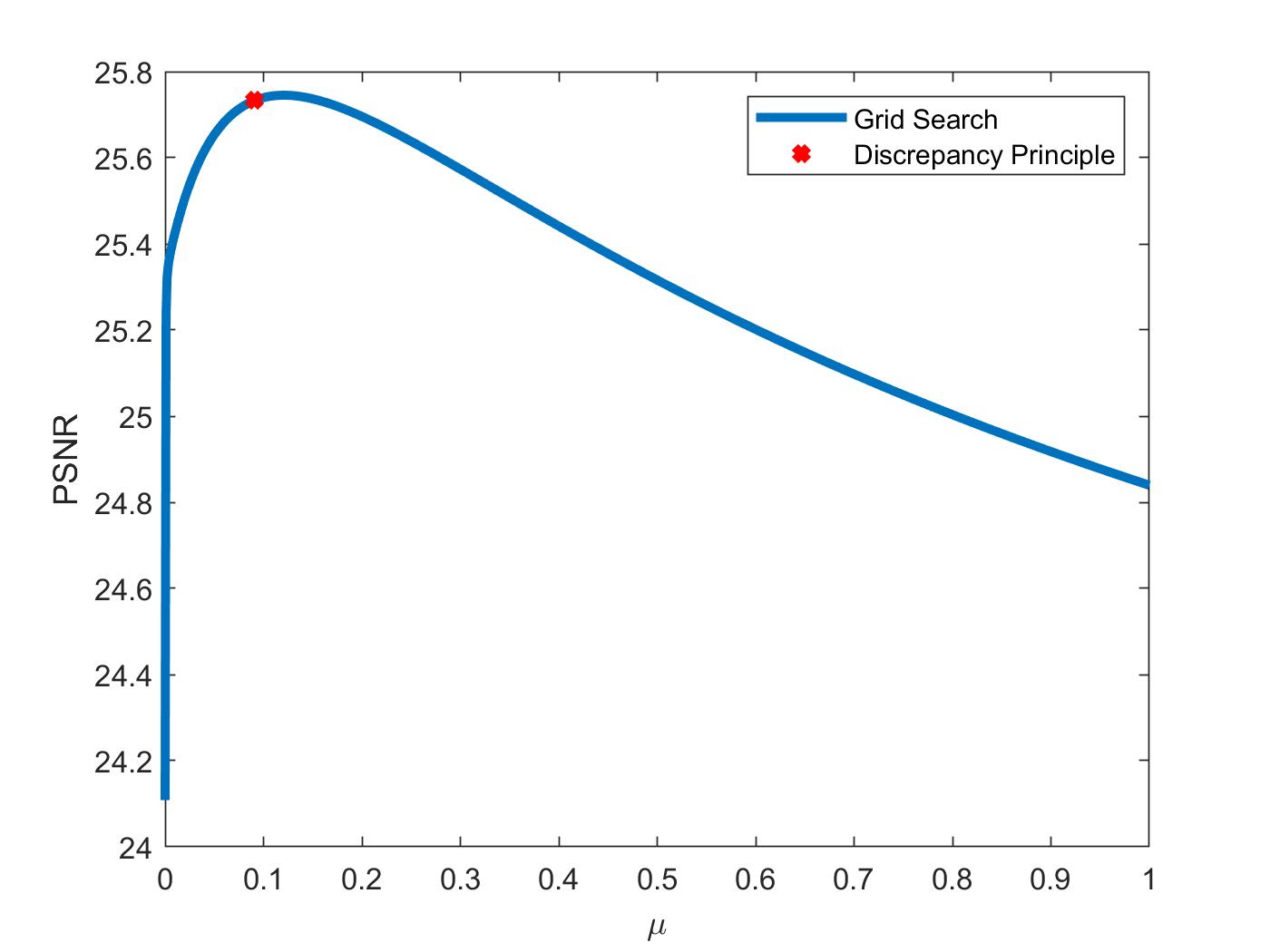}
    \caption{The solid blue line shows the PSNR values computed by solving \eqref{eq:intensity_penalized} for several values of $\mu$ within a specific range. Tha data used are the HB dataset with  $T=500$ frames (Figure \ref{subfig:HB}) and the $\ell_1$-norm regularization penalty. The red cross shows the PSNR value $\hat{\mu}$ obtained by applying the Discrepancy Principle. We note that such value is very close to one maximizing the PSNR metric}
    \label{fig:grid_search_DP}
\end{figure}

We can find $\hat{\mathbf{x}}'_\mu$ by iterating
\begin{equation}
    {\mathbf{x}'}_\mu^{n+1} = \textbf{\text{prox}}_{\overline{h},\tau}({\mathbf{x}'}_\mu^{n} - \tau \nabla \overline{g}({\mathbf{x}'}_\mu^{n})), \quad n=1,2,.. ,
    \label{x_m'}
\end{equation}
where 
\begin{equation}
    \overline{g}(\mathbf{x}) := \frac12\|\mathbf{\Psi x} \|_2^2 + \frac{\mu}{2} \|\nabla \mathbf{x} + \mathbf{c}\|_2^2, \qquad
    \overline{h}(\mathbf{x}) := \frac{\alpha}{2} \left(\|\mathbf{I_\Omega x}\|_2^2 + \| \mathbf{I}_{\mathbf{\hat{x}}_\mu}\mathbf{x}\|^2_2\right).
    \label{eq:overline_h}
\end{equation}
For $\mathbf{z}\in\mathbb{R}^{L^2}$, the proximal operator $\textbf{\text{prox}}_{\overline{h}, \tau}(\mathbf{z})$ can be obtained following the computations in Appendix \ref{appendixA}:
 \begin{equation}
         (\textbf{\text{prox}}_{\overline{h}, \tau}(\mathbf{z}))_i= {\text{prox}}_{\overline{h}, \tau}(\mathbf{z}_i) = \frac{\mathbf{z}_i}{1 + \alpha \tau \left( \mathbf{I_\Omega}(i,i) + \mathbf{I}_{\hat{\mathbf{x}}_\mu}(i,i) \right)},
    \end{equation}
while 
\begin{equation}
    \nabla \overline{g}(\mathbf{x}') =  (\mathbf\Psi^\intercal \mathbf\Psi +\mu \nabla^\intercal \nabla) \mathbf{x}' +  \nabla^\intercal \nabla \hat{\mathbf{x}}_\mu,
\end{equation}
and the step $\tau\in (0, \frac{1}{L_{\overline{g}}}]$, with $L_{\overline{g}} = \| \mathbf\Psi^\intercal \mathbf\Psi +\mu \nabla^\intercal \nabla\|_2$ the Lipschitz constant of $\nabla \overline{g}$. A pseudo-code explaining the procedure we follow to find the optimal $\hat{\mu}$ can be found in Algorithm \ref{Algorithm:NewtMeth}. Finally, in Figure \ref{fig:grid_search_DP}, a numerical example is available to show the good estimation of the parameter $\hat{\mu}$.

\begin{algorithm}[!h]
\caption{Discrepancy Principle}
\label{Algorithm:NewtMeth}
\begin{algorithmic}
\REQUIRE $\overline{\mathbf{y}}\in\mathbb{R}^{M^2}, \mathbf{x}^0\in\mathbb{R}^{L^2},\mathbf{b}^0\in\mathbb{R}^{M^2},{\mu}_0, \beta>0$, $\alpha\gg 1$
\REPEAT 
\STATE Find $\hat{\mathbf{x}}_{\mu_n}, \hat{\mathbf{b}}$ \qquad\qquad\qquad using Algorithm \ref{Algorithm:AMA_intensity}
\STATE Find ${\hat{\mathbf{x}}'}_{\mu_n}$ \qquad\qquad\qquad\quad solving \eqref{min x'}
\STATE Compute $f(\mu_n),f'(\mu_n)$\quad from \eqref{f} and \eqref{f_til}
\STATE $\mu_{n+1} \gets \mu_n - \frac{f(\mu_n)}{f'(\mu_n)}$
\UNTIL convergence
\RETURN $\hat{\mu}$
\end{algorithmic}
\end{algorithm}

% {\large{The algorithm:}}
% \begin{algorithmic}
% \State $\mu_0 \gets 0$
% \For{$n=0:n_{max}-1$ or till convergence}     
%     \State compute $\mathbf{x_{\mu_n}}$  \qquad\qquad\qquad\% eq. \ref{x_m}
%     \State compute $\mathbf{x'_{\mu_n}}$ \qquad\qquad\qquad\% eq. \ref{x_m'}
%     \State compute
%     $f(\mathbf{x}_\mu),f'(\mathbf{x}_\mu)$ \qquad\% eq. \ref{f}, eq. \ref{f_til}
%     \State $\mu_{n+1} \gets \mu_n - \frac{f(\mathbf{x}_\mu)}{f'(\mathbf{x}_\mu)}$
% \EndFor\\
% \end{algorithmic}

\section{Results}
\label{sec:Results}
In this section we compare the method COL0RME with state-of-the-art methods that exploit the temporal fluctuations/blinking of fluorophores, while applying them to simulated and real data. More precisely we compare: COL0RME-CEL0 (using the CEL0 regularization in the support estimation), COL0RME-$\ell_1$ (using the $\ell_1$-norm regularization in the support estimation), SRRF\cite{srrf}, SPARCOM\cite{SPARCOM} and LSPARCOM\cite{LSPARCOM}. 
%\mdfsec{We further performed preliminary comparisons also with the ESI, 3B and bSOFI approaches, but as they were unsatisfying 
%or very hard to tune in order to get meaningful results, we omit them in the following.}
\mdfsec{We further performed preliminary comparisons also with the ESI, 3B and bSOFI approaches using available codes provided by the authors on the web\footnote{ESI: \href{https://github.com/biophotonics-bielefeld/ESI}{https://github.com/biophotonics-bielefeld/ESI}, 3B: \href{http://www.coxphysics.com/3b}{http://www.coxphysics.com/3b}, bSOFI implemented in SOFI Simulation Tool software package: \href{https://github.com/lob-epfl/sofitool}{https://github.com/lob-epfl/sofitool}}, but \mdfsecrev{we did not successfully} obtain satisfactory results, so we omit them in the following.}
% Laure : We further performed preliminary comparisons also with the ESI, 3B and bSOFI approaches using available codes provided by the authors on the web (mettre les liens en footnotes),  but we do not succeed to obtain satisfactory results, so we omit them in the following.

\subsection{Simulated Data}
\label{sec: Simulated Data}
\begin{figure}[h!]
     \centering
     \begin{subfigure}[b]{0.24\textwidth}
         \centering
         \includegraphics[width=1.1\textwidth]{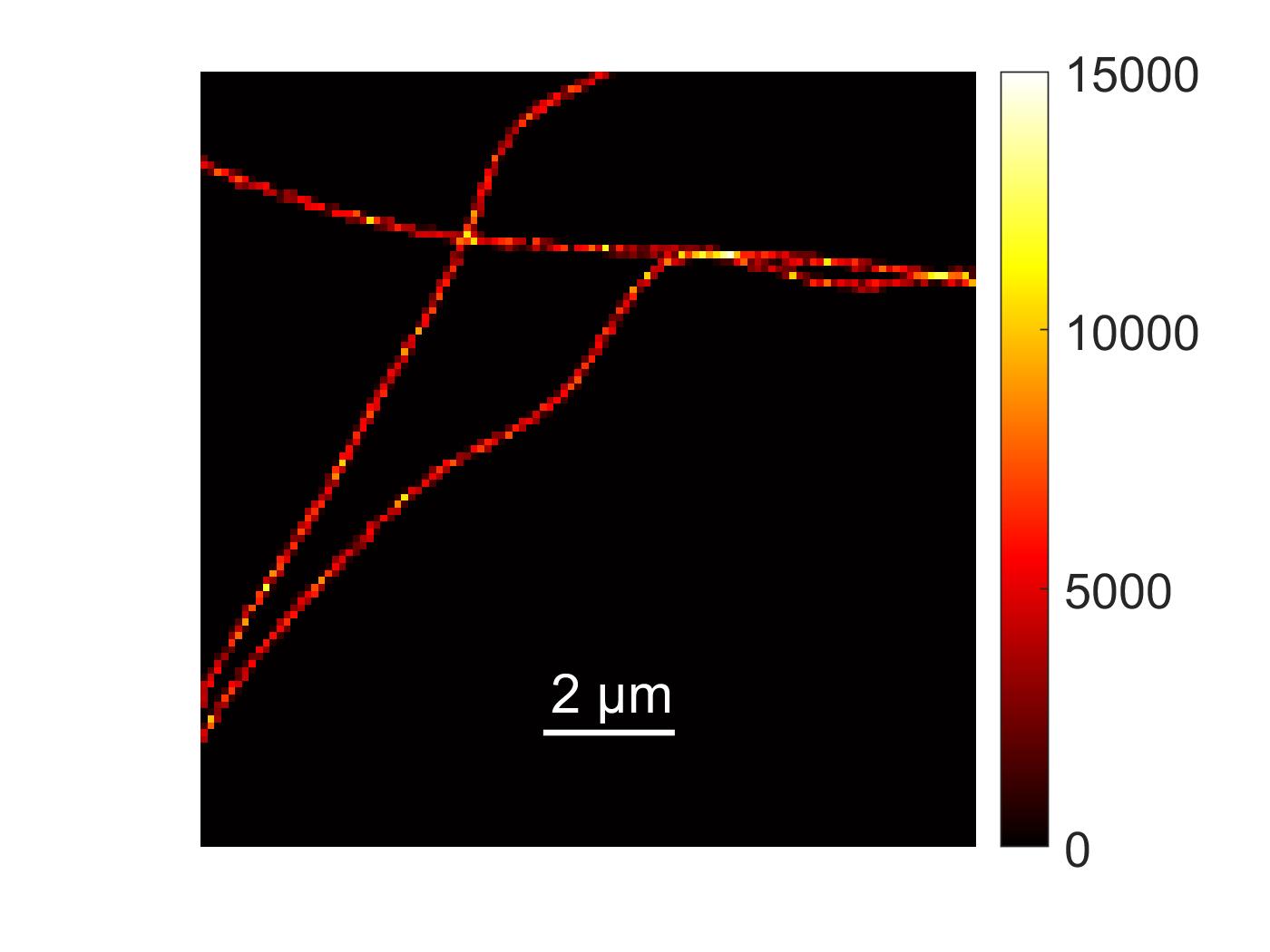}
         \caption{GT image}
         \label{GT}
     \end{subfigure}
    %  \hfill
    %  \begin{subfigure}[b]{0.24\textwidth}
    %      \centering
    %      \includegraphics[width=1.1\textwidth]{Figures/DL_noBg14_SNR15_K500_N34.jpg}
    %      \caption{$\bar{\mathbf{y}}$ (NB)}
    %  \end{subfigure}
    %  \hfill
     \begin{subfigure}[b]{0.24\textwidth}
         \centering
         \includegraphics[width=1.1\textwidth]{Figures/DLwithCB_lowBg14_SNR15_K500_N34.jpg}
         \caption{$\bar{\mathbf{y}}$ (LB)}
     \end{subfigure}
    %   \hfill
     \begin{subfigure}[b]{0.24\textwidth}
         \centering
         \includegraphics[width=1.1\textwidth]{Figures/DLwithCB_highBg14_SNR15_K500_N34.jpg}
         \caption{$\bar{\mathbf{y}}$ (HB)}
         \label{subfig:HB}
     \end{subfigure}
        \caption{ The Ground truth (GT) intensity image, as well as, the diffraction limited images $\bar{\mathbf{y}}=\frac{1}{T}\sum_{t=1}^T \mathbf{y_t}$ for the two datasets with a 4x zoom, for a sequence of T=500 frames}
        \label{fig:datasets}
\end{figure}

To evaluate the method COL0RME we choose images of tubular structures that simulate standard microscope acquisitions with standard fluorescent dyes. In particular, the spatial pattern (see Figure \ref{GT}) is taken from the MT0 microtubules training dataset uploaded for the SMLM Challenge of 2016\footnote{\href{http://bigwww.epfl.ch/smlm/datasets/index.html}{http://bigwww.epfl.ch/smlm/datasets/index.html}}. The temporal fluctuations are obtained by using the SOFI Simulation Tool \cite{SOFItool}. This simulation software, implemented in \textsc{Matlab}, generates realistic stacks of images, similar to the ones obtained from real microscopes, as it makes use of parameters of the microscope setup and some of the sample's main properties. \mdfsec{However, differently from the fluctuating\footnote{the emission of a single fluorophore over time can be described by a Poisson distribution} microscopic data presented in section \ref{sec: Real Data}, the blinking generated by the SOFI Simulation Tool have a more distinctive "on-off" behaviour.}

For the experiments presented in this paper, we generate initially a video of $700$ frames, however we evaluate the methods using the first $T=100$, $T=300$, $T=500$ and $T=700$ frames, so as to examine further the trade-off between temporal and spatial resolution. The frame rate is fixed at 100 frames per second (fps) and the pixel size is $100$ nm. Regarding the optical parameters, we set the numerical aperture equal to 1.4 and the emission \mdfsecrev{wavelength} to 525 nm, while the FWHM of the PSF is equal to $228.75$nm. The fluorophore parameters are set as follows: $20$ms for on-state average lifetime, $40$ms for off-state average lifetime and $20$s for average time until bleaching. The emitter density is equal to 10.7 emitters/pixel/frame, while 500 photons are emitted, on average, by a single fluorescent molecule in every frame. 

\begin{figure}
    \centering
     \begin{subfigure}[b]{0.35\textwidth}
         \centering
         \includegraphics[width=0.9\textwidth]{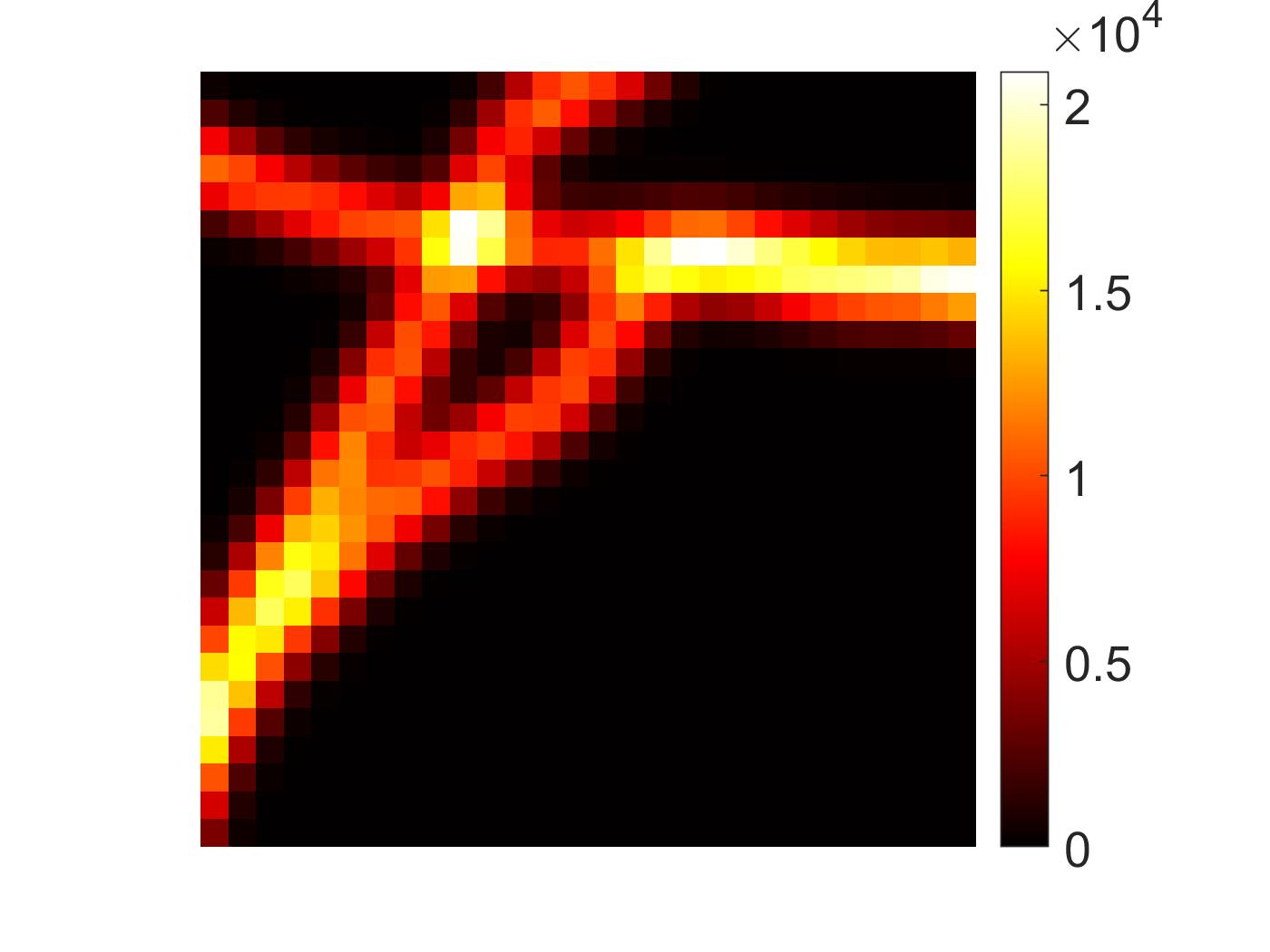}
         \caption{}
         \label{before_noise}
     \end{subfigure}
     \begin{subfigure}[b]{0.35\textwidth}
         \centering
         \includegraphics[width=0.9\textwidth]{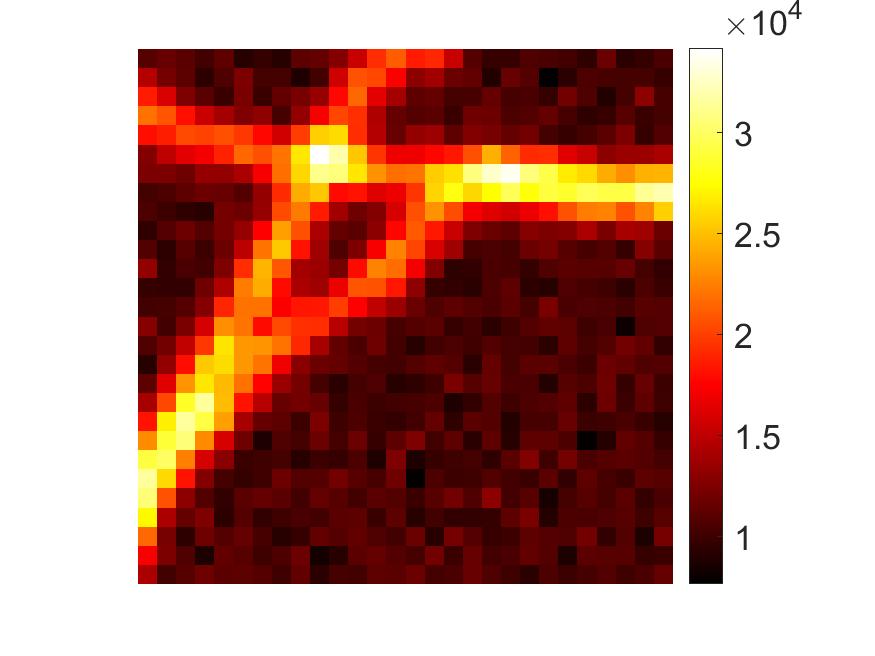}
         \caption{}
         \label{after_noise}
     \end{subfigure}
    \caption{One frame of the HB dataset, before and after the addition of background and the simulated noise degradation. (a) A convoluted and down-sampled image $\mathbf{\Psi} \mathbf{x}^{GT}_t$ obtained from a ground truth frame $\mathbf{x}^{GT}_t$, (b) A frame of the final noisy sequence: $\mathbf{y}_t$. Note the different colormaps to better capture the presence of noise and background}
    \label{fig:back_noise}
\end{figure}

We create two datasets with the main difference between them being the background level, as in real scenarios the background is usually present. More precisely we create: the low-Background (LB) dataset, where the background is equal to $50$ photons/pixel/frame and, the most realistic of the two, the high-Background (HB) dataset, where the background is equal to $2500$ photons/pixel/frame.  \mdfsecrev{In both datasets, we proceed as follows: initially, Poisson noise is added to simulate the photon noise (see (2)); subsequently, the number of photons recorded by each camera pixel is converted into an electric charge in accordance with the quantum efficiency and gain of the camera that have been set to 0.7 and 6 respectively (thus resulting in an overall gain of 4.2); finally, Gaussian noise is added.} In order to give a visual inspection of the background and noise, in Figure \ref{fig:back_noise}, one frame of the HB dataset is presented before and after the background/noise addition. As we want, also, to provide a quantitative assessment, we measure the quality of the reconstruction of the final sequence of $T$ frames ($\mathbf{y}_t, t=1,2,\dots,T$) using the Signal-to-Noise-Ration (SNR) metric, given by the following formula:
 \begin{equation}
     \text{SNR}_\text{dB} = 10 \log_{10} \left( \frac{\frac{1}{TM^2}\sum\limits_{i=1}^{TM^2}\left(\mathbf{R}_i\right)^2}{\frac{1}{TM^2}\sum\limits_{i=1}^{TM^2}\left(\mathbf{R}_i - \mathbf{K}_i\right)^2}\right),
 \end{equation}
where $\mathbf{R}\in \mathbb{R}^{TM^2}$ is the reference image and $\mathbf{K}\in \mathbb{R}^{TM^2}$ the image we want to evaluate, both of them in a vectorized form. As reference, we choose the sequence of convoluted and down-sampled ground truth frames (see one frame of the reference sequence in Figure \ref{before_noise}). The SNR values for a sequence of $T=500$ frames for the LB and HB dataset are $15.57$dB and $-6.07$dB, respectively. A negative value is computed for the HB dataset due to the very high background used in this case.

 The diffraction limited image (the average image of the stack) of each dataset as well as the ground truth intensity image are available in Figure \ref{fig:datasets}. In the LB dataset, due to the high signal values, the background is not visible. Further, as the observed microscopic images and the reconstructed ones belong to different grids, coarse and fine grid respectively, their intensity values are not comparable and we can not use the same colorbar to represent them. The intensity of one pixel in the coarse grid is the summation of the intensities of $ q \times q $ pixels in the fine grid, where $q$ is the super-resolution factor. For this reason, we use two different colorbars. %presented in Figure \ref{fig:colorbars}.

% \begin{figure}
%     \centering
%     \hspace{0.1cm}
%       \begin{subfigure}[b]{0.45\textwidth}
%          \centering
%          \includegraphics[width=\textwidth]{Figures/colorbarDL}
%          \caption{Colorbar for the coarse-grid representations }
%      \end{subfigure}
%      \hfill
%      \begin{subfigure}[b]{0.43\textwidth}
%          \centering
%          \includegraphics[width=\textwidth]{Figures/colorbarGT}
%           \vspace{0.05cm}
%          \caption{Colorbar for the fine-grid representations}
%      \end{subfigure}  
%      \caption{Colorbars used for the representation of the images}
%      \label{fig:colorbars}
% \end{figure}

The comparison of the method COL0RME with other state-of-the-art methods that take advantage of the blinking fluorophores is available bellow. Regarding the method COL0RME-CEL0 and COL0RME-$\ell_1$, a regularization parameter equal to $\lambda = 5 \times 10^{-4} \times \lambda_{max}^{CEL0}$ and $\lambda = 5 \times 10^{-4} \times \lambda_{max}^{\ell_1}$, respectively, is used in the support estimation. The hyper-parameters $\alpha$ and $\beta$ are set as follows: $\alpha = 10^6$, $\beta = 20$. For the method COL0RME-CEL0 the algorithmic restarting approach is used for a better support estimation. \mdf{It stops when there are not additional pixels added to the estimated support or if a maximum number of $10$ restarts is reached. Such number was empirically determined by preliminary simulations}. For the method SRRF we are using the NanoJ SRRF plugin for ImageJ\footnote{\href{https://github.com/HenriquesLab/NanoJ-SRRF}{https://github.com/HenriquesLab/NanoJ-SRRF}}. Concerning the method SPARCOM, we make use of the \textsc{Matlab} code available online\footnote{\href{https://github.com/KrakenLeaf/SPARCOM}{https://github.com/KrakenLeaf/SPARCOM}}. As regularization penalty we choose the $\ell_1$-norm with a regularization parameter equal to $10^{-10}$ and we avoid the post-processing step (the convolution with a small Gaussian function) for most precise localization. Finally we test the method LSPARCOM, using the code that is available online\footnote{\href{https://github.com/gilidar/LSPARCOM}{https://github.com/gilidar/LSPARCOM}} and the tubulin (TU) training set that is provided. 
 
\begin{figure}[h!]
\centering
\setlength\tabcolsep{1.5pt}
\renewcommand{\arraystretch}{1.5}
\begin{tabular}{ccc}
GT image & \small{COL0RME-CEL0} & \small{COL0RME-$\ell_1$}\\
\adjustbox{valign=m,vspace=1pt}{\begin{tikzpicture}[spy using outlines={rectangle,magnification=1.7,size=1.4cm,connect spies}]
\node {\pgfimage[width=4cm]{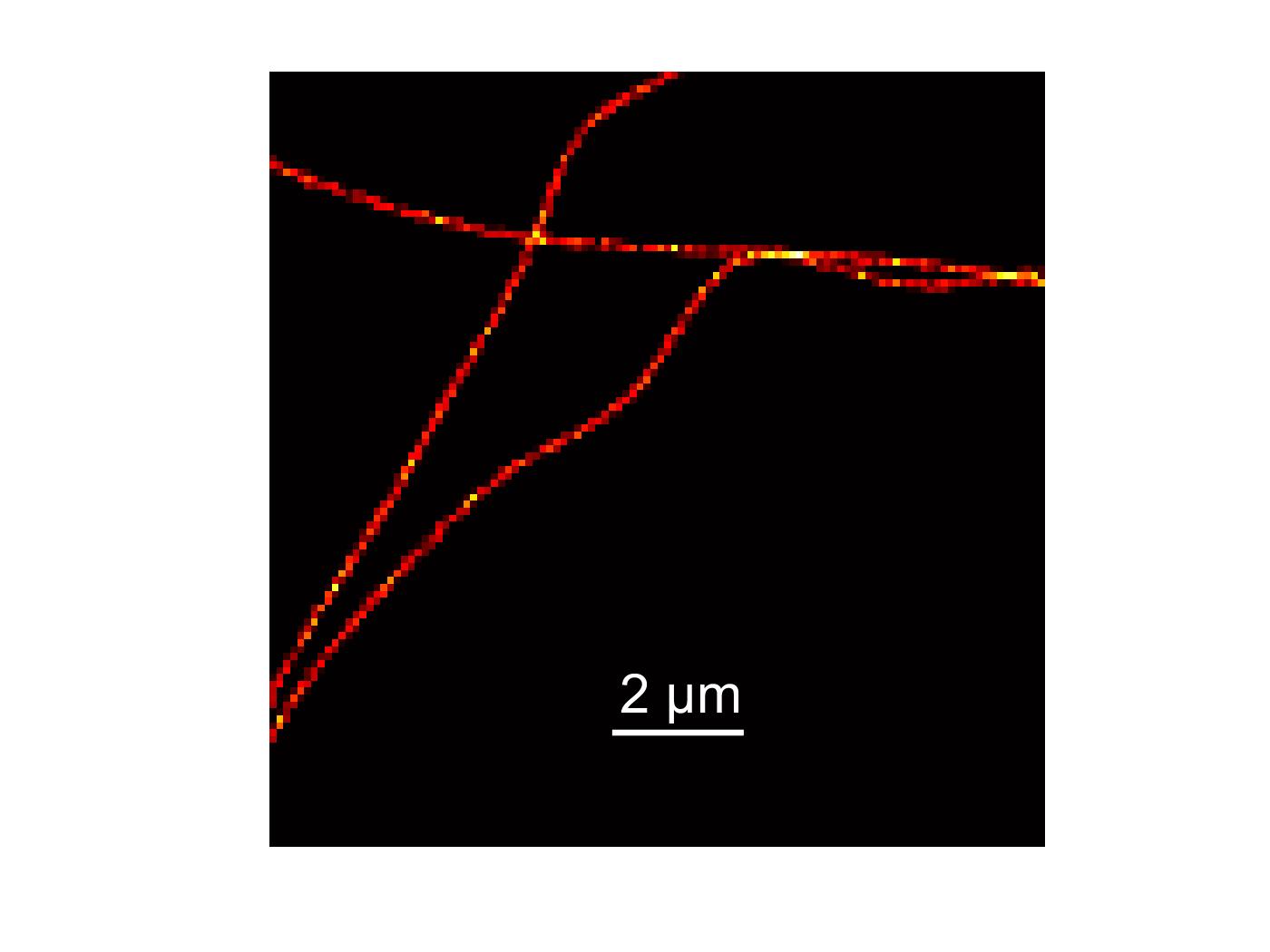}};
\spy[green] on (+0.82,0.6) in node [left] at (1.53,-1.7);
\spy[yellow] on (-0.7,-0.4) in node [left] at (0,-1.7);
\end{tikzpicture}}
& \adjustbox{valign=m,vspace=1pt}{\begin{tikzpicture}[spy using outlines={rectangle,magnification=1.7,size=1.4cm,connect spies}]
\node {\pgfimage[width=4cm]{Figures/COL0RME_int_lowBg14_SNR15_K500_N34_CEL0_0_RP2.jpg}};
\spy[green] on (+0.82,0.6) in node [left] at (1.53,-1.7);
\spy[yellow] on (-0.7,-0.4) in node [left] at (0,-1.7);
\end{tikzpicture}}
& \adjustbox{valign=m,vspace=1pt}{\begin{tikzpicture}[spy using outlines={rectangle,magnification=1.7,size=1.4cm,connect spies}]
\node {\pgfimage[width=4cm]{Figures/COL0RME_int_lowBg14_SNR15_K500_N34_L1_RP2.jpg}};
\spy[green] on (+0.82,0.6) in node [left] at (1.53,-1.7);
\spy[yellow] on (-0.7,-0.4) in node [left] at (0,-1.7);
\end{tikzpicture}}\\
\adjustbox{valign=m,vspace=1pt}{\includegraphics[width=.28\textwidth]{Figures/colorbarGT-eps-converted-to}}
& \adjustbox{valign=m,vspace=1pt}{\includegraphics[width=.28\textwidth]{Figures/colorbarGT-eps-converted-to}}
& \adjustbox{valign=m,vspace=1pt}{\includegraphics[width=.28\textwidth]{Figures/colorbarGT-eps-converted-to}}\\
SRRF & SPARCOM & LSPARCOM\\
\adjustbox{valign=m,vspace=1pt}{\begin{tikzpicture}[spy using outlines={rectangle,green,magnification=1.7,size=1.4cm,connect spies}]
\node {\pgfimage[width=4cm]{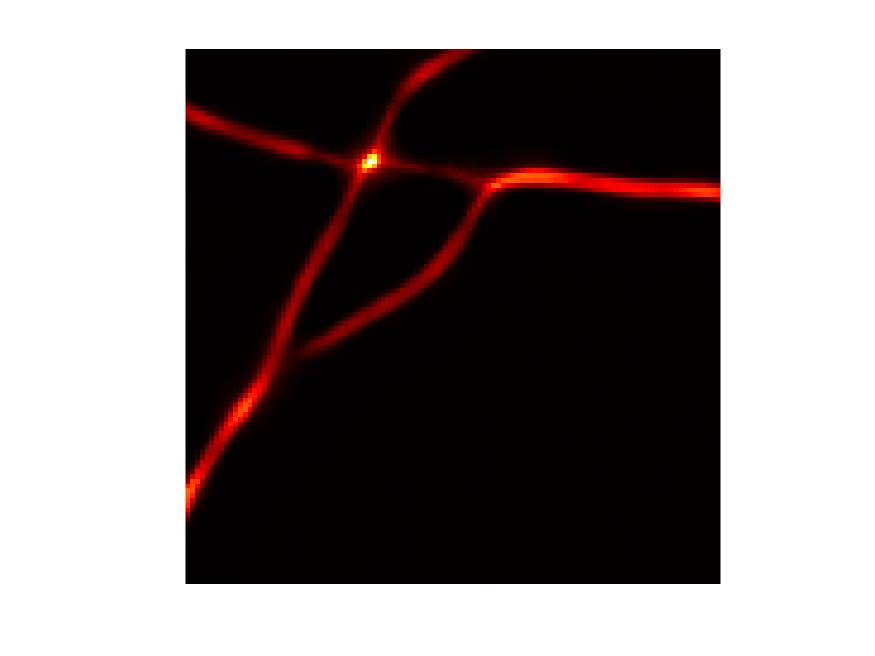}};
\spy[green] on (+0.82,0.6) in node [left] at (1.53,-1.7);
\spy[yellow] on (-0.7,-0.4) in node [left] at (0,-1.7);
\end{tikzpicture}}
& \adjustbox{valign=m,vspace=1pt}{\begin{tikzpicture}[spy using outlines={rectangle,green,magnification=1.7,size=1.4cm,connect spies}]
\node {\pgfimage[width=4cm]{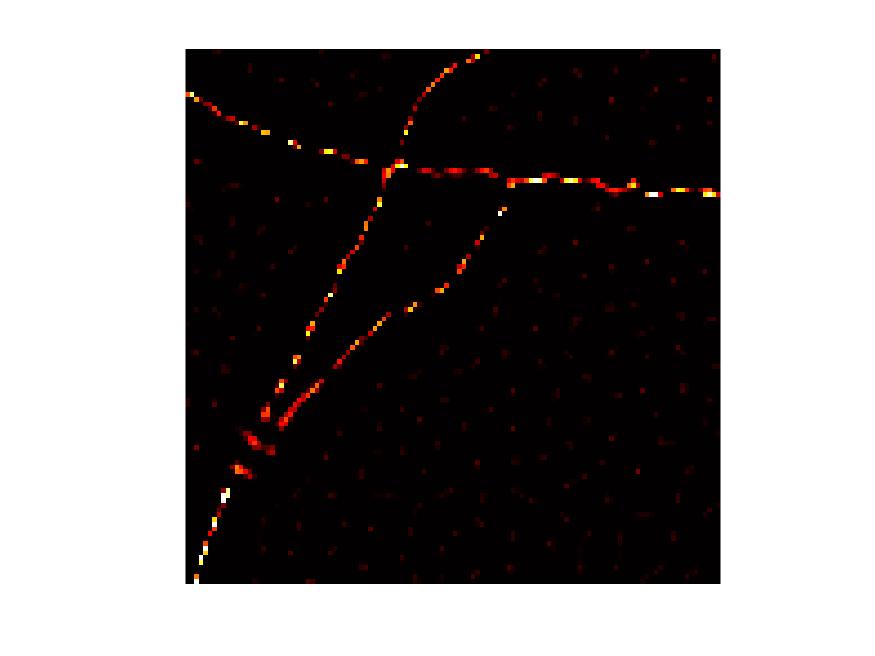}};
\spy[green] on (+0.82,0.6) in node [left] at (1.53,-1.7);
\spy[yellow] on (-0.7,-0.4) in node [left] at (0,-1.7);
\end{tikzpicture}}
& \adjustbox{valign=m,vspace=1pt}{\begin{tikzpicture}[spy using outlines={rectangle,green,magnification=1.7,size=1.4cm,connect spies}]
\node {\pgfimage[width=4cm]{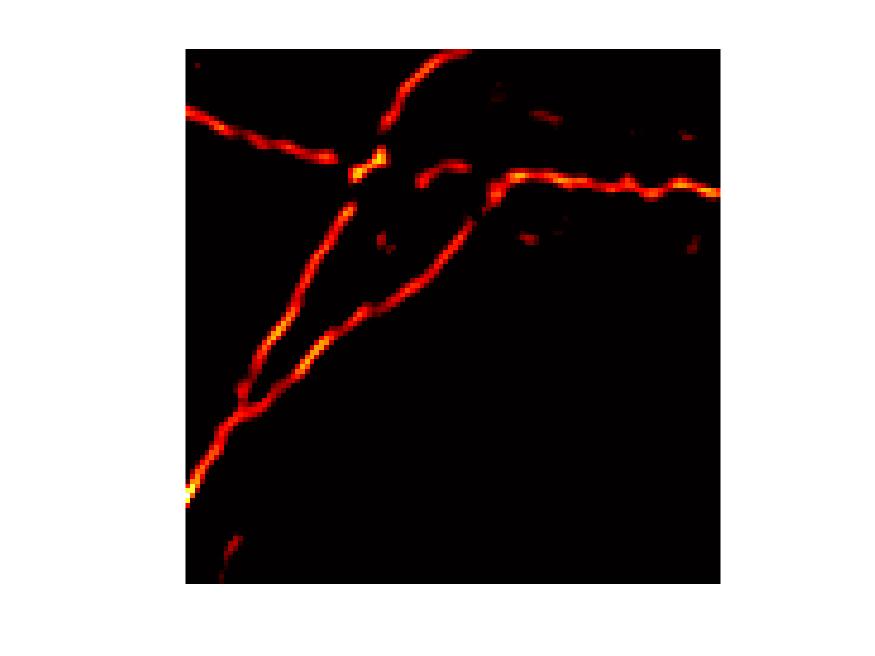}};
\spy[green] on (+0.82,0.6) in node [left] at (1.53,-1.7);
\spy[yellow] on (-0.7,-0.4) in node [left] at (0,-1.7);
\end{tikzpicture}}\\
\adjustbox{valign=m,vspace=1pt}{\includegraphics[width=.28\textwidth]{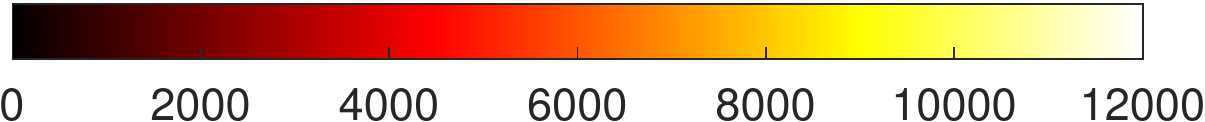}}
& \adjustbox{valign=m,vspace=1pt}{\includegraphics[width=.28\textwidth]{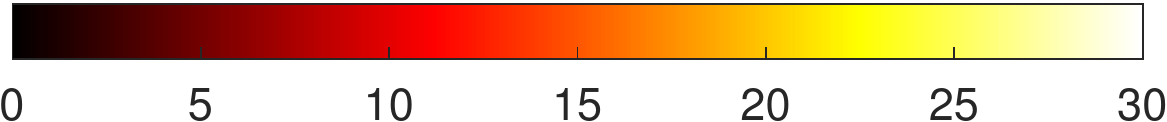}}
& \adjustbox{valign=m,vspace=1pt}{\includegraphics[width=.28\textwidth]{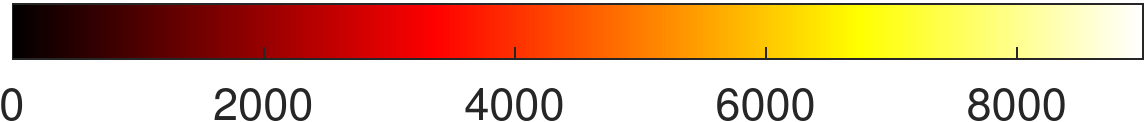}}
\end{tabular}

 \caption{Results for the low-background (LB) dataset with $T=500$. Note that the methods SRRF, SPARCOM and LSPARCOM do not estimate real intensity values. Between the compared methods only COL0RME is capable of estimating them, while the other methods estimate the mean of a radiality image sequence (SRRF) and normalized autocovariances (SPARCOM, LSPARCOM)} 
 %(see Figure \ref{fig:colorbars} for the colorbar of the GT and the COL0RME images)
 \label{fig: results LB}
 \end{figure}

\begin{figure}[h!]
\centering
\setlength\tabcolsep{1.5pt}
\renewcommand{\arraystretch}{1.5}
\begin{tabular}{ccc}
GT image & \small{COL0RME-CEL0} & \small{COL0RME-$\ell_1$}\\
\adjustbox{valign=m,vspace=1pt}{\begin{tikzpicture}[spy using outlines={rectangle,magnification=1.7,size=1.4cm,connect spies}]
\node {\pgfimage[width=4cm]{Figures/scalebar_GT_highBg14_SNR15_K500_N34.jpg}};
\spy[green] on (+0.82,0.6) in node [left] at (1.53,-1.7);
\spy[yellow] on (-0.7,-0.4) in node [left] at (0,-1.7);
\end{tikzpicture}}
& \adjustbox{valign=m,vspace=1pt}{\begin{tikzpicture}[spy using outlines={rectangle,magnification=1.7,size=1.4cm,connect spies}]
\node {\pgfimage[width=4cm]{Figures/COL0RME_int_highBg14_SNR15_K500_N34_CEL0_0_RP2.jpg}};
\spy[green] on (+0.82,0.6) in node [left] at (1.53,-1.7);
\spy[yellow] on (-0.7,-0.4) in node [left] at (0,-1.7);
\end{tikzpicture}}
& \adjustbox{valign=m,vspace=1pt}{\begin{tikzpicture}[spy using outlines={rectangle,magnification=1.7,size=1.4cm,connect spies}]
\node {\pgfimage[width=4cm]{Figures/COL0RME_int_highBg14_SNR15_K500_N34_L1_RP2.jpg}};
\spy[green] on (+0.82,0.6) in node [left] at (1.53,-1.7);
\spy[yellow] on (-0.7,-0.4) in node [left] at (0,-1.7);
\end{tikzpicture}}\\
\adjustbox{valign=m,vspace=1pt}{\includegraphics[width=.28\textwidth]{Figures/colorbarGT-eps-converted-to}}
& \adjustbox{valign=m,vspace=1pt}{\includegraphics[width=.28\textwidth]{Figures/colorbarGT-eps-converted-to}}
& \adjustbox{valign=m,vspace=1pt}{\includegraphics[width=.28\textwidth]{Figures/colorbarGT-eps-converted-to}}\\
SRRF & SPARCOM & LSPARCOM\\
\adjustbox{valign=m,vspace=1pt}{\begin{tikzpicture}[spy using outlines={rectangle,green,magnification=1.7,size=1.4cm,connect spies}]
\node {\pgfimage[width=4cm]{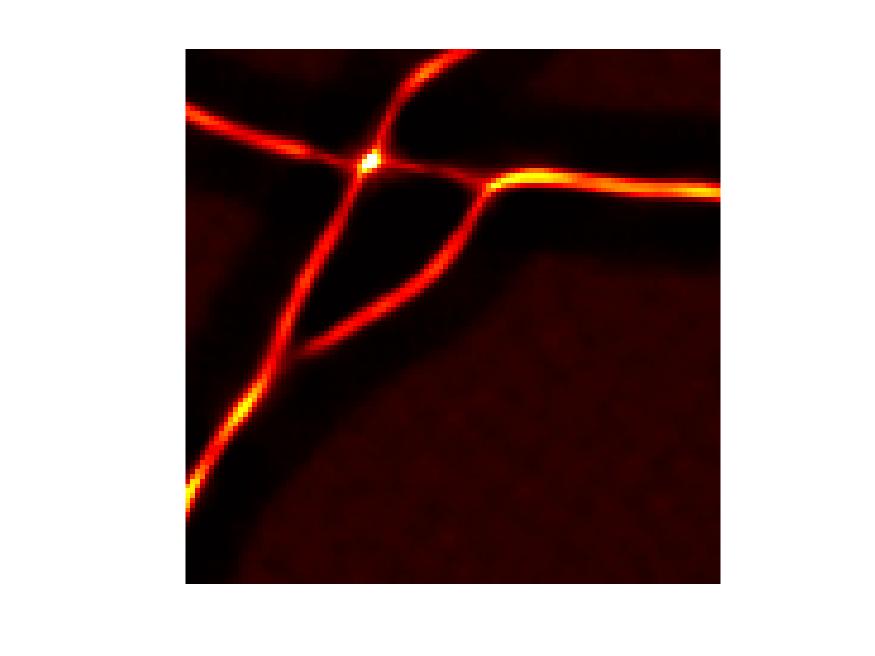}};
\spy[green] on (+0.82,0.6) in node [left] at (1.53,-1.7);
\spy[yellow] on (-0.7,-0.4) in node [left] at (0,-1.7);
\end{tikzpicture}}
& \adjustbox{valign=m,vspace=1pt}{\begin{tikzpicture}[spy using outlines={rectangle,green,magnification=1.7,size=1.4cm,connect spies}]
\node {\pgfimage[width=4cm]{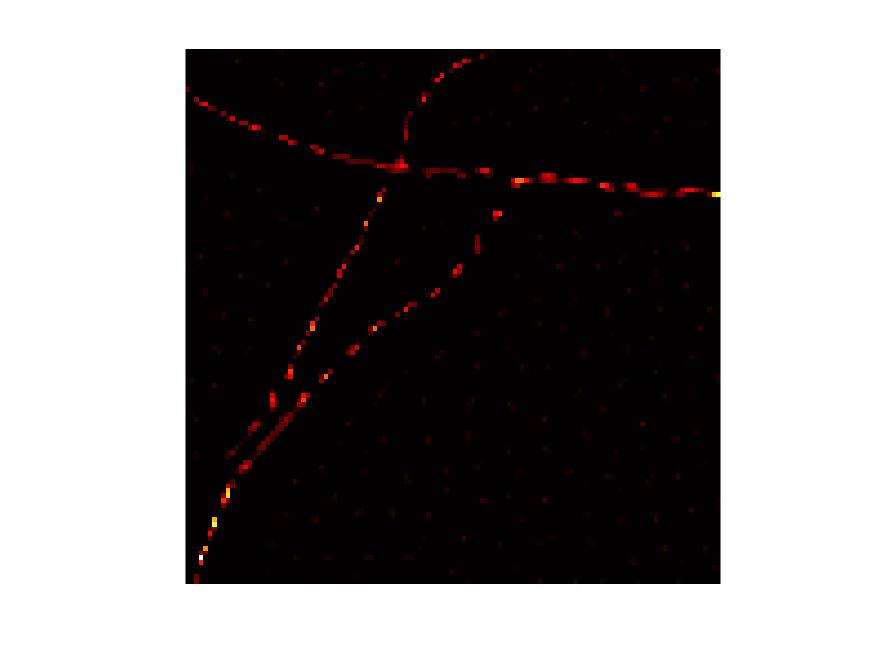}};
\spy[green] on (+0.82,0.6) in node [left] at (1.53,-1.7);
\spy[yellow] on (-0.7,-0.4) in node [left] at (0,-1.7);
\end{tikzpicture}}
& \adjustbox{valign=m,vspace=1pt}{\begin{tikzpicture}[spy using outlines={rectangle,green,magnification=1.7,size=1.4cm,connect spies}]
\node {\pgfimage[width=4cm]{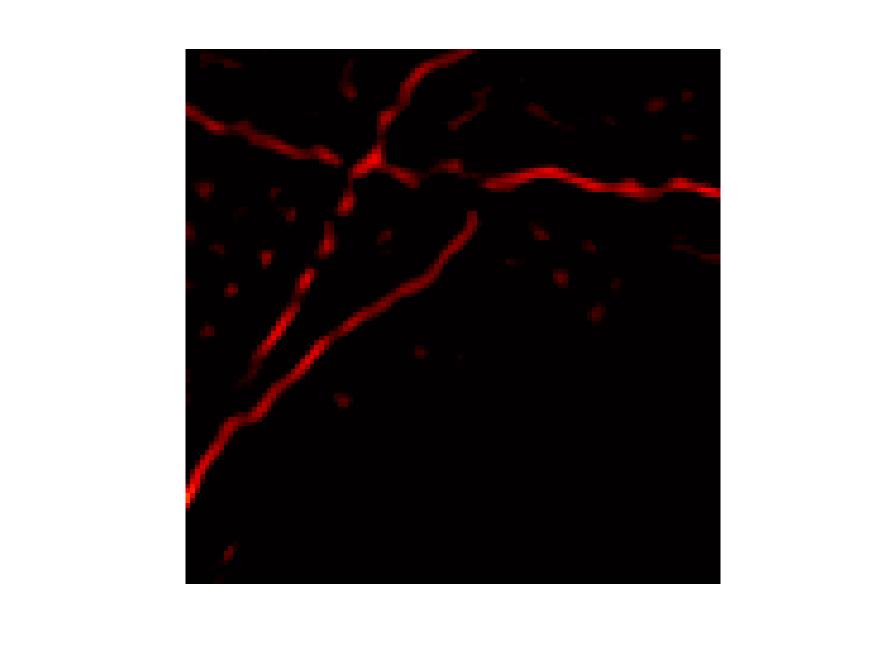}};
\spy[green] on (+0.82,0.6) in node [left] at (1.53,-1.7);
\spy[yellow] on (-0.7,-0.4) in node [left] at (0,-1.7);
\end{tikzpicture}}\\
\adjustbox{valign=m,vspace=1pt}{\includegraphics[width=.28\textwidth]{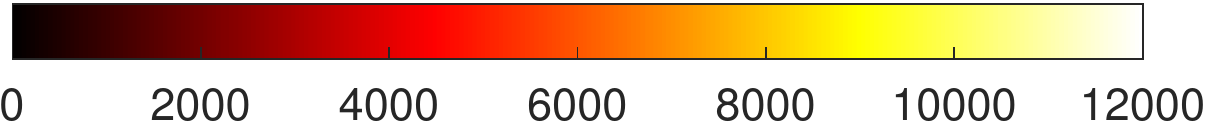}}
& \adjustbox{valign=m,vspace=1pt}{\includegraphics[width=.28\textwidth]{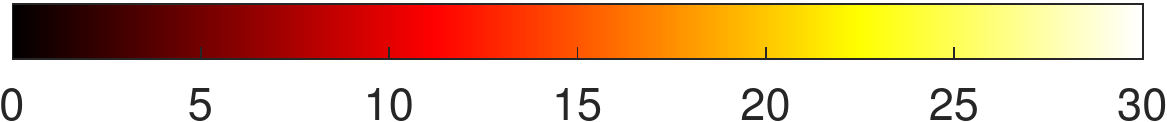}}
& \adjustbox{valign=m,vspace=1pt}{\includegraphics[width=.28\textwidth]{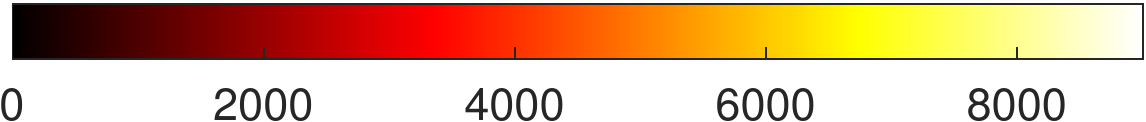}}
\end{tabular}
 \caption{Results for the high-background(HB) dataset with $T=500$} 
%  The colorbar of the GT and the COL0RME images is available in Figure \ref{fig:colorbars} 
  \label{fig: results HB}
 \end{figure}

In Figure \ref{fig: results LB} we compare the reconstructions of the methods COL0RME-CEL0, COL0RME-$\ell_1$, SRRF, SPARCOM and LSPARCOM for the LB dataset  and in Figure \ref{fig: results HB} for the HB dataset, for a sequence of T = $500$ frames. Results for different stack sizes, are available in the Supplementary Figures S1, S2 and S3. \mdffirst{Quantitative metrics like the Jaccard Index (JI) for the localization precision and the Peak-Signal-to-Noise-ration (PSNR) for the evaluation of the estimated intensities, are only available for the methods COL0RME-CEL0 and COL0RME-$\ell_1$ (see Figures \ref{fig: JI}, \ref{psnr}). For the rest of the methods, the JI values are very small due to background and noise artifacts in the reconstructions that lead to the appearance of many false positives, while the PSNR metric is not possible to be computed as the methods SRRF, SPARCOM and LSPARCOM do not reconstruct the intensity level.} In both datasets, LB and HB dataset, and for a sequence of T= $500$ frames, the better reconstruction, visually, is the one of the method COL0RME-$\ell_1$, as it is able to achieve a more clear separation of the filaments in the critical regions (yellow and green zoom boxes). The method COL0RME-CEL0 achieves also a good result, eventhough the separation of the filaments, that are magnified in the green box, is not so obvious. The same happens also when the method SPARCOM is being used. Finally, the reconstruction of the methods SRRF and LSPARCOM, is slightly misleading.

\subsection{Real Data}
\label{sec: Real Data}

To show the effectiveness of our method for handling real-world  data, we apply COL0RME to an image sequence acquired from a Total Internal Reflection Fluorescence (TIRF) microscope. The TIRF microscope offers a good observation of the activities happening next to the cell membrane, as it uses an evanescent wave to illuminate and excite fluorescent molecules only in this restricted region of the specimen\cite{axelrod2001total}. Further, the TIRF microscope does not require specific fluorescent dyes, allows live cell imaging by using a low illumination laser, with really low out-of-focus contribution and produces images with a relatively good, in comparison with other fluorescence microscopy techniques, SNR. To enhance the resolution of the images acquired from a TIRF microscope, super-resolution approaches that exploit the temporal fluctuations of blinking/fluctuating fluorophores, like COL0RME, can be applied.

\begin{figure}
\centering
\setlength\tabcolsep{1.5pt}
\begin{tabular}{ccc}
    $\bar{\mathbf{y}}$ & \mdfsecrev{{\small{COL0RME-CEL0}} ($\mathbf{x}$)} & \mdfsecrev{{\small{COL0RME-CEL0}} ($\mathbf{b}$)} \\
    \adjustbox{valign=m,vspace=1pt}{\begin{tikzpicture}[spy using outlines={rectangle}]
    \node {\includegraphics[width=0.3\textwidth]{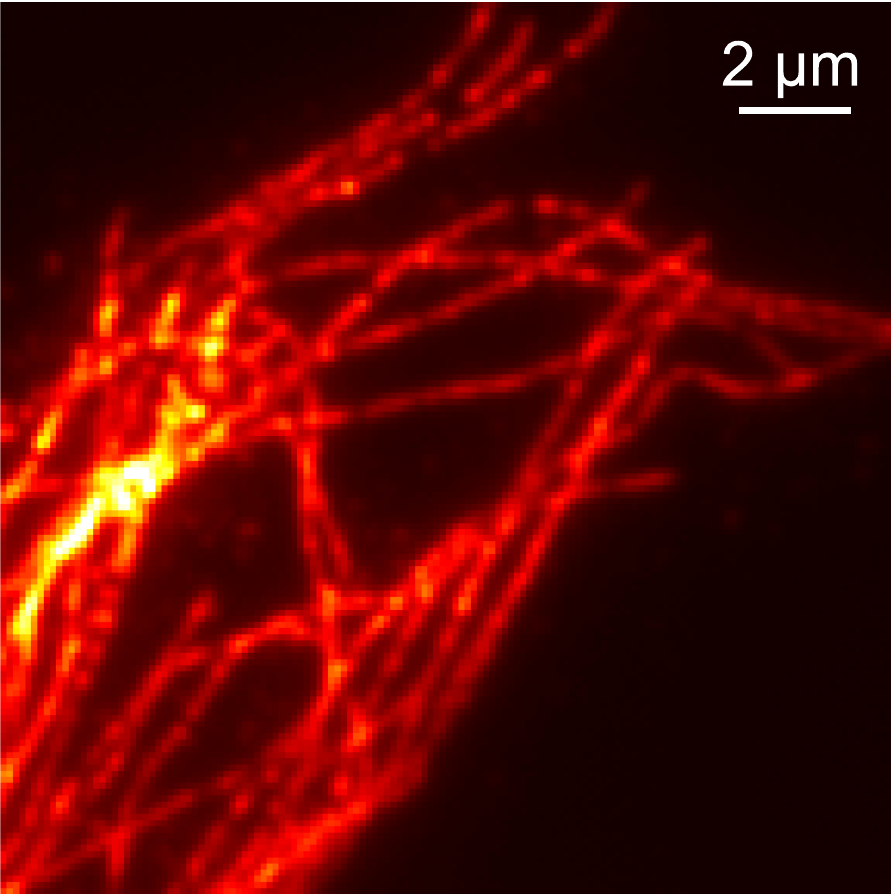}};
    \spy[green,magnification=1.7,size=0.8cm] on (-1.55,0.2) in node [left] at (-1.2,1.65);
    \spy[yellow,magnification=1.7,size=1.4cm] on (0.85,0.75) in node [left] at (2.1,-1.4);
    \end{tikzpicture}}
    & \adjustbox{valign=m,vspace=1pt}{\begin{tikzpicture}[spy using outlines={rectangle}]
    \node {\includegraphics[width=0.3\textwidth]{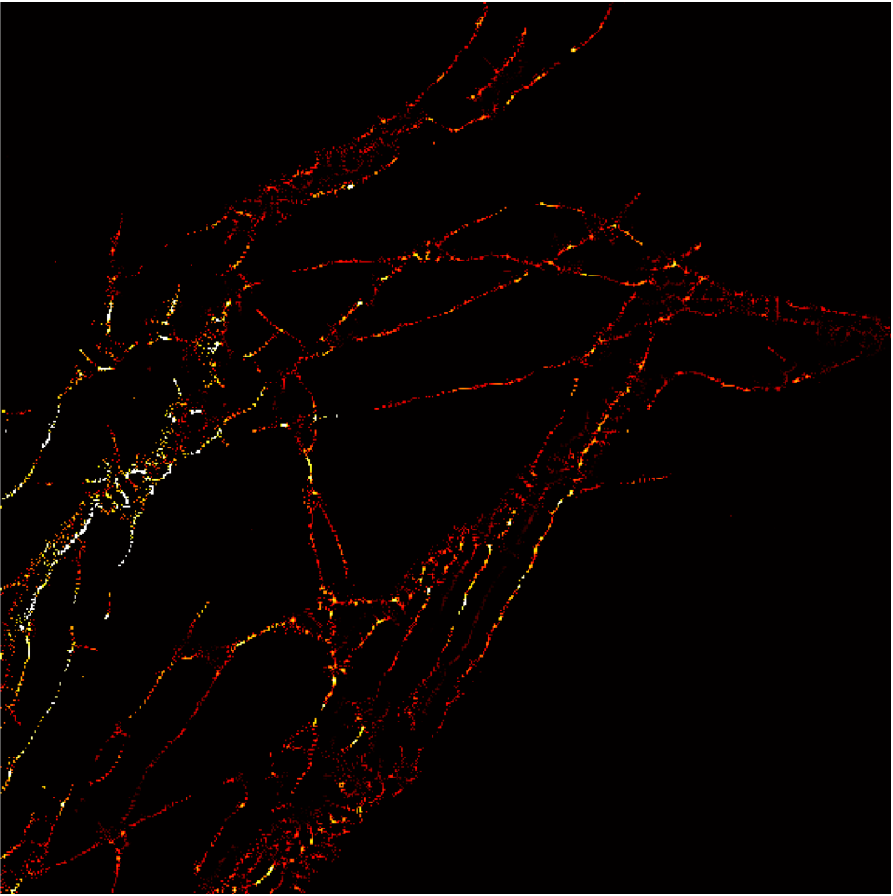}};
    \spy[green,magnification=1.7,size=0.8cm] on (-1.55,0.2) in node [left] at (-1.2,1.65);
    \spy[yellow,magnification=1.7,size=1.4cm] on (0.85,0.75) in node [left] at (2.1,-1.4);
    \end{tikzpicture}}
    &\adjustbox{valign=m,vspace=1pt}{\includegraphics[width=0.3\textwidth]{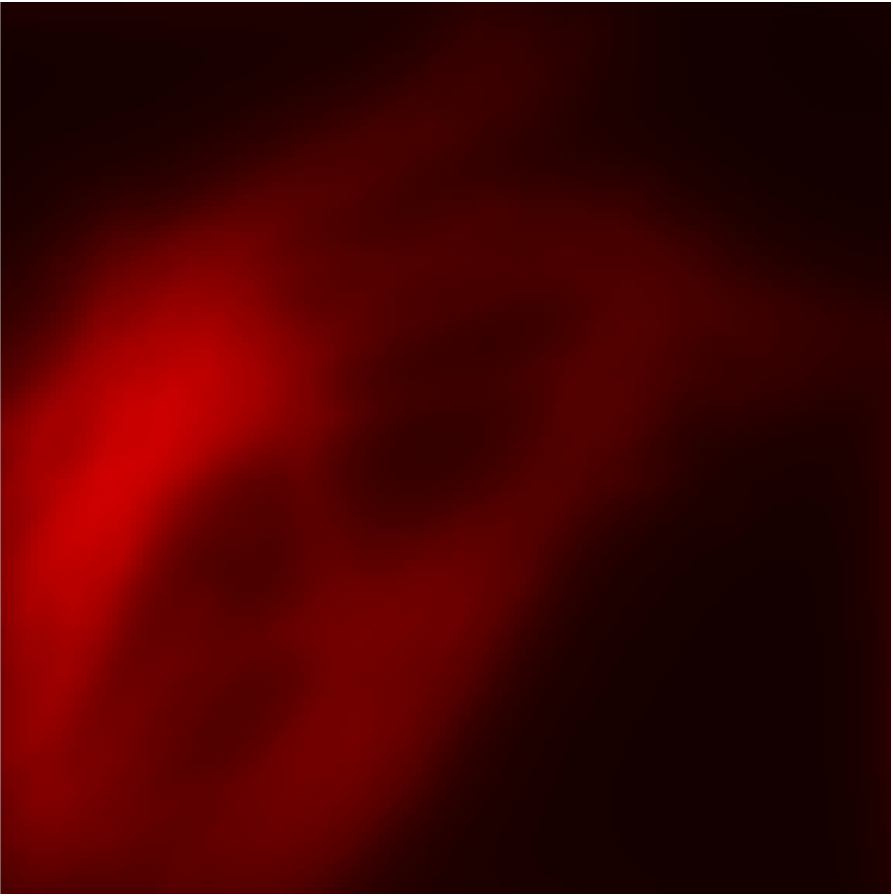}}\\
    \adjustbox{valign=m,vspace=1pt}{\includegraphics[width=.28\textwidth]{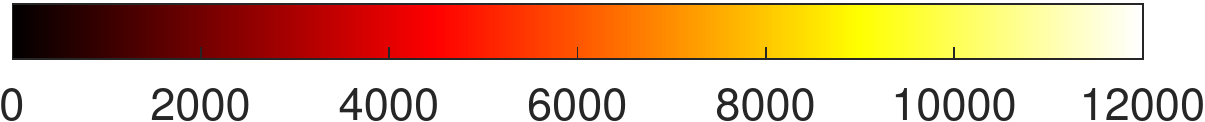}} 
    & \adjustbox{valign=m,vspace=1pt}{\includegraphics[width=.28\textwidth]{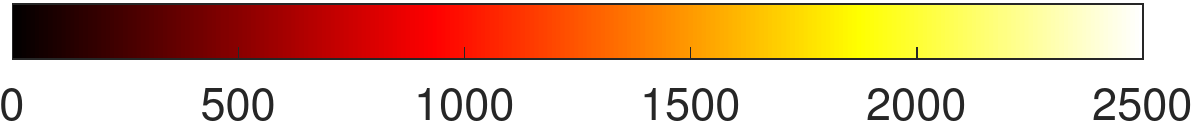}} 
    & \adjustbox{valign=m,vspace=1pt}{\includegraphics[width=.28\textwidth]{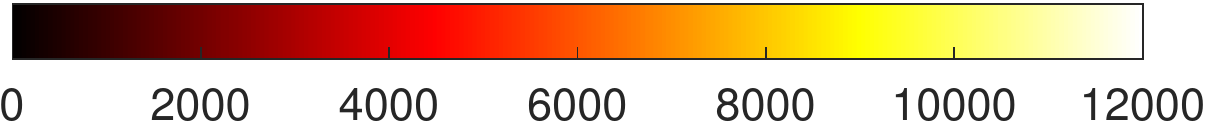}}\\
    $\mathbf{y}_t$ & \mdfsecrev{{\small{COL0RME-$\ell_1$}} ($\mathbf{x}$)} & \mdfsecrev{{\small{COL0RME-$\ell_1$}} ($\mathbf{b}$)} \\
    \adjustbox{valign=m,vspace=1pt}{\begin{tikzpicture}[spy using outlines={rectangle}]
    \node {\includegraphics[width=0.3\textwidth]{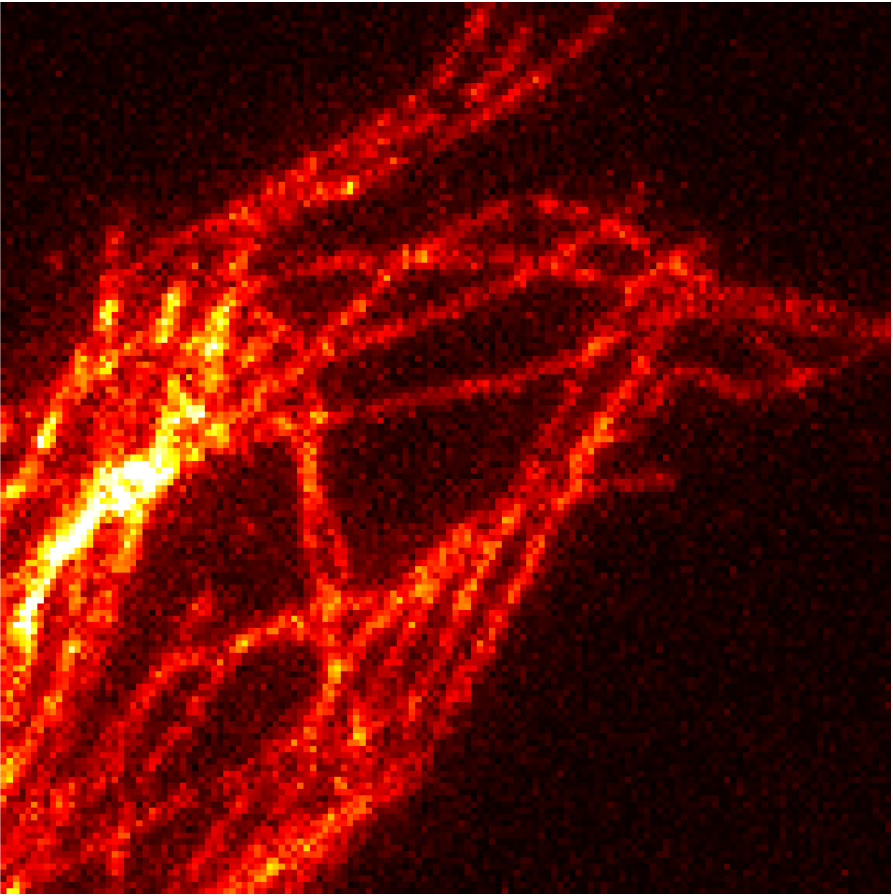}};
    \spy[green,magnification=1.7,size=0.8cm] on (-1.55,0.2) in node [left] at (-1.2,1.65);
    \spy[yellow,magnification=1.7,size=1.4cm] on (0.85,0.75) in node [left] at (2.1,-1.4);
    \end{tikzpicture}}
    & \adjustbox{valign=m,vspace=1pt}{\begin{tikzpicture}[spy using outlines={rectangle}]
    \node {\includegraphics[width=0.3\textwidth]{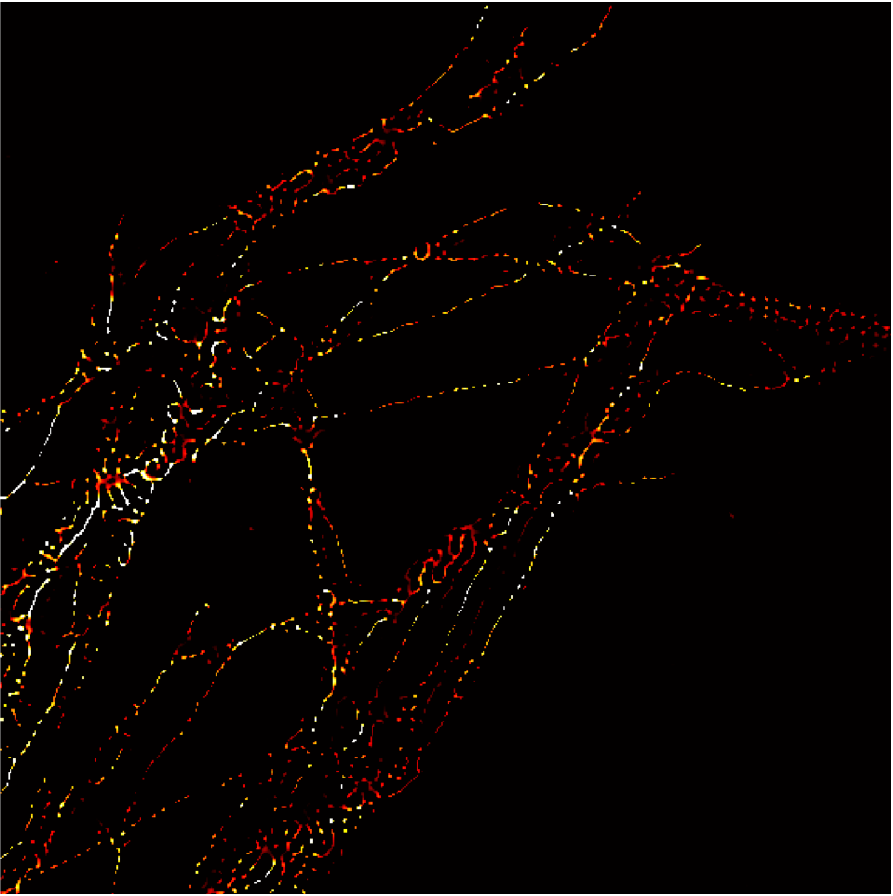}};
    \spy[green,magnification=1.7,size=0.8cm] on (-1.55,0.2) in node [left] at (-1.2,1.65);
    \spy[yellow,magnification=1.7,size=1.4cm] on (0.85,0.75) in node [left] at (2.1,-1.4);  
    \end{tikzpicture}}   
    &\adjustbox{valign=m,vspace=1pt}{\includegraphics[width=0.3\textwidth]{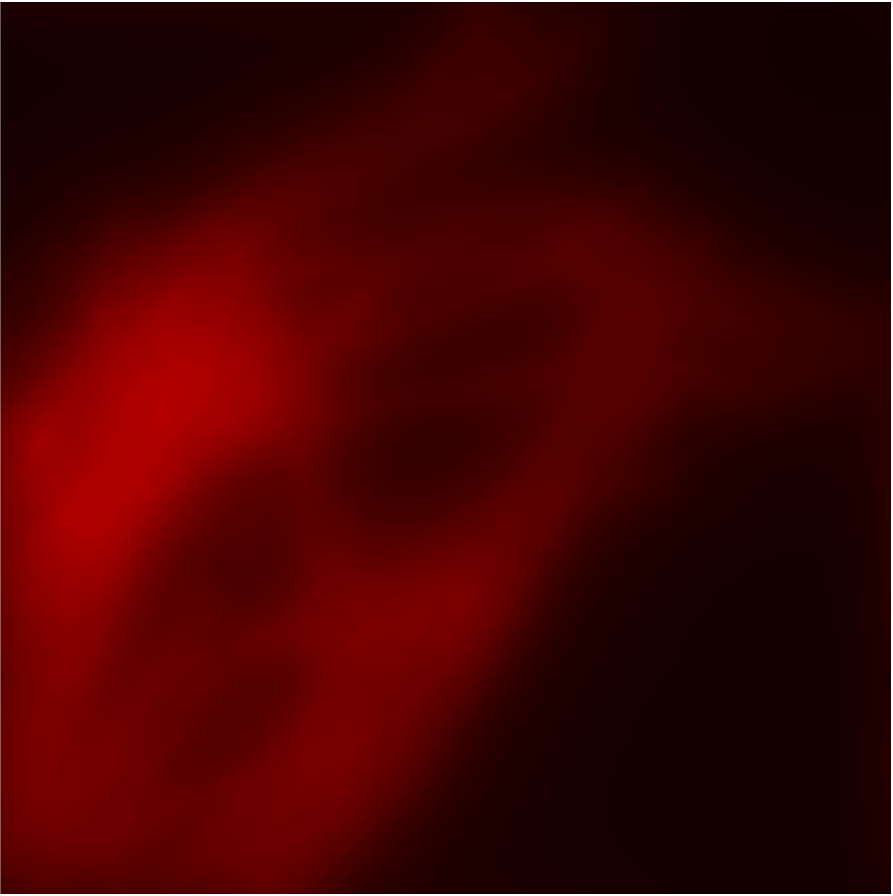}}\\
    \adjustbox{valign=m,vspace=1pt}{\includegraphics[width=.28\textwidth]{Figures/colorbar_COL0RME_realData_coarsegrid-eps-converted-to}} 
    & \adjustbox{valign=m,vspace=1pt}{\includegraphics[width=.28\textwidth]{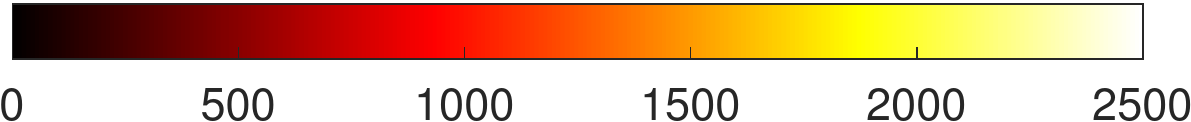}} 
    & \adjustbox{valign=m,vspace=1pt}{\includegraphics[width=.28\textwidth]{Figures/colorbar_COL0RME_realData_coarsegrid-eps-converted-to}}\\
    \end{tabular}
\caption{Real TIRF data, $T=500$ frames. Diffraction limited image or the mean of the stack $\bar{\mathbf{y}}$ \mdfsecrev{(4x zoom)}, \mdffirst{A frame $\mathbf{y}_t$ from the stack \mdfsecrev{(4x zoom)},} The intensity and background estimation of the methods COL0RME-CEL0 and COL0RME-$\ell_1$}
\label{fig: real_L1_cel0} 
\end{figure}

The data we are using have been obtained from a Multi-Angle TIRF microscope, with a fixed angle close to the critical one. A sequence of $500$ frames has been acquired, with an acquisition time equal to $25$s. \mdf{Tubulins in endothelial cells are being imaged, while they are colored with the Alexa Fluor 488. The variance of fluctuations over time for a typical pixel is measured and is belonging to the range $5\times10^{5}-7\times10^{5}$.} The diffraction limited image, or with other words the mean stack image $\bar{\mathbf{y}}$ is shown in Figure \ref{fig: real_L1_cel0}, together with  \mdffirst{one frame $\mathbf{y}_t$ extracted from the entire stack.} The FWHM of the PSF has been measured experimentally and is equal to $292.03$nm, while the CCD camera has a pixel of size $106$nm.  

The results of the method COL0RME-CEL0 and COL0RME-$\ell_1$ 
%for a sequence of 500 frames 
and more precisely the intensity and the background estimation, can be found in Figure \ref{fig: real_L1_cel0}. \mdf{Experiments using different stack sizes have been done showing 
%similarly to the simulated data presented in this paper (Figures \ref{fig: results LB},\ref{fig: results HB}) but also in the supplementary material (Figures S1, S2, S3) 
that the more frames we use (up to a point that we do not have many molecules bleached), the more continuous filaments we find. However, by acquiring only 500 frames we have a good balance between temporal and spatial resolution. For this reason we present here only results using a stack of 500 frames.}
For the method COL0RME-CEL0 the regularization parameter $\lambda$ is equal to $\lambda = 5 \times 10^{-4} \times \lambda_{max}^{CEL0}$ and the algorithmic restarting approach has been used \mdf{(stopping criteria: when, in a certain restarting, there are not additional pixels added to the global support, but with maximum 10 restarts)}. Regarding the method COL0RME-$\ell_1$ the regularization parameter $\lambda$ is equal to $\lambda = 5 \times 10^{-6} \times \lambda_{max}^{\ell_1}$, a relatively small value so as to be sure that we will include all the pixels that contain fluorescent molecules. Even if we underestimate $\lambda$ and find more false positives in the support estimation, after the second step of the algorithm, the final reconstruction is corrected, as explained in \ref{parameter_lambda}. The hyper-parameters $\alpha$ and $\beta$ are equal to: $\alpha = 10^6$, $\beta = 20$. Using any of the two regularizers the spatial resolution is enhanced, as it can be also observed from the yellow zoom boxes. \mdf{However, the reconstruction obtained by both COL0RME-CEL0 and COL0RME-$\ell_1$ is to some degree punctuated due to mainly limitations arising from experimental difficulties to get a staining sufficiently homogeneous for this imaging resolution.} Furthermore, there are a few filaments that do not seem to be well reconstructed, especially using the COL0RME-CEL0 method, e.g. the one inside the green box. 

\begin{figure}
\centering
\setlength\tabcolsep{1.5pt}
\begin{tabular}{ccc}
    $\bar{\mathbf{y}}$ & COL0RME-CEL0 & COL0RME-$\ell_1$ \\
    \adjustbox{valign=m,vspace=1pt}{\begin{tikzpicture}[spy using outlines={rectangle,green,magnification=1.7,size=1.4cm}]
    \node {\includegraphics[width=0.3\textwidth]{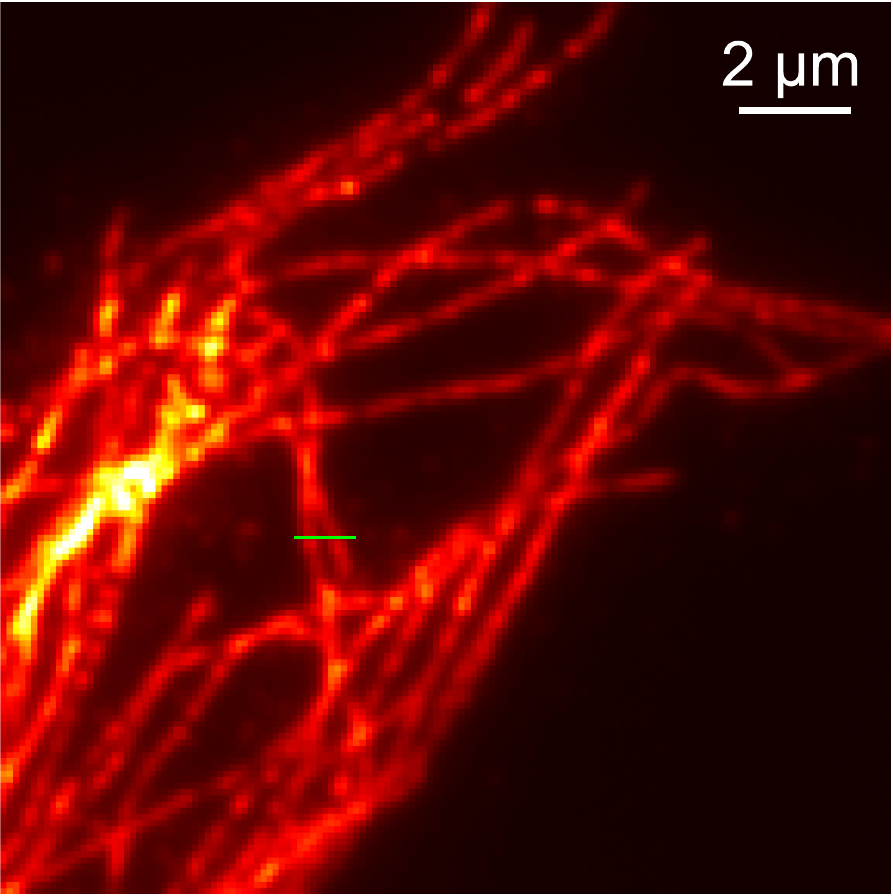}};
    \spy[green] on (0.85,0.75) in node [left] at (2.2,-2.9);
    \spy[blue] on (-0.4,-0.75) in node [left] at (0.7,-2.9);
    \spy[lime] on (-0.7,0.3) in node [left] at (-0.8,-2.9);
    \end{tikzpicture}}
    &\adjustbox{valign=m,vspace=1pt}{\begin{tikzpicture}[spy using outlines={rectangle,green,magnification=1.7,size=1.4cm}]
    \node {\includegraphics[width=0.3\textwidth]{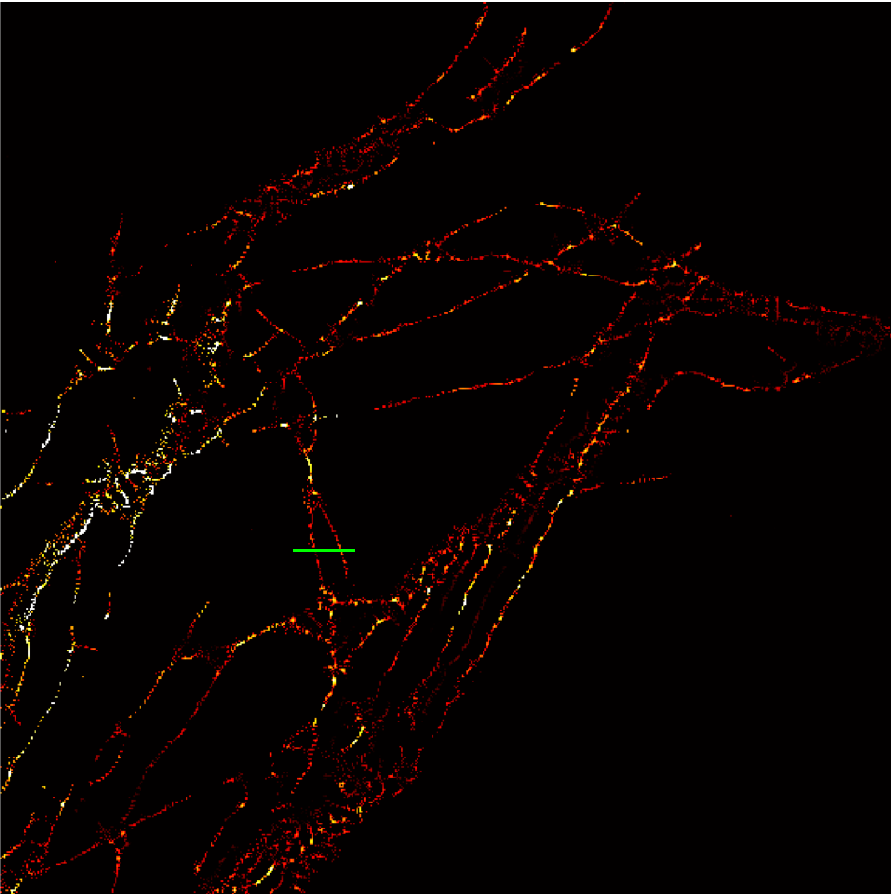}};
    \spy[green] on (0.85,0.75) in node [left] at (2.2,-2.9);
    \spy[blue] on (-0.4,-0.75) in node [left] at (0.7,-2.9);
    \spy[lime] on (-0.7,0.3) in node [left] at (-0.8,-2.9);
    \end{tikzpicture}}
    &\adjustbox{valign=m,vspace=1pt}{\begin{tikzpicture}[spy using outlines={rectangle,green,magnification=1.7,size=1.4cm}]
    \node {\includegraphics[width=0.3\textwidth]{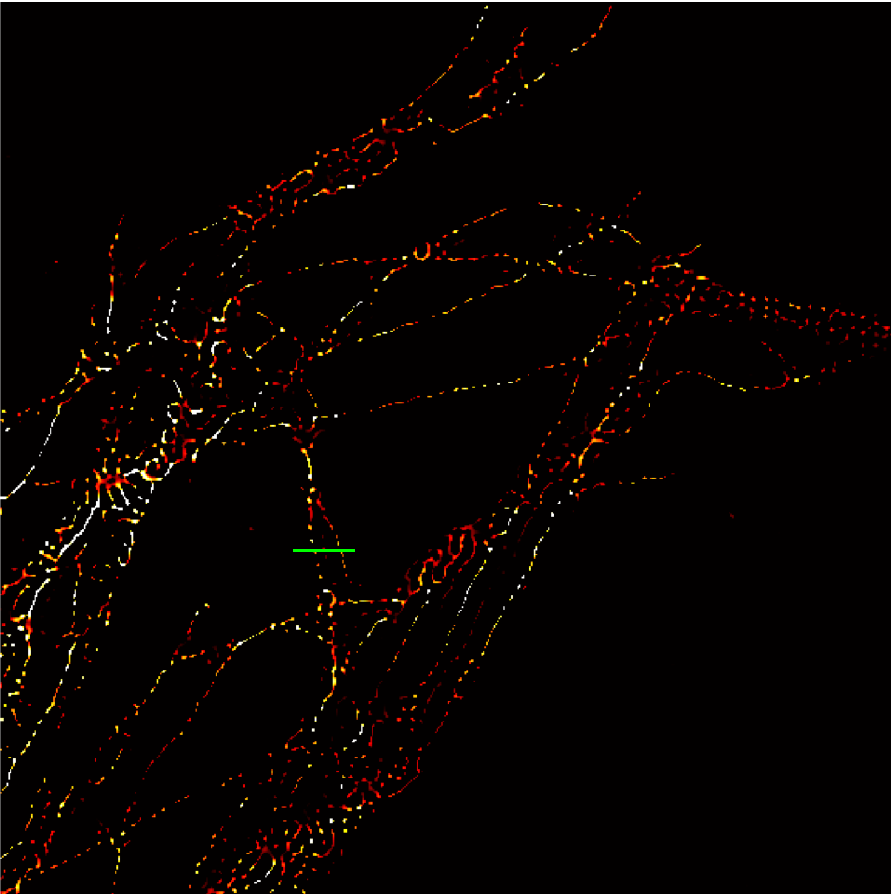}};
    \spy[green] on (0.85,0.75) in node [left] at (2.2,-2.9);
    \spy[blue] on (-0.4,-0.75) in node [left] at (0.7,-2.9);
    \spy[lime] on (-0.7,0.3) in node [left] at (-0.8,-2.9);
    \end{tikzpicture}} \\
    \adjustbox{valign=m,vspace=1pt}{\includegraphics[width=.28\textwidth]{Figures/colorbar_COL0RME_realData_coarsegrid-eps-converted-to}} 
    & \adjustbox{valign=m,vspace=1pt}{\includegraphics[width=.28\textwidth]{Figures/colorbar_COL0RME_realData-eps-converted-to}} 
    & \adjustbox{valign=m,vspace=1pt}{\includegraphics[width=.28\textwidth]{Figures/colorbar_COL0RME_realData-eps-converted-to}}\\
    SRRF & SPARCOM & LSPARCOM \\
    \adjustbox{valign=m,vspace=1pt}{\begin{tikzpicture}[spy using outlines={rectangle,green,magnification=1.7,size=1.4cm}]
    \node {\includegraphics[width=0.3\textwidth]{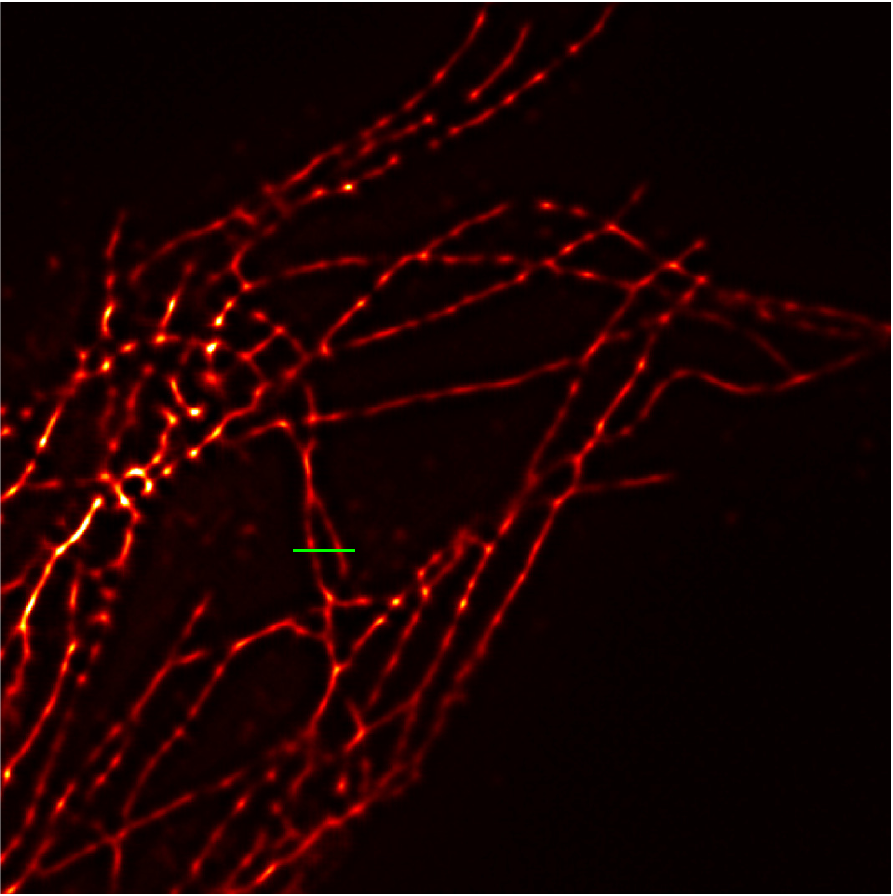}};
    \spy[green] on (0.85,0.75) in node [left] at (2.2,-2.9);
    \spy[blue] on (-0.4,-0.75) in node [left] at (0.7,-2.9);
    \spy[lime] on (-0.7,0.3) in node [left] at (-0.8,-2.9);
    \end{tikzpicture}}
    &\adjustbox{valign=m,vspace=1pt}{\begin{tikzpicture}[spy using outlines={rectangle,green,magnification=1.7,size=1.4cm}]
    \node {\includegraphics[width=0.3\textwidth]{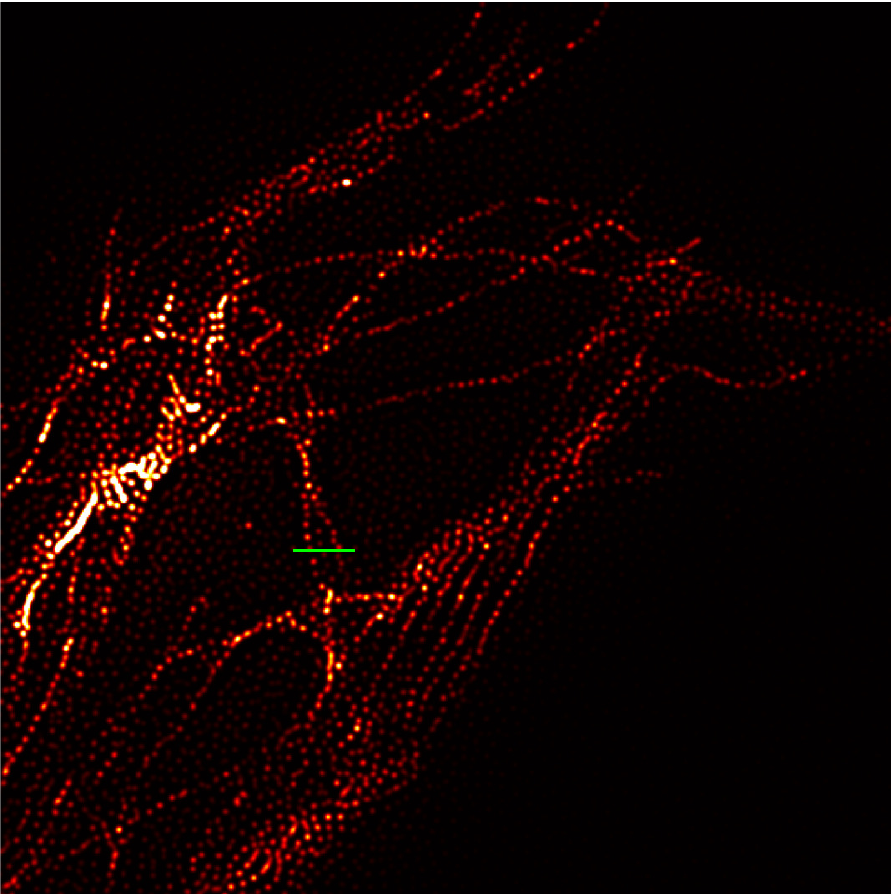}};
    \spy[green] on (0.85,0.75) in node [left] at (2.2,-2.9);
    \spy[blue] on (-0.4,-0.75) in node [left] at (0.7,-2.9);
    \spy[lime] on (-0.7,0.3) in node [left] at (-0.8,-2.9);
    \end{tikzpicture}} 
    &\adjustbox{valign=m,vspace=1pt}{\begin{tikzpicture}[spy using outlines={rectangle,green,magnification=1.7,size=1.4cm}]
    \node {\includegraphics[width=0.3\textwidth]{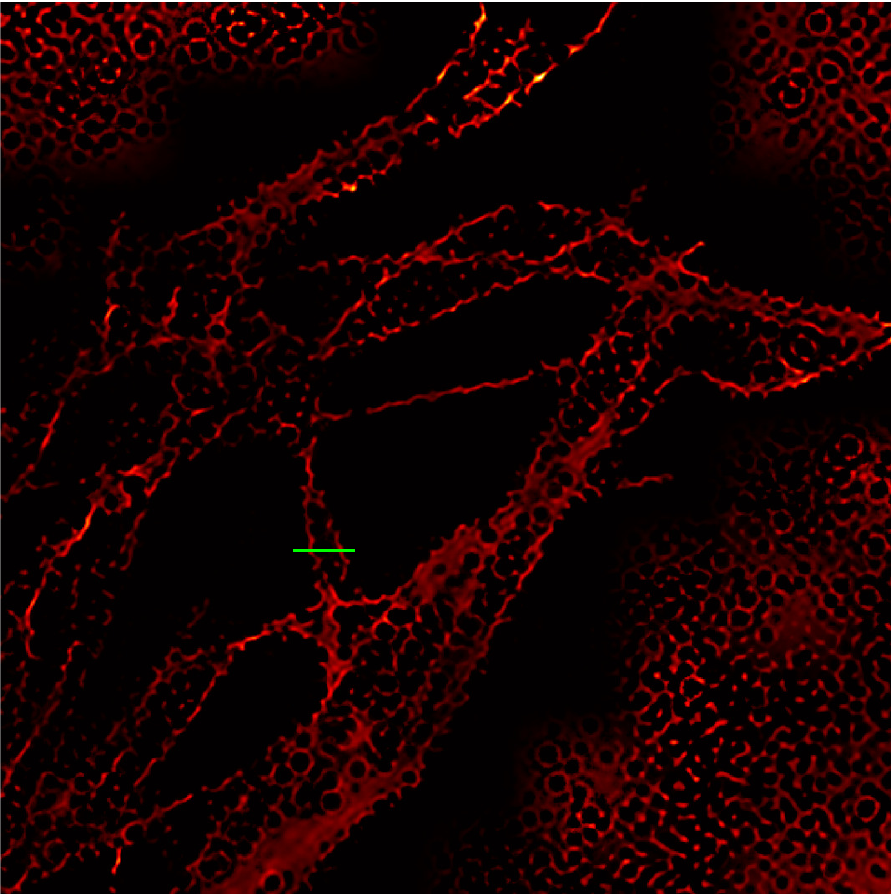}};
    \spy[green] on (0.85,0.75) in node [left] at (2.2,-2.9);
    \spy[blue] on (-0.4,-0.75) in node [left] at (0.7,-2.9);
    \spy[lime] on (-0.7,0.3) in node [left] at (-0.8,-2.9);
    \end{tikzpicture}}\\
    \adjustbox{valign=m,vspace=1pt}{\includegraphics[width=.28\textwidth]{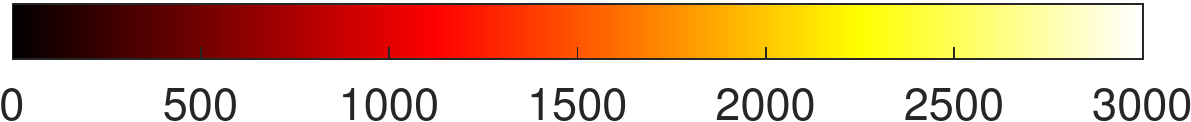}}
    & \adjustbox{valign=m,vspace=1pt}{\includegraphics[width=.28\textwidth]{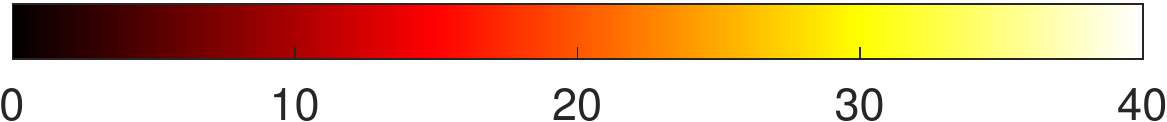}}
    & \adjustbox{valign=m,vspace=1pt}{\includegraphics[width=.28\textwidth]{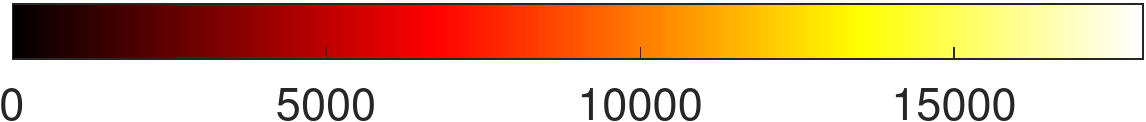}}\\
\end{tabular}
\centering
\includegraphics[width=.42\textwidth]{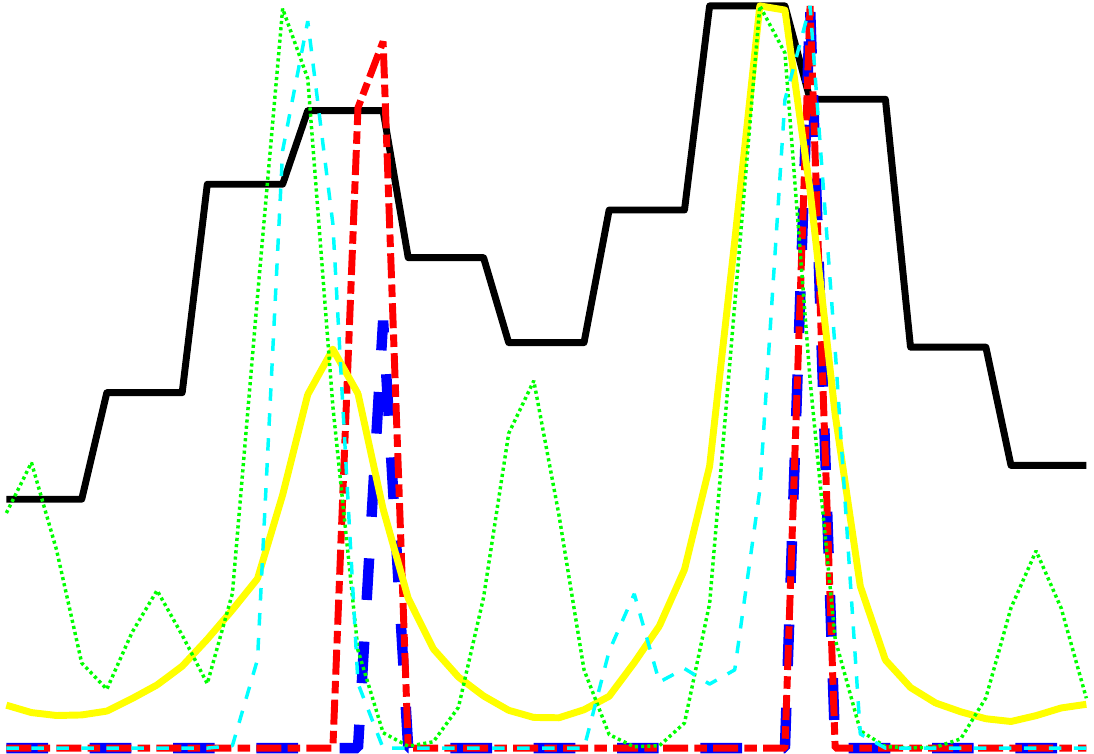}
\hspace{0.2cm}
\includegraphics[width=.28\textwidth]{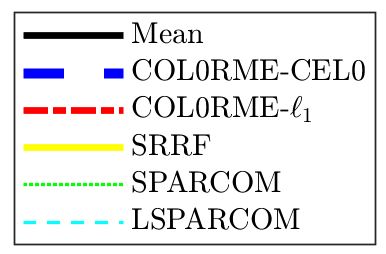}
\caption{Real TIRF data, $T=500$ frames. Diffraction limited image $\bar{\mathbf{y}}$ (4x zoom), Comparisons between the method that exploit the temporal fluctuations, Normalized cross-section along the green line presented in the diffraction limited and reconstructed images, but also in the blue zoom-boxes. Discription of colorbars: real intensity values for $\overline{\mathbf{y}}$ and COL0RME in two different grids, mean of the radiality image sequence for SRRF, normalized autocovariances for SPARCOM and LSPARCOM}
\label{fig: real_compare}
\end{figure}

Finally, the comparison of the methods COL0RME-CEL0 and COL0RME-$\ell_1$ with the other state-of-the-art methods, is available in Figure \ref{fig: real_compare}. The parameters used for the methods SRRF, SPARCOM and LSPARCOM, are explained in the section \ref{sec: Simulated Data}. Here, we further use the post-processing step (convolution with a small Gaussian function) in the method SPARCOM, as the result was dotted. The methods COL0RME-CEL0 and COL0RME-$\ell_1$ seem to have the most precise localization, \mdf{by reconstructing thin filaments, as shown in the cross-section plotted in Figure \ref{fig: real_compare}, though a bit punctuated.} The most appealing visually is the result of the method SRRF, where the filaments have a more continuous structure, however from the cross-section, we can see that the resolution is not so much improved compared to the other methods . SPARCOM and LSPARCOM do not perform very well in this real image sequence due to, mainly, background artifacts.

\section{Discussion and Conclusions}

In this paper, we propose and discuss the model and the performance of COL0RME, a method for super-resolution microscopy imaging based on the sparse analysis of the stochastic fluctuations of molecules' intensities. Similarly to other methods exploiting temporal fluctuations, COL0RME relaxes all the requirements for special equipment (microscope and fluorescent dyes) and allows for live-cell imaging, due to the good temporal resolution and the low power lasers employed.
In comparison with competing methods, COL0RME achieves higher spatial resolution than other methods exploiting fluctuations while having a sufficient temporal resolution. COL0RME is based on two different steps: a former one where accurate molecule localization and noise estimation are achieved by solving non-smooth convex/non-convex optimization problems in the covariance domain and the latter where intensity information is retrieved in correspondence with the estimated support only. Our numerical results show that COL0RME outperforms competing approaches in terms of localization precision. 
%Moreover, differently from other methods, spatially-variant background information can further be estimated.  
To the best of our knowledge, \CL\ is the only super-resolution method exploiting temporal fluctuations which is capable of retrieving intensity-type information, signal and spatially-varying background, which are of fundamental interest in biological data analysis. For both steps, automatic parameter selection strategies are detailed. 
Let us remark that such strategy of intensity estimation could be applied to the other competing super-resolution methods in the literature.
Several results obtained on both simulated and real data are discussed, showing the superior performance of COL0RME in comparison with analogous methods such as SPARCOM, LSPARCOM and SRRF.
Possible extensions of this work shall address the use of intensity information estimated by COL0RME for 3D reconstruction in, e.g., MA-TIRF acquisitions. \mdfsec{Furthermore, a systematic study to assess quantitatively the spatial resolution achieved by COL0RME under different scenarios (different background levels, different PSNRs, number of frames) is envisaged.}

\begin{appendix}\appheader
\section{Appendix. Proximal computations}\label{appendixA}
Given the function $h: \mathbb{R}^{L^2} \rightarrow \mathbb{R}$, defined in (\ref{eq: h}), the proximal mapping of $h$ is a an operator given by:

\begin{align}
    \textbf{\text{prox}}_{h, \tau}(\mathbf{w}) 
    &= \argmin_{\mathbf{u}}\left(\frac{1}{2\tau}\|\mathbf{u - w}\|_2^2 + h(\mathbf{w})\right) \nonumber \\ 
    &= \argmin_{\mathbf{u}}\left(\frac{1}{2 \tau}\|\mathbf{u - w}\|_2^2 + \frac{\alpha}{2} \left(\|\mathbf{I_\Omega u}\|_2^2 + \sum_{i=1}^{L^2}\ [\phi(\mathbf{u}_i)]^2\right)\right).
    \label{eq: prox_h}
\end{align}

The optimal solution $\hat{\mathbf{u}}$ ($\hat{\mathbf{u}} = \textbf{\text{prox}}_{h,\tau}(\mathbf{w})$), as the problem \eqref{eq: prox_h} is convex, is attained when: 
\begin{align}
    \mathbf{0} &\in \nabla\left(\frac{1}{2\tau}\|\hat{\mathbf{u}}- \mathbf{w}\|_2^2+ \frac{\alpha}{2} \left(\|\mathbf{I_\Omega} \hat{\mathbf{u}}\|_2^2+ \sum_{i=1}^{L^2}\ [\phi(\hat{\mathbf{u}}_i)]^2\right)\right), \nonumber \\
    \mathbf{0} &\in \frac{1}{\tau}(\hat{\mathbf{u}} - \mathbf{w})+\alpha\left(\mathbf{I_\Omega} \hat{\mathbf{u}}+ \ [\phi(\hat{\mathbf{u}}_i) \phi'(\hat{\mathbf{u}}_i)]_{\{i=1,...,{L^2}\}}\right).
\end{align}

Starting from (\ref{phi}) we can compute $\phi': \mathbb{R}\rightarrow\mathbb{R}_+$, as:
\begin{equation}
    \phi'(z) := 
    \begin{cases}
      0 & \text{if } z \geq 0,\\
      1 & \text{if } z < 0,
    \end{cases}\qquad\forall z\in\mathbb{R}.  
    \label{eq:phi'}
\end{equation}

Given \eqref{eq:phi'}, we can write:
\begin{align}
    \mathbf{0} &\in \frac{1}{\tau}(\hat{\mathbf{u}} - \mathbf{w})+\alpha\left(\mathbf{I_\Omega} \hat{\mathbf{u}} + [\phi(\hat{\mathbf{u}}_i)]_{\{i=1,...,{L^2}\}}\right).
    %\mathbf{w} &= \hat{\mathbf{u}} + \alpha \tau \left(\mathbf{I_\Omega} \hat{\mathbf{u}} + [\phi(\hat{\mathbf{u}}_i)]_{\{i=1,...,{L^2}\}}\right)
    \label{eq: w}
\end{align}

Exploiting component-wise, as problem \eqref{eq: prox_h} is separable with respect to both $\mathbf{x}$ and $\mathbf{w}$, and assuming $\hat{\mathbf{u}}_i \geq 0$, the derivative computed at \eqref{eq: w} vanishes for:
\begin{equation}
   \hat{\mathbf{u}}_i = \frac{1}{1 + \alpha \tau \mathbf{I_\Omega}(i,i)} \mathbf{w}_i,
\end{equation} 
and it holds for $\mathbf{w}_i\geq 0$. Similarly, for the case $\hat{\mathbf{u}}_i < 0$, this analysis yields:
\begin{equation}
    \hat{\mathbf{u}}_i = \frac{1}{1 + \alpha \tau ( \mathbf{I_\Omega}(i,i)+1)} \mathbf{w}_i,
\end{equation}
for $\mathbf{w}_i < 0$.

So finally, the proximal operator is given by:
\begin{equation} 
 \left( \textbf{\text{prox}}_{h, \tau}(\mathbf{w}) \right)_i = {\text{prox}}_{h, \tau}(\mathbf{w}_i) = 
    \begin{cases}
      \frac{\mathbf{w}_i}{1 + \alpha \tau \mathbf{I_\Omega}(i,i)} & \text{if } {\mathbf{w}_i} \geq 0,\\
      \frac{\mathbf{w}_i}{1 + \alpha \tau ( \mathbf{I_\Omega}(i,i)+1)} & \text{if } \mathbf{w}_i < 0.
    \end{cases}       
\end{equation}

In a similar way, we compute the proximal mapping of the function $\overline{h}: \mathbb{R}^{L^2} \rightarrow \mathbb{R}$, defined in (\ref{eq:overline_h}), as follows:
\begin{align}
    \textbf{\text{prox}}_{\overline{h},\tau}(\mathbf{z}) 
    &= \argmin_{\mathbf{u}}\left(\frac{1}{2\tau}\|\mathbf{u - z}\|_2^2 + \overline{h}(\mathbf{u})\right) \nonumber \\
    &= \argmin_{\mathbf{u}}\left(\frac{1}{2\tau}\|\mathbf{u - z}\|_2^2 + \frac{\alpha}{2} \left(\|\mathbf{I_\Omega u}\|_2^2 + \|\mathbf{I}_{\hat{\mathbf{x}}_\mu} \mathbf{u}\|_2^2\right)\right).
    \label{eq: prox_overline_h}
\end{align}
The optimal solution $\hat{\mathbf{u}}$ of \eqref{eq: prox_overline_h} ($\hat{\mathbf{u}} = \textbf{\text{prox}}_{\overline{h},\tau}(\mathbf{z})$) is attained when: 
\begin{align}\label{eq: prox_comp}
    \mathbf{0} &\in \nabla\left(\frac{1}{2\tau}\|\hat{\mathbf{u}} - \mathbf{z}\|_2^2+ \frac{\alpha}{2} \left(\|\mathbf{I_\Omega} \hat{\mathbf{u}}\|_2^2+ \|\mathbf{I}_{\hat{\mathbf{x}}_\mu} \hat{\mathbf{u}}\|_2^2\right)\right), \nonumber \\
    %\mathbf{0} &\in \frac{1}{t}(\hat{\mathbf{u}} - \mathbf{x})+\frac{\alpha}{2}( 2\mathbf{I_\Omega} ^\intercal \mathbf{I_\Omega} \hat{\mathbf{u}}+ 2{\mathbf{I}_{\mathbf{x}_\mu}}^\intercal \mathbf{I}_{\mathbf{x}_\mu} \hat{\mathbf{u}}) \nonumber \\
    \mathbf{0} &\in \frac{1}{\tau}\left(\hat{\mathbf{u}} - \mathbf{z}\right)+\alpha\left(\mathbf{I_\Omega} \hat{\mathbf{u}}+ \mathbf{I}_{\hat{\mathbf{x}}_\mu} \hat{\mathbf{u}}\right). 
    %\mathbf{x}_i &= (1 + \alpha t (\mathbf{I_\Omega}(i,i) + \mathbf{I}_{\mathbf{x}_\mu}(i,i) ))\hat{\mathbf{u}}_i \nonumber \\ 
\end{align}
By eliminating $\hat{\mathbf{u}}$ in the expression \eqref{eq: prox_comp}, we compute element-wise the proximal operator:
\begin{equation}
     (\textbf{\text{prox}}_{\overline{h}, \tau}(\mathbf{z}))_i= {\text{prox}}_{\overline{h}, \tau}(\mathbf{z}_i) = \frac{\mathbf{z}_i}{1 + \alpha \tau \left( \mathbf{I_\Omega}(i,i) + \mathbf{I}_{\hat{\mathbf{x}}_\mu}(i,i)\right)}.
\end{equation}

\section{Appendix. The minimization problem to estimate $\hat{\mathbf{x}}'_\mu$}\label{appendixB}

Starting from the penalized optimization problem (\ref{eq:intensity_penalized}) and having $\mathbf{b}$ fixed, we aim to find a relation that contains the optimal $\hat{\mathbf{x}}_\mu$. While there are only quadratic terms, we proceed as following : 
\begin{align} 
    \mathbf{0} & \in \nabla\left(~\frac12 \| \mathbf{\Psi} \hat{\mathbf{x}}_\mu - (\overline{\mathbf{y}} - \mathbf{b})\|_2^2 + \frac{\mu}{2} \|\nabla\hat{\mathbf{x}}_\mu\|_2^2 + \frac{\alpha}{2} \left(\|\mathbf{I_\Omega} \hat{\mathbf{x}}_\mu\|_2^2 + \sum_{i=1}^{L^2}\ [\phi((\hat{\mathbf{x}}_\mu)_i)]^2\right)\right),\nonumber\\
    \mathbf{0} & \in \mathbf\Psi^\intercal \left( \mathbf\Psi \hat{\mathbf{x}}_\mu - (\overline{\mathbf{y}} - \mathbf{b})\right) + \mu \nabla^\intercal \nabla \hat{\mathbf{x}}_\mu + \alpha \left( \mathbf{I_\Omega} \hat{\mathbf{x}}_\mu + [\phi((\hat{\mathbf{x}}_\mu)_i)\phi'((\hat{\mathbf{x}}_\mu)_i)]_{\{i=1,...,{L^2}\}}\right).
\end{align}

Given \eqref{eq:phi'} we can write:
\begin{equation}
        \mathbf{0} \in \mathbf\Psi^\intercal \left( \mathbf\Psi \hat{\mathbf{x}}_\mu - \overline{\mathbf{y}} - \mathbf{b}\right) + \mu \nabla^\intercal \nabla \hat{\mathbf{x}}_\mu + \alpha \left( \mathbf{I_\Omega} \hat{\mathbf{x}}_\mu + [\phi((\hat{\mathbf{x}}_\mu)_i)]_{\{i=1,...,{L^2}\}}\right).
\end{equation}

Our goal is to compute $\hat{\mathbf{x}}'_\mu$, the partial derivative of $\hat{\mathbf{x}}_\mu$ w.r.t. $\mu$. So, we derive as follows: 
\begin{align}
        \frac{\partial}{\partial\mu} & \left(\mathbf{0} \in \mathbf\Psi^\intercal \left( \mathbf\Psi \hat{\mathbf{x}}_\mu - \overline{\mathbf{y}} - \mathbf{b}\right) + \mu \nabla^\intercal \nabla \hat{\mathbf{x}}_\mu + \alpha \left( \mathbf{I_\Omega} \hat{\mathbf{x}}_\mu + [\phi((\hat{\mathbf{x}}_\mu)_i)]_{\{i=1,...,{L^2}\}}\right)\right), \nonumber\\
        \mathbf{0} &\in \mathbf\Psi^\intercal \mathbf\Psi \hat{\mathbf{x}}'_\mu + \mu \nabla^\intercal \nabla \hat{\mathbf{x}}'_\mu + \nabla^\intercal \nabla \hat{\mathbf{x}}_\mu + \alpha \left( \mathbf{I_\Omega} \hat{\mathbf{x}}'_\mu + [\phi'((\hat{\mathbf{x}}_\mu)_i) (\hat{\mathbf{x}}' _\mu)_i]_{\{i=1,...,{L^2}\}}\right). \label{eq:aux}
\end{align}

We define the matrix $\mathbf{I}_{\hat{\mathbf{x}}_\mu}$ such as:
\[
\mathbf{I}_{\hat{\mathbf{x}}_\mu}(i,i) = \begin{cases}
 0 & \text{if ${(\hat{\mathbf{x}}_\mu)_i} \geq 0$},\\
          1 & \text{if ${(\hat{\mathbf{x}}_\mu)_i} < 0$}  .
\end{cases}
\]

Now the vector $[\phi'((\hat{\mathbf{x}}_\mu)_i) (\hat{\mathbf{x}}' _\mu)_i]_{\{i=1,...,{L^2}\}}$, using further the equation (\ref{eq:phi'}), can be simply written as: $\mathbf{I}_{\hat{\mathbf{x}}_\mu} \hat{\mathbf{x}}'_\mu$ and then \eqref{eq:aux} becomes:
\[
\mathbf{0} \in \mathbf\Psi^\intercal \mathbf\Psi \hat{\mathbf{x}}'_\mu + \mu \nabla^\intercal \nabla \hat{\mathbf{x}}'_\mu + \nabla^\intercal \nabla \hat{\mathbf{x}}_\mu + \alpha \left( \mathbf{I_\Omega} \hat{\mathbf{x}}'_\mu + \mathbf{I}_{\hat{\mathbf{x}}_\mu}\hat{\mathbf{x}}'_\mu\right).
\]

The minimization problem we should solve in order to find $\hat{\mathbf{x}}'_\mu$ thus is:
\begin{equation}
     \hat{\mathbf{x}}'_\mu = \argmin_{\mathbf{x}\in \mathbb{R}^{L^2}} \frac12\|\mathbf{\Psi x} \|_2^2 + \frac{\mu}{2} \|\nabla\mathbf{x} + \frac{1}{\mu} \nabla\hat{\mathbf{x}}_\mu\|_2^2+ \frac{\alpha}{2} \left(\|\mathbf{I_\Omega x}\|_2^2 + \| \mathbf{I}_{\hat{\mathbf{x}}_\mu}\mathbf{x}\|^2_2\right).
\end{equation}

\section{Appendix. Algorithmic restart.}\label{appendixc}
\begin{figure}
    \centering
    \includegraphics[width = 0.3\textwidth]{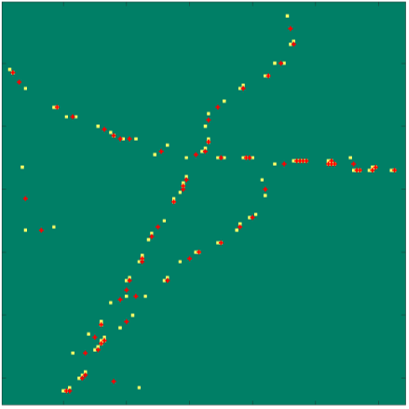}
    \caption{The yellow pixels belong to the support estimated in the previous restarting, while the red pixels belong to the initialization that is used in the current restarting}
    \label{fig:restarting}
\end{figure}

Every  initialization is based on the solution obtained at the previous restarting. There are many ways to choose the new initialization, deterministic and stochastic ones. In this paper we chose a deterministic way based on the following idea: for every pixel belonging to the solution of the previous restarting we find its closest neighbor. Then, we define the middle point between the two and we include it in the initialization of the current restarting. A small example is given in the Figure \ref{fig:restarting}. The yellow points belong to the support estimation of the previous restarting. Starting from them we define the red points, used for the initialization of the current restarting. 

\end{appendix}

\begin{Backmatter}

\paragraph{Acknowledgements}
The authors would like to thank E. van Obberghen-Shilling and D.Grall from the Institut de Biologie Valrose (iBV) who kindly prepared and provided the experimental samples. Furthermore, we would like to thank the anonymous reviewers for their valuable comments and suggestions.

\paragraph{Funding Statement}
The work of VS and LBF has been supported by the French government, through the 3IA Côte d’Azur Investments in the Future project managed by the National Research Agency (ANR) with the reference number ANR-19-P3IA-0002. LC acknowledges the support received by the academy on complex systems of UCA JEDI, the one received by the EU H2020 RISE program NoMADS, GA 777826, and the one received by the GdR ISIS grant SPLIN. The work of JHG was supported by the French Agence Nationale de la Recherche in the context of the project Investissements d’Avenir UCAJEDI (ANR-15-IDEX-01). Support for development of the microscope was received from IBiSA (Infrastructures en Biologie Santé et Agronomie) to the MICA microscopy platform.

\paragraph{Competing Interests}
None

\paragraph{Data Availability Statement}
Replication data and code can be found in: \href{https://github.com/VStergiop/COL0RME}{https://github.com/VStergiop/COL0RME}.
%A statement about how to access data, code and other materials allowing users to understand, verify and replicate findings --- e.g. Replication data and code can be found in Harvard Dataverse: \verb+\url{https://doi.org/link}+.

\paragraph{Ethical Standards}
The research meets all ethical guidelines, including adherence to the legal requirements of the study country.

\paragraph{Author Contributions}
% Methodology: V.S; L.C; J.H.G; L.B.F. Software: V.S, J.H.G. Data curation: S.S. Data visualisation: S.S. Supervision: L.C; L.B.F. Writing original draft: V.S. All authors approved the final submitted draft.
VS, LC, JHG and LBF conceived and designed the study. SS conducted data gathering. VS and JHG implemented the software. VS carried out the experiments. LC and LBF supervised the work. VS, LC, JHG and LBF wrote the article. All authors approved the final submission.

\paragraph{Supplementary Material}
A supplementary document intended for publication has been provided with the submission.

% Bibliography
\bibliographystyle{vancouver}
\bibliography{biblio}

\end{Backmatter}

\end{document}